\documentclass[english,12pt]{smfart}
\usepackage{amssymb}
\usepackage{amsfonts}
\usepackage{amscd}
\usepackage[all]{xypic}
\usepackage{amsmath}

\usepackage{mathrsfs}

\CompileMatrices

\swapnumbers

\def\sq{{\square}}
\def\ord{{\rm ord}}
\def\ac{{\overline{\rm ac}}}

\def\Supp{{\rm Supp}}

\def\Lp{{{\mathbf L}_{\rm PR}}}
\def\Eu{{\rm Eu}}
\def\Var{{\rm Var}}
\def\Def{{\rm Def}}
\def\Cons{{\rm Cons}}
\def\RDef{{\rm RDef}}
\def\GDef{{\rm GDef}}

\def\SA{{\rm SA}}
\def\LPas{\cL_{\rm DP}}
\def\LPre{\cL_{\rm DP,P}}

\def\Kdim{{\rm Kdim}\, }

\def\ordjac{{\rm ordjac}}

\def\11{{\mathbf 1}}
\def\AA{{\mathbb A}}

\def\CC{{\mathbb C}}

\def\FF{{\mathbb F}}
\def\GG{{\mathbb G}}

\def\LL{{\mathbb L}}

\def\NN{{\mathbb N}}

\def\QQ{{\mathbb Q}}
\def\RR{{\mathbb R}}

\def\ZZ{{\mathbb Z}}

\def\bfL{{\mathbf L}}

\def\cA{{\mathcal A}}

\def\cC{{\mathscr C}}

\def\cL{{\mathcal L}}
\def\cM{{\mathcal M}}

\def\cO{{\mathcal O}}
\def\cP{{\mathcal P}}

\def\cU{{\mathcal U}}

\def\cX{{\mathcal X}}
\def\cY{{\mathcal Y}}
\def\cZ{{\mathcal Z}}
\def\llp{\mathopen{(\!(}}
\def\llb{\mathopen{[\![}}
\def\rrp{\mathopen{)\!)}}
\def\rrb{\mathopen{]\!]}}

\def\dl{{\rm deg}_{\LL}}
\mathchardef\alphag="7C0B
\mathchardef\betag="7C0C
\mathchardef\gammag="7C0D
\mathchardef\deltag="7C0E
\mathchardef\varepsilong="7C22
\mathchardef\varphig="7C27
\mathchardef\psig="7C20
\mathchardef\zetag="7C10
\mathchardef\epsilong="7C0F
\mathchardef\rhog="7C1A
\mathchardef\taug="7C1C
\mathchardef\upsilong="7C1D
\mathchardef\iotag="7C13
\mathchardef\thetag="7C12
\mathchardef\pig="7C19
\mathchardef\sigmag="7C1B
\mathchardef\etag="7C11
\mathchardef\omegag="7C21
\mathchardef\kappag="7C14
\mathchardef\lambdag="7C15
\mathchardef\mug="7C16
\mathchardef\xig="7C18
\mathchardef\chig="7C1F
\mathchardef\nug="7C17
\mathchardef\varthetag="7C23
\mathchardef\varpig="7C24
\mathchardef\varrhog="7C25
\mathchardef\varsigmag="7C26
\mathchardef\Omegag="7C0A
\mathchardef\Thetag="7C02
\mathchardef\Sigmag="7C06
\mathchardef\Deltag="7C01
\mathchardef\Phig="7C08
\mathchardef\Gammag="7C00
\mathchardef\Psig="7C09
\mathchardef\Lambdag="7C03
\mathchardef\Xig="7C04
\mathchardef\Pig="7C05
\mathchardef\Upsilong="7C07

\newtheorem{theorem}[subsubsection]{Theorem}
\newtheorem{lem}[subsubsection]{Lemma}
\newtheorem{cor}[subsubsection]{Corollary}
\newtheorem{prop}[subsubsection]{Proposition}

\theoremstyle{definition}
\newtheorem{definition}[subsubsection]{Definition}
\newtheorem{example}[subsubsection]{Example}

\newtheorem{def-prop}[subsubsection]{Proposition-Definition}
\newtheorem{def-theorem}[subsubsection]{Theorem-Definition}
\newtheorem{def-lem}[subsubsection]{Lemma-Definition}

\theoremstyle{remark}
\newtheorem{remark}[subsubsection]{Remark}

\newtheorem{claim}[subsubsection]{Claim}

\theoremstyle{plain}

\numberwithin{equation}{subsection}

\def\boxit#1#2{\setbox1=\hbox{\kern#1{#2}\kern#1}%
\dimen1=\ht1 \advance\dimen1 by #1
\dimen2=\dp1 \advance\dimen2 by #1
\setbox1=\hbox{\vrule height\dimen1 depth\dimen2\box1\vrule}%
\setbox1=\vbox{\hrule\box1\hrule}%
\advance\dimen1 by .4pt \ht1=\dimen1
\advance\dimen2 by .4pt \dp1=\dimen2 \box1\relax}

\newcommand{\sur}[2]{\genfrac{}{}{0pt}{}{#1}{#2}}
\renewcommand{\theequation}{\thesubsection.\arabic{equation}}

\mathchardef\alphag="7C0B
\mathchardef\betag="7C0C
\mathchardef\gammag="7C0D
\mathchardef\deltag="7C0E
\mathchardef\varepsilong="7C22
\mathchardef\varphig="7C27
\mathchardef\psig="7C20
\mathchardef\zetag="7C10
\mathchardef\epsilong="7C0F
\mathchardef\rhog="7C1A
\mathchardef\taug="7C1C
\mathchardef\upsilong="7C1D
\mathchardef\iotag="7C13
\mathchardef\thetag="7C12
\mathchardef\pig="7C19
\mathchardef\sigmag="7C1B
\mathchardef\etag="7C11
\mathchardef\omegag="7C21
\mathchardef\kappag="7C14
\mathchardef\lambdag="7C15
\mathchardef\mug="7C16
\mathchardef\xig="7C18
\mathchardef\chig="7C1F
\mathchardef\nug="7C17
\mathchardef\varthetag="7C23
\mathchardef\varpig="7C24
\mathchardef\varrhog="7C25
\mathchardef\varsigmag="7C26
\mathchardef\Omegag="7C0A
\mathchardef\Thetag="7C02
\mathchardef\Sigmag="7C06
\mathchardef\Deltag="7C01
\mathchardef\Phig="7C08
\mathchardef\Gammag="7C00
\mathchardef\Psig="7C09
\mathchardef\Lambdag="7C03
\mathchardef\Xig="7C04
\mathchardef\Pig="7C05
\mathchardef\Upsilong="7C07
\newcommand{\Pm}{\mathcal{P}}

\DeclareMathOperator*{\Spec}{Spec}

\def\ord{{\rm ord}}

\begin{document}

\title[Constructible motivic functions and motivic integration]{Constructible motivic functions and motivic integration}

\author{Raf Cluckers}
\address{Katholieke Universiteit Leuven, Departement wiskunde,
Celestijnenlaan 200B, B-3001 Leu\-ven, Bel\-gium. Current address:
\'Ecole Normale Sup\'erieure, D\'epartement de
ma\-th\'e\-ma\-ti\-ques et applications, 45 rue d'Ulm, 75230 Paris
Cedex 05, France} \email{cluckers@ens.fr}
\urladdr{www.dma.ens.fr/$\sim$cluckers/}

\author{Fran\c cois Loeser}

\address{{\'E}cole Normale Sup{\'e}rieure,
D{\'e}partement de math{\'e}matiques et applications,
45 rue d'Ulm,
75230 Paris Cedex 05, France
(UMR 8553 du CNRS)}
\email{Francois.Loeser@ens.fr}
\urladdr{www.dma.ens.fr/$\sim$loeser/}
%\dedicatory{Preliminary notes (\today)}

%\begin{abstract}
%\end{abstract}

\maketitle

\renewcommand{\partname}{}
%\part{e2e2e2}

\section{Introduction}

\subsection{}In this paper, intended to be the first in a series,
we lay new general foundations for motivic integration and give
answers to some important issues in the subject. Since its
creation by Maxim Kontsevich \cite{k}, motivic integration
developed
quickly
and has spread out in many directions. In a nutshell, in motivic
integration, numbers are replaced by geometric objects, like
virtual varieties, or motives. But, classicaly, not only numbers
are defined using integrals, but also interesting classes of
functions. The previous constructions of motivic integration were
all quite geometric, and it was quite unclear how they could be
generalized to handle integrals depending on parameters. The new
approach we present here, based on cell decomposition, allows us
to develop a very general  theory of motivic integration taking
parameters in account. More precisely, we define a natural class
of functions - constructible motivic functions - which is stable
under integration.

The basic idea underlying our approach is to construct more
generally push-forward morphisms $f_!$ which are functorial - they
satisfy $(f \circ g)_! = f_! \circ g_!$ - so that performing
motivic integration corresponds to taking the push-forward to the
point. This strategy has many technical advantages.
In essence,
it allows to reduce the construction of $f_!$ to the case of
closed immersions and projections, and in the latter case we can
perform induction on the relative dimension, the basic case being
that of relative dimension 1, for which we can make use of the
Cell Decomposition Theorem of Denef-Pas \cite{Pas}.

\subsection{}Our main construction being inspired by analogy with  integration along the Euler
characteristic for constructible functions over the reals, let us
first present a brief overview of this theory, for which we refer
to \cite{mp}, \cite{schap}, and \cite{viro} for more details. We
shall put some emphasis on formulation in terms of Grothendieck
rings. Let us denote by $\SA_{\RR}$ the category of real
semialgebraic sets, that is, objects of $\SA_{\RR}$ are
semialgebraic sets and morphisms are semialgebraic maps. Since
every real semialgebraic set admits  a semialgebraic
triangulation, the Euler characteristic of real semialgebraic sets
may be defined as the unique $\ZZ$-valued additive invariant on
the category of real semialgebraic sets which takes value one on
closed simplexes. More precisely, let us define $K_0 (\SA_{\RR})$,
the Grothendieck ring of real semialgebraic sets,  as the quotient
of the free abelian group on symbols $[X]$, for $X$ real
semialgebraic, by the relations $[X] = [X']$ if $X$ and $X'$ are
isomorphic, and $[X \cup Y] = [X] + [Y] - [X \cap Y]$, the product
being induced by the cartesian product of semialgebraic sets.
Then, existence of semialgebraic triangulations easily implies the
following statement:

\begin{prop}The Euler characteristic morphism
$[X] \mapsto \Eu (X)$ induces a ring isomorphism
$$
K_0 (\SA_{\RR}) \simeq \ZZ.
$$
\end{prop}

A constructible function on a semialgebraic set $X$ is a function
$\varphi : X \rightarrow \ZZ$ that can be written as a finite sum
$\varphi = \sum_{i \in I} m_i \11_{X_i}$ with $m_i$ in $\ZZ$,

 $X_i$
semialgebraic subsets of $X$, and $\11_{X_i}$ the characteristic
function of $X_i$.
 The set $\Cons (X)$ of constructible functions on $X$ is a ring.
If $f : X \rightarrow Y$ is a morphism of semialgebraic sets, we
have a natural pullback morphism $f^* : \Cons (Y ) \rightarrow \Cons
(X)$ given by $\varphi \mapsto \varphi \circ f$. Now let us explain
how the construction of a push-forward morphism $f_* : \Cons (X)
\rightarrow \Cons (Y)$ is related to integration with respect to
Euler characteristic.

Let $\varphi = \sum_{i \in I} m_i \11_{X_i}$ be in $\Cons (X)$.
One sets
$$
\int_X \varphi := \sum_{i \in I} m_i \Eu (X_i).
$$
It is quite easy to check that this quantity depends only on
$\varphi$. Now if $f : X \rightarrow Y$ is a morphism, one checks
that defining $f_*$ by
$$
f_* (\varphi) (y) = \int_{f^{-1} (y)} \varphi_{|f^{-1} (y)},
$$
indeed yields a morphism $f_*   : \Cons (X) \rightarrow \Cons (Y),$
and that furthermore $(f \circ g)_* = f_* \circ g_*$. For our
purposes it will be more enlightening to express the preceding
construction in terms of relative Grothendieck rings.

For $X$ a semialgebraic set, let us consider the category $\SA_X$
of semialgebraic sets over $X$. Hence objects of $\SA_X$ are
morphisms $Y \rightarrow X$ in $\SA_{\RR}$ and a morphism $(h : S
\rightarrow X) \rightarrow (h': S' \rightarrow X)$ in $\SA_X$ is a
morphism $g : S \rightarrow S'$ such that $h'\circ g = h$. Out of
$\SA_X$, one constructs a Grothendieck ring $K_0 (\SA_X)$
similarly as before, and we have the following statement, which
should be folklore, though we could not find it the literature.

\begin{prop}Let $X$ be a semialgebraic set.
\begin{enumerate}
\item[(1)]The mapping $[h : S \rightarrow X] \mapsto h_* (\11_S)$
induces an isomorphism
\begin{equation*}
K_0 (\SA_X) \simeq \Cons (X).
\end{equation*}
\item[(2)]Let $f : X \rightarrow Y$ be a morphism in $\SA_{\RR}$.
Under the above isomorphism $f_* : \Cons (X) \rightarrow \Cons
(Y)$ corresponds to the morphism $K_0 (\SA_X) \rightarrow K_0
(\SA_Y)$ induced by composition with $f$.
\end{enumerate}
\end{prop}

\subsection{}Let us now explain more about
our framework. Fix a field $k$ of characteristic $0$. We want to
integrate (functions defined on) subobjects of $k \llp t \rrp^m$. For
technical reasons it is wiser to consider more generally
integration on subobjects of $k \llp t \rrp^m \times k^n \times \ZZ^r$.
This will allow considering parameters lying in the valued field,
the residue field, and the value group. In fact, we shall restrict
ourselves to considering a certain class of reasonably tame
objects, that of definable subsets in a language $\LPas$.
Typically these objects are defined by formulas involving usual
symbols $0, 1, +, -, \times$ for the $k \llp t \rrp$ and $k$ variables,
and $0, 1, +, -, \leq$ for the $\ZZ$-variables, and also symbols
$\ord$ for the valuation and $\ac$ for the first non trivial
coefficient of elements of $k\llp t \rrp$, and the usual logical symbols
(see \ref{sec:pas} for more details). Furthermore we shall not
only consider the set of points in $k \llp t \rrp^m \times k^n \times
\ZZ^r$ satisfying a given formula $\varphi$, but also look to the
whole family of subsets of $K \llp t \rrp^m \times K^n \times \ZZ^r$,
for $K$ running over all fields containing $k$, of points that
satisfy $\varphi$. This is what we call definable subassignments.
Definable subassignments form a category and are our basic objects
of study.

Let us fix such a definable subassignment $S$. Basically,
constructible motivic functions on $S$ are built from

- classes $[Z]$ in a suitable Grothendieck ring of definable
subassignments $Z$ of $S \times \AA_k^d$ for some $d$ ;

- symbols $\LL^{\alpha}$, where $\LL$ stands for the class of the
relative affine line over $S$ and $\alpha$ is some definable
$\ZZ$-valued function on $S$;

- symbols $\alpha$ for $\alpha$ a definable $\ZZ$-valued function
on $S$.

Constructible motivic functions on $S$ form a ring $\cC (S)$. Any
definable subassignment $C$ of $S$ has a characteristic function
$\11_C$ in $\cC (S)$.

\subsection{}\label{ex}We explain now on an example how one can recover the motivic
volume by considering the push-forward of constructible functions.
We shall consider the points of the affine elliptic curve $x^2 = y
(y - 1) (y- 2)$ with nonnegative valuation, namely the definable
subassignment $C$ of $\AA^2_{k \llp t \rrp}$  defined by the conditions
$$x^2 = y (y - 1) (y- 2),  \quad \ord (x) \geq 0 \quad \text{and} \quad
\ord (y) \geq 0.$$ Since the affine elliptic curve $E$ defined by
$\xi^2 = \eta (\eta - 1) (\eta - 2)$ in $\AA^2_k$ is smooth, we
know that the motivic volume $\mu(C)$ should be equal to
$\frac{[E]}{\LL}$, cf.~\cite{arcs}. Let us consider the projection
$p : \AA^2_{k \llp t \rrp} \rightarrow \AA^1_{k \llp t \rrp}$ given by $(x, y)
\mapsto y$. In our formalism $p_! ([\11_C])$ is equal to a sum $A
+  B_0 + B_1 + B_2$ with
\begin{gather*}
A = [\xi^2 = \ac (y)(\ac (y) - 1) (\ac (y) - 2)]
[\11_{C(A)}]\\
B_0 =    [\xi^2 = 2 \ac (y)]
[\11_{C_0}]  \LL^{\ord (y) /2}\\
B_1 = [\xi^2 = -  \ac (y - 1)]
[\11_{C_1}]  \LL^{\ord (y - 1) /2}\\
\intertext{and} B_2 = [\xi^2 = 2 \ac (y - 2)] [\11_{C_2}]
\LL^{\ord (y - 2) /2},
\end{gather*}
with $C (A) = \{y \,  \vert\,  \ord (y) = \ord (y - 1) = \ord (y -
2) = 0\}$ and $C_i = \{y \, \vert \, \ord (y - i) >0 \, {\rm
and}\, \ord (y- i) \equiv 0 \mod 2\}$. So $p_! ([\11_C])$  looks
already  like a quite general constructible motivic function.

Let us show how one can recover the motivic volume of $\mu (C)$ by
computing the integral of $p_! ([\11_C])$ on $\AA^1_{k \llp t
\rrp}$. Let $\pi_i$ denote the projection of $\AA^i_{k \llp t \rrp}$
on the point. One computes $\pi_{1!} (A) = \frac{[E] - 3}{\LL}$,
while summing up the corresponding geometric series leads to that $
\pi_{1!} (B_0) = \pi_{1!} (B_1) = \pi_{1!} (B_2) = \LL^{-1}$, so
that finally $\pi_{1!} (p_! ([1_C])) = \frac{[E]}{\LL}$. Hence the
computation
 fits with the requirements $\pi_{1!} \circ p_! =
(\pi_1 \circ p)_! = \pi_{2!}$ and $\pi_{2!} ([1_C]) = \mu (C)$.

As we will see in the main construction of the push forward operator
denoted with subscript $_!$, in this example $p_!$ is calculated
with ``the line element'' determined by the forms ``$dx$'' and
``$dy$'' on $C$ and $\pi_{1!}$ calculates an integral over the line
with respect to the form ``$dy$'' (see Theorem \ref{mt}). In our
context the line element is of course non-archimedean, see section
\ref{canvf}.

\subsection{}Such a computation is maybe a bit surprising at first sight, since one could
think that is not possible to recover the motive of an elliptic
curve by projecting on to the line and computing the volume of the
fibers, which consist of  0, 1 or 2 points. The point is that our
approach is not so naive and keeps track of the elliptic curve
which remains encoded at the residue field level. Our main construction
can be considered as a vast amplification of that example and one
may understand that the main difficulty in the construction is
proving that our construction of $f_!$ is independent of the way
we may decompose $f$ into  a composition of morphisms.

In fact, we do not integrate functions in $\cC (S)$, but rather
their classes in a graded object $C (S) = \oplus_d C^d (S)$. The
reason for that is that we have to take in account dimension
considerations. For instance we could factor the identity morphism
from the point to itself as the composition of an embedding in the
line with the projection of the line on the point. But then a
problem arises: certainly the point should be of measure 1 in
itself, but as a subset of  the line it should be of measure 0! To
circumvent this difficulty, we filter $\cC (S)$ by ``$k
\llp t \rrp$-dimension of support''. Typically, if $\varphi$ has ``$k
\llp t \rrp$-dimension of support'' equal to $d$, we denote by
$[\varphi]$ its class in the graded piece $C^d (S)$\footnote{That
notation was already used in \ref{ex} without explanation.}. We
call elements of $C (S)$ constructible motivic Functions (with
capital F).
 One further difficulty is that arbitrary elements of
$C(S)$ may not be integrable, that is, the corresponding integral
could diverge. So we need to define at the same time the integral
(or the push-forward) and the integrability condition. Also, as in
the usual construction of Lebesgue integral, it is technically
very useful to consider first only ``positive constructible
functions'' on $S$.  They form a semiring $\cC_+ (S)$ and we may
consider the corresponding  graded object $C_+ (S)$. An important
difference with the classical case, is that in general the
canonical morphisms $\cC_+ (S) \rightarrow \cC (S)$ and $C_+ (S)
\rightarrow C (S)$ are not injective.

The main achievement of the present paper is the following:  we
establish existence and uniqueness of
 a) a
subgroup ${\rm I}_{S'}C_+ (S)$ of $C_+(S)$ consisting of
$S'$-integrable  positive Functions on $S$, b) a  push-forward morphism $f_!
: {\rm I}_{S'} C_+ (S) \rightarrow C_+(S')$, under a certain system of
natural axioms, for every morphism $f : S \rightarrow S'$  of definable
subassignments.

\subsection{}Once  the main result is proven, we
can grasp  its rewards. Firstly, it may be directly generalized to
the relative setting of integrals with parameters. In particular we
get that motivic integrals parametrized by a definable subassignment
$S$ take their  values in $\cC_+ (S)$ or in $\cC (S)$. Also, our use
of the quite abstract definable subassignments allows us to work at
a level of generality that encompasses both
 ``classical''  motivic integration as developed in \cite{arcs}
and  the ``arithmetical'' motivic integration of \cite{JAMS}. More
precisely, we show that the present theory may be specialized both
to ``classical'' motivic integration and ``arithmetical'' motivic
integration, but with the bonus that no more completion process is
needed. Indeed, there is a canonical forgetful morphism $\cC ({\rm
point}) \rightarrow K_0 ({\rm Var}_k)_{\rm rat}$, with $K_0 ({\rm
Var}_k)_{\rm rat}$ the localization of the Grothendieck ring of
varieties over $k$ with respect to $\LL$ and $1 - \LL^{-n}$, $n \geq
1$, that sends the motivic volume of a definable object as defined
here, to a representative of the ``classical'' motivic volume in
$K_0 ({\rm Var}_k)_{\rm rat}$. So in the definable setting,
``classical'' motivic volume takes values in  $K_0 ({\rm
Var}_k)_{\rm rat}$ (and not in any completion of it). Such a
result lies in the fact that in our machinery, the only infinite
process that occurs is summation of geometric series in powers of
$\LL^{-1}$. A similar statement holds in the arithmetic case.

Another important feature is that no use at all is made of
desingularization results. On the other side we rely very strongly
on the Cell Decomposition Theorem  of Denef-Pas. This makes in some sense things much
worse in positive characteristic, since then desingularization  is
a reasonable conjecture while there is even no sensible guess of
what cell decomposition could be in that case!

\subsection{}Let us now describe  briefly the content
of the paper. Our basic objects of study are the various
categories of definable subassignments in the Denef-Pas language
that we review in section \ref{sec1}. An essential feature of
these definable subassignments is that they admit a good dimension
theory with respect to the valued field variables that we call
$K$-dimension. This is established in section \ref{sec:dim}. As a
first step in constructing motivic integrals, we develop in
section \ref{secpre} a general machinery for summing over the
integers. This is done in the framework of functions definable in
the Presburger language. We prove a general rationality statement
theorem \ref{presrat} which we formulate in terms of  a Mellin
transformation. This allows to express punctual summability of a
series in terms of polar loci of its Mellin transform and thus to
define the sum of the  series by evaluation of the Mellin
transform at $1$. This construction is the main device that allows
us to avoid any completion process in our integration theory, in
contrast with previous approaches. In the following section
\ref{sec2}, we define constructible motivic functions and we
extend  the constructions and the results of the previous section
to this framework. After the short section \ref{sec4} which is
devoted to the construction of motivic constructible Functions (as
opposed to functions) and their relative variants, section
\ref{sec5} is devoted to cell decomposition, which is, as we
already stressed, a basic tool in our approach. We need a
definition of cells slightly more flexible than the one of
Denef-Pas for which we give the appropriate Cell Decomposition
Theorem \`a la Denef-Pas, and we also introduce bicells. We prove
a fundamental structure result, Theorem \ref{normal}, for
definable functions with values in the valued field which may show
interesting for its own right. Section \ref{df1} is devoted to
introducing basic notions of differential calculus, like
differential forms, volume forms and order of jacobian in the
definable setting. In section \ref{dim1}, which appears to be
technically quite involved, we construct motivic integrals in
relative dimension 1 (with respect to the valued field variable).
In particular we prove a fundamental change of variable formula in
relative dimension 1, whose proof uses Theorem \ref{normal}, and
which will be of essential use in the rest of the paper.

We are then able to state our main result, Theorem \ref{mt}, in
section \ref{sec7}, and section \ref{sec8} is devoted to its
proof. The idea of the proof is quite simple. By a graph
construction one reduces to the case of definable injections and
projections. Injections being quite easy to handle, let us
consider projections. We already now how to integrate with respect
to  $\ZZ$-variables and also with respect to one valued field
variable, integration with respect to  residue field variables
being essentially tautological. So to be able to deal with the
general case, we need to prove various statements of Fubini type,
that will allow us to interchange the order in which we perform
integration with respect to various variables. The most difficult
case is that of two valued field variables, that requires careful
analysis of what happens on various types of bicells. Let us note
that van den Dries encounters a similar difficulty in his
construction of Euler characteristics in the o-minimal framework
\cite{vdd3}. Once the main theorem is proved, we can derive the
main properties and applications. In section \ref{sec9}, we prove
a general change of variable formula and also the fundamental fact
that a positive Function that is bounded above by an integrable
Function is also integrable. We then develop the integration
formalism for Functions in $C (X)$ - that is with no positivity
assumption - in section \ref{sec11}. In section  \ref{sec10} we
consider integrals with parameters and extend all previous resuts
to this framework. As a side result, we prove the very general
rationality theorem \ref{ratth}.

The last part of the paper is devoted to generalization to the
global setting and to comparison results. In section \ref{secg}, we
consider integration on definable subsassignments of varieties. This
is done by replacing functions by volume forms, as one can expect.
More precisely, if $f$ is a morphism between global definable
subassignments $S$ and $S'$, we construct a morphism $f_!^{\rm top}$
sending $f$-integrable volume forms on $S$ to volume forms on $S'$,
which corresponds  to integrating Functions in top dimension in the
affine case. This provides the right framework for a general Fubini
Theorem for fiber integrals (Theorem \ref{MTFI}). We then show in
section \ref{compa} how our construction relates with the previous
constructions of motivic integration. In the paper  \cite{miami} we
explain how it specializes to $p$-adic integration and we also give
some applications to Ax-Kochen-Er{\v s}ov Theorems for integrals
depending on parameters. The main results of this paper have been
anounced in the notes \cite{cr1} and \cite{cr2}. The present %October
%2007 
version of the paper does not differ from the original version
except for very minor changes. Since our paper  was originally put
on the arxiv as math.AG/0410203,  we have been able to extend our
work to the exponential setting and to prove a general ``Transfer
principle"  \`a la Ax-Kochen-Er{\v s}ov in this context
 \cite{expnote}, \cite{exp}. Also Hrushovski and Kazhdan \cite{hk} developed a general theory of
 integration for valued fields based on Robinson's quantifier elimination for algebraically closed
 valued fields.

 In writing the paper we  tried our best keeping
it accessible to a wide audience including  algebraic geometers and
model theorists. In particular, only basic familiarity with the
first chapters of textbooks like those of Hartshorne \cite{hart} and
Marker \cite{marker} is required. We also attempted  to stay within
the realm  of geometry as much as possible. For example, we use,
with the hope it would appeal to geometers, the terminology of
``definable subassignments", introduced in \cite{JAMS}, which is
certainly familiar to logicians under other guises. By the
foundational nature of the paper, many constructions and proofs are
somewhat lengthy and technical, so we made an effort to make the
main results and properties directly accessible and usable by the
reader without having to digest all details. In particular,
potential users might gain additional motivation by having a first
look to  sections \ref{sec7}, \ref{sec9}, \ref{sec10} and \ref{secg}
as early as possible. Also, one should stress that, for many
applications, it is enough to consider integration in maximal
dimension.

\subsection*{}%{Acknowledgments}
\hspace{0.5cm}
The present work would not exist without Jan Denef, whose  insight and work,
in particular concerning the ubiquity of cell decomposition,
did have a strong   influence  on our approach.
%We would like to acknowledge our debt to Jan Denef,
We also thank him warmly for  his crucial encouragements when we
started this project in February 2002. During the preparation of
this work, we benefited from the support of many colleagues and
friends. In particular,  we would like to adress special  thanks to
Antoine Chambert-Loir, Clifton Cunningham, Lou van den Dries, Tom
Hales and  Udi Hrushovski for the interest they have shown in our
work, and for comments and useful discussions that helped to improve
the paper. The first author has been supported as a postdoctoral
fellow by the Fund for Scientific Research - Flanders (Belgium)
(F.W.O.) during the preparation of this paper.

$$
\star \quad \star \quad \star
$$

\newpage

\tableofcontents

\part{Preliminary constructions}

\section{Definable subassignments}\label{sec1}

In this section, we extend the notion of definable subassignments,
introduced in \cite{JAMS}, to the context of $\LPas$-definable
sets, with $\LPas$ a language of Denef-Pas.

\subsection{Languages of Denef-Pas}\label{sec:pas}
Let $K$ be a valued field, with valuation map $\ord:K^\times\to
\Gamma$ for some additive ordered group $\Gamma$, $R$ its
valuation ring, $k$ the residue field. We denote by $x\mapsto \bar
x$ the projection $R \rightarrow k$ modulo the maximal ideal $M$
of $R$. An angular component map (modulo $M$) on $K$ is a
multiplicative map $\ac:K^\times\to k^\times$ extended by putting
$\ac(0)=0$ and satisfying $\ac(x)=\bar x$ for all $x$ with
$\ord(x)=0$.

If $K=k\llp t \rrp$ for some field $k$, there exists a natural
valuation map $K^{\times} \to\ZZ$ and a natural angular component
map sending $x=\sum_{i\geq \ell} a_it^i$ in $K^{\times}$ with
$a_i$ in $ k$ and $a_{\ell}\not=0$ to $\ell$ and $a_{\ell}$,
respectively.

Fix an arbitrary expansion $\bfL_{\rm Ord}$ of the language of
ordered groups $(+,-,0,\leq)$ and an arbitrary expansion $\bfL_{\rm
Res}$ of the language of rings $\bfL_{\rm Rings}=(+,-,\cdot,0,1)$.
A language $\LPas$ of Denef-Pas is a three-sorted language of the
form $$(\bfL_{\rm Val},\bfL_{\rm Res},\bfL_{\rm Ord},\ord,\ac),$$
with as sorts:
\begin{itemize}
 \item[(i)] a ${\rm Val}$-sort for the valued field-sort,
 \item[(ii)] a ${\rm Res}$-sort for the residue
 field-sort, and
 \item[(iii)] an $\rm{Ord}$-sort for the value
 group-sort,
\end{itemize}
where the language $\bfL_{\rm Val}$ for the ${\rm Val}$-sort is the
language of rings $\bfL_{\rm Rings}$, and the languages $\bfL_{\rm
Res}$ and $\bfL_{\rm Ord}$ are used for the ${\rm Res}$-sort and
the ${\rm Ord}$-sort, respectively. We only consider structures
for $\LPas$ consisting of tuples $(K,k,\Gamma)$ where $K$ is a
valued field with value group $\Gamma$, residue field $k$, a
valuation map $\ord$, and an angular component map $\ac$, together
with an interpretation of $\bfL_{\rm Res}$ and $\bfL_{\rm Ord}$ in
$k$ and $\Gamma$, respectively.

When $\bfL_{\rm Res}$ is $\bfL_{\rm Rings}$ and $\bfL_{\rm Ord}$ is
the Presburger language
$$
\bfL_{\rm PR} = \{+, -, 0, 1, \leq\} \cup \{\equiv_n\ \mid n\in
\NN,\ n
>1\},
$$
with $\equiv_n$ the equivalence relation modulo $n$ and $1$ a
constant symbol (with the natural interpretation if $\Gamma=\ZZ$),
we write $\LPre$ for $\LPas$.

As is standard for first order languages, $\LPas$-formulas are
(meaningfully) built up from the $\LPas$-symbols together with
variables, the logical connectives $\wedge$ (and), $\vee$ (or),
$\neg$ (not), the quantifiers $\exists$, $\forall$, the equality
symbol $=$, and parameters\footnote{For first order languages,
function symbols need to have a Cartesian product of sorts as
domain, while the symbol $\ord$ has the valued field-sort minus
the point zero as domain. Our use of the symbol $\ord$ with
argument $x$ in a $\LPas$-formula is in fact an abbreviation for a
function with domain the ${\rm Val}$-sort which extends the
valuation (the reader may choose the value of $0$), conjoined with
the condition $x\not=0$.}.
%Also in a standard way, any
%$\LPas$-formula yields a subset of a Cartesian product of $K,k$,
%and $\Gamma$ for any $\LPas$-structure $(K,k,\Gamma)$, where this
%Cartesian product is up to order determined by the free variables
%(these are the variables not bound by quantifiers) in the formula.

Let us now recall the statement of the Denef-Pas Theorem on
elimination of valued field quantifiers. Fix a language $\LPas$ of
Denef-Pas. Denote by $H_{\ac, 0}$ the $\LPas$-theory of the above
described structures whose valued field is Henselian and whose
residue field is of characteristic zero.

\begin{theorem}[Denef-Pas]\label{pqe}
The theory $H_{\ac, 0}$  admits elimination of quantifiers in the
valued field sort. More precisely, every $\LPas$-formula $\phi (x,
\xi, \alpha)$ (without parameters), with $x$ variables in the ${\rm
Val}$-sort, $\xi$ variables in the ${\rm Res}$-sort and $\alpha$
variables in the ${\rm Ord}$-sort, is $H_{\ac, 0}$-equivalent to a
finite disjunction of formulas of the form
$$
\psi (\ac\, f_1 (x), \dots, \ac\, f_k (x), \xi) \wedge \vartheta
(\ord\, f_1 (x), \dots, \ord\, f_k (x), \alpha),
$$
with $\psi$ a $\bfL_{\rm Res}$-formula, $\vartheta$ a $\bfL_{\rm
Ord}$-formula and $f_1$, \dots, $f_k$ polynomials in $\ZZ [X]$.
\end{theorem}

Theorem \ref{pqe} is not exactly expressed this way in \cite{Pas}.
The present statement can be found in \cite{vdd2} (3.5) and (3.7).
We will mostly use the following corollary, which is standard in
model theory.

\begin{cor}
Let $(K,k,\Gamma)$ be a model of $H_{\ac, 0}$, $S\subset K$ be a
subring, $T_S$ be the diagram of $S$ in the language $\LPas\cup
S$, that is, $T_S$ is the set of atomic $\LPas\cup S$-formulas and
negations of atomic formulas $\varphi$ such that $S \models
\varphi$, and $H_S$ be the union of $H_{\ac, 0}$ and $T_S$. Then
Theorem \ref{pqe} holds with $H_{\ac, 0}$ replaced by $H_S$,
$\LPas$ replaced by $\LPas\cup S$, and $\ZZ[X]$ replaced by
$S[X]$.
\end{cor}

\subsection{General subassignments}\label{sec:sub}

Let $F : \cC \rightarrow {\rm Sets}$ be a functor from a category
$\cC$ to the category of sets.
Any data, which associates to each object $C$ of $\cC$ a subset
$h (C)$ of $F (C)$, will be called a {\em subassignment} of $F$.
The point in this definition is that $h$ is not assumed to be a
subfunctor of $F$.

For $h$ and $h'$ two subassignments of $F$, we shall denote by $h
\cap h'$ and $h \cup h'$, the subassignments $ C \mapsto h (C)
\cap h' (C)$ and $ C \mapsto h (C) \cup h' (C)$, respectively.
Similarly, we denote by $h' \setminus h$ the subassignment $C
\mapsto h' (C) \setminus h (C)$.

We also write $h' \subset h$ if $h' (C) \subset h (C)$ for every
object $C$ of $\cC$. In the case where
 $h'\subset h$ are
subassignments of $F$ we will also call $h'$ a subassignment of
$h$ (although $h$ itself need not to be a functor).

There is a trivial notion of a morphism between subassignments:
for $h_1$ and $h_2$ subassignments of some functors $F_1$,
$F_2:\cC \rightarrow {\rm Sets}$, a morphism $f:h_1\to h_2$ is
just the datum, for every object $C$ of $\cC$, of a function
$f(C)$ (or $f$ for short) from $h_1(C)$ to $h_2(C)$. If $h_i'$ is
a subassignment of $h_i$, $i=1,2$, one defines the subassignments
$f(h_1')$ and $f^{-1} (h_2')$  in the obvious way. We can also
define the Cartesian product $h_1\times h_2$ of $h_1$ and $h_2$ by
$(h_1\times h_2)(C):=h_1(C)\times h_2(C)$ for every object $C$ of
$\cC$; it is a subassignment of the functor $F_1\times F_2$ which
sends an  object $C$ of $\cC$ to $F_1(C)\times F_2(C)$. Similarly,
one can perform other operations of set theory, for example:

The graph of a morphism $f:h_1\to h_2$ with $h_i$ a subassignment
of $F_i$ is the subassignment of $F_1\times F_2$ sending an object $C$ of $\cC$
to
\[
\{(x,y)\in h_1(C)\times h_2(C)\mid f(x)=y\}.
\]

If $h_i$ for $i=1,2,3$ are  subassignments of $F_i:\cC \rightarrow
{\rm Sets}$ and $f_j:h_j\to h_3$ morphisms for $j=1,2$, the fiber
product $h_1\otimes_{h_3} h_2$ is the subassignment of $F_1\times
F_2$ sending an object $C$ of  $\cC$ to
\[
\{(x,y)\in h_1(C)\times h_2(C)\mid f_1(x)=f_2(y)\}.
\]

\subsection{Definable subassignments}\label{sec:def:sub}
Let $k$ be a field. We denote by ${\rm Field}_{k}$ the category of
all fields containing $k$. More precisely, to avoid any
set-theoretical issue, we shall fix a Grothendieck universe $\cU$
containing $k$ and we define ${\rm Field}_{k}$ as the small category
of all fields in $\cU$ containing $k$.

It goes without saying that if $\ell$ is any other field than $k$,
then ${\rm Field}_{\ell}$ stands for the category of all fields
containing $\ell$. Although $k$ plays the role of standard base
field throughout the paper, the definitions of section \ref{sec1}
make sense over any base field $\ell$ instead of $k$.

We consider $W := \AA^m_{k \llp t \rrp} \times \AA^n_k \times
\ZZ^r$, $m,n,r\geq 0$. It defines a functor $h_W$ from the category
${\rm Field}_k$ to the category of sets by setting $h_W (K) = K\llp
t \rrp^m \times K^n \times \ZZ^r$. We shall write $h[m, n, r]$ for
$h_W$, where the base field $k$ is implicit in the notation $h[m, n,
r]$; thus to avoid confusion, we only use this notation when the
base field is clear. However, we will usually explicitly write
$h_{\Spec k}$ instead of $h[0,0,0]$; it is the functor which assigns
to each $K$ in ${\rm Field}_k$ the one point set.

Fix a language  $\LPas$ of Denef-Pas\footnote{Except in Theorem
\ref{pqe}, we always assume without writing that $(K\llp t \rrp,K,\ZZ)$
is a structure for $\LPas$ for all fields $K$ under consideration
(usually $K$ runs over a category of the form ${\rm Field}_k$). Starting from
section \ref{sec2} the language will be $\LPre$,
which satisfies this condition.}. Any formula $\varphi$ in $\LPas$
with coefficients in $k \llp t \rrp$ in the valued field sort and
coefficients in $k$ in the residue field sort, with $m$ free
variables in the valued field sort, $n$ in the residue field sort
and $r$ in the value group sort, defines a subassignment
$h_{\varphi}$  of $h[m,n,r]$ by assigning to $K$ in ${\rm
Field}_k$ the subset of $h[m,n,r](K)$ defined by $\varphi$,
namely,
\[
h_{\varphi}(K)=\{x\in h[m,n,r](K)\mid (K,K \llp t \rrp,\ZZ)\models
\varphi(x)\}.
\]
We call $h_\varphi$ a definable subassignment of $h[m,n,r]$.

If the coefficients of the formula $\varphi$ in the valued field
sort all lie in some subring $S$ of $k \llp t \rrp$ and the
coefficients in the residue field sort are still allowed to be in
$k$, we call $h_{\varphi}$ a $\LPas(S)$-definable subassignment and
we write $\LPas(S)$ to denote the language $\LPas$ with such
coefficients.

We denote by $\emptyset$  the empty definable subassignment which
sends each $K$ in ${\rm Field}_{k}$ to the empty set $\emptyset$.
Here again, the base field $k$ is implicit in the notation
$\emptyset$, and we only use it when the base field is clear.

More generally, if $\cX$ is a variety, that is, a separated and
reduced scheme of finite type,  over $k \llp t \rrp$ and $X$ is a
variety over $k$, we consider $W' := \cX \times X \times \ZZ^r$
and the functor $h_{W'}$ from ${\rm Field}_k$ to the category of
sets which to $K$ assigns $h_{W'} (K) = \cX (K\llp t \rrp) \times X (K)
\times \ZZ^r$. We will define definable subassignments of $h_{W'}$
by a glueing procedure. Assume first $X$ is affine and embedded as
a closed subscheme in $\AA^{n}_{k}$ and similarly for $\cX$ in
$\AA^{m}_{k\llp t \rrp}$. We shall say a subassignment of $h_{W'}$ is a
definable subassignment if it is a definable subassignment of
$h[m,n,r]$. Clearly, this definition is independent of the choice
of the embedding of $X$ and $\cX$ in affine spaces.

In general, a subassignment $h$ of $h_{W'}$ will be a definable
subassignment if there exist finite covers $(X_{i})_{i \in I}$ of
$X$ and $(\cX_{j})_{j \in J}$ of $\cX$ by affine open subschemes
(defined over $k$ and $k \llp t \rrp$ respectively; such covers always
exist) and definable subassignments $h_{ij}$ of $h_{\cX_{i}\times
X_{j}\times \ZZ^r}$, for $i$ in $I$ and $j$ in $J$, such that $h =
\cup_{i,j} h_{ij}$. If $\cX$ as well as its cover $(\cX_{j})_j$ is
defined over some subring $S$ of $k \llp t \rrp$, and if the $h_{ij}$
are $\LPas(S)$-definable subassignments, we call $h$ a
$\LPas(S)$-definable subassignment.

For $i=1,2$, let $h_i$ be a definable subassignment of $h_{W_i}$
with $W_i=\cX_i \times X_i \times \ZZ^{r_i}$, $\cX_i$ a variety
over $k \llp t \rrp$, and $X_i$ a variety over $k$. A definable morphism
$f:h_1\to h_2$ is a morphism $h_1\to h_2$ (as in section
\ref{sec:sub}) whose graph is a definable subassignment of
$h_{W_1\times W_2}$. If moreover $h_1$, $h_2$, and the graph of
$f$ are $\LPas(S)$-definable subassignments for some subring $S$
of $k \llp t \rrp$, we call $f$ a $\LPas(S)$-definable morphism.

The set-theoretical operations defined above for general
subassignments also work at the level of definable subassignments,
for example, fiber products of definable subassignments are again
definable subassignments.

Sometimes we call a definable morphism a definable function,
especially when the image is $h_{\ZZ^r}$ for some $r$.

\subsection{}\label{fstar}
Using our fixed language $\LPas$ of Denef-Pas, we define the
category of (affine) definable subassignments ${\Def}_k$ (also
written $\Def_k(\LPas)$), as the category whose objects are pairs
$(Z,h [m, n, r])$ with $Z$ a definable subassignment of $h [m, n,
r]$, a morphism between $(Z,h [m, n, r])$ and $(Z', h [m', n',
r'])$ being a definable morphism $Z \rightarrow Z'$, that is, a
morphism of subassignments whose graph is a definable
subassignment of $h [m + m', n +n', r +r']$. Similarly one defines
the category  of (global) definable subassignments $\GDef_k$ (also
written $\GDef_k(\LPas)$), as the category whose objects are pairs
$(Z, h_W)$ with $Z$ a definable subassignment of $h_W$, where $W$
is of the form $ \cX \times X \times \ZZ^r$ with $\cX$ a
$k\llp t \rrp$-variety and $X$ a $k$-variety, a morphism between
$(Z,h_W)$ and $(Z', h_{W'})$ being a definable morphism $Z
\rightarrow Z'$.

More generally if $Z$ is in $\Def_k$, resp.~$\GDef_k$, one
considers the category $\Def_Z$, resp.~$\GDef_Z$, of objects over
$Z$, that is,  objects are definable morphisms $Y \rightarrow Z$
in $\Def_k$, resp.~$\GDef_k$, and a morphism between $Y
\rightarrow Z$ and $Y' \rightarrow Z$ is just a morphism $Y
\rightarrow Y'$ making the obvious diagram commute.

For every morphism $f  : Z  \rightarrow Z'$ in $\Def_k$,
composition with $f$ defines a functor $f_! : \Def_Z \rightarrow
\Def_{Z'}$. Also, fiber product defines a functor $f^{*}:
\Def_{Z'}  \rightarrow \Def_{Z}$,
namely, by sending $Y\to Z'$ to $Y\otimes_{Z'} Z\to Z$.
We use similar notations when $f : Z \rightarrow Z'$ is a morphism
in $\GDef_k$.

Let $Y$ and $Y'$ be in $\Def_k$ (resp.~$\GDef_k$). We write $Y
\times Y'$ for the object $Y \otimes_{h_{\Spec k}} Y'$ of
$\Def_k$ (resp.~$\GDef_k$). We shall also write $Y [m, n, r]$ for
$Y \times h [m, n, r]$.
Nevertheless, the notation $h[m,n,r]$ will only be used for the
definable subassignment defined above, and never for a product
$h\times h[m,n,r]$ when $h$ itself is a definable subassignment.

Note that $h_{\Spec k}$ and $\emptyset$ are respectivly the final
and the initial objects of $\Def_k$ and $\GDef_k$.

For a subring $S$ of $k \llp t \rrp$, we define $\Def_k(\LPas (S))$
as the subcategory of $\Def_k$ whose objects are pairs $(Z,h [m, n,
r])$ with $Z$ a $\LPas (S)$-definable subassignment of $h [m, n,
r]$, and whose morphisms are $\LPas (S)$-definable morphisms.
Similarly we define ${\rm GDef}_k(\LPas (S))$, $\Def_{Z}(\LPas
(S))$, and ${\rm GDef}_{Z}(\LPas (S))$ for some $Z$ is in $\Def_k$
or in  $\GDef_k$, respectively\footnote{When the goal is to
interpolate $p$-adic integrals and $\FF_q\llb t \rrb$-integrals, one
can use the ring $S=\cO\llb t \rrb$ with $\cO$ some ring of
integers. This interpolation as well as a transfer result between
$p$-adic and $\FF_q\llb t \rrb$-integrals is announced in
\cite{miami} \cite{expnote} and detailed in \cite{exp}.}.

\subsection{Extension of scalars}
Let $W=\cX \times X \times \ZZ^r$, with $\cX$ a $k\llp t \rrp$-variety
and $X$ a $k$-variety, let $K$ be in ${\rm Field}_k$, and let $W'= \cX'
\times X' \times \ZZ^r$, with $\cX'=\cX\otimes\Spec K\llp t \rrp$ and
$X'=X\otimes\Spec K$. The extension of scalars functor sends
$Z=(Z_0,h_W)$ in $\GDef_k$ to $Z\otimes h_{\Spec
K}:=(Z_0',h_{W'})$ in $\GDef_K$, where $Z_0'$ is the definable
subassignment of $h_{W'}$ which is given by the same
$\LPas$-formulas as $Z_0$ on affine covers of $\cX'$ and $X'$
which are defined over $k\llp t \rrp$ and $k$, respectively. Using
graphs, one defines similarly the image of morphisms in $\GDef_k$
under extension of scalars, getting a functor of extension of
scalars $\GDef_k\to \GDef_K$.

\subsection{Points on definable subassignments}\label{nnn}
For $Z$ in $\GDef_k$, a point $x$ on $Z$ is by definition a tuple
$x=(x_0,K)$ such that $x_0$ is in $Z (K)$ and $K$ is in ${\rm
Field}_k$. For a point $x=(x_0,K)$ on $Z$ we write $k(x)=K$ and we
call $k(x)$ the residue field of $x$.

Let $f : X \rightarrow Y$ be a morphism in $\Def_k$, with $X=(X_0,h
[m, n, r])$ and $Y=(Y_0,h [m', n', r'])$. Let $\varphi(x,y)$ be the
formula which describes the graph of $f$, where $x$ runs over
$h[m,n,r]$ and $y$ runs over $h[m',n',r']$. For every point
$y=(y_0,k(y))$
 of $Y$, we may consider its fibre $X_{y}$, which is
the object in $\Def_{k(y)}$ defined by the formula
$\varphi(x,y_0)$\footnote{Note that the formula $\varphi(x,y_0)$ has
coefficients in $k(y)$ and $k(y)\llp t \rrp$, which is allowed in
$\Def_{k(y)}$. So, in $\Def_{k(y)}$, the base field is $k(y)$
instead of $k$.
 }.
Taking fibers at $y$ gives rise to a functor $i_y^{\ast} : {\rm
Def}_{Y} \rightarrow \Def_{ k (y)}$.

Fibers of a morphism $f : X \rightarrow Y$ in $\GDef_k$ are
defined similarly via affine covers and we shall use similar notations
as for morphisms in $\Def_k$.

\subsection{$T$-subassignments}\label{theories}

Let $T$ be a theory given by sentences in $\LPas$ with
coefficients in $k$ and $k\llp t \rrp$ (a sentence is a formula without
any free variables). We denote by ${\rm Field}_k(T)$ the category
of fields $F$ over $k$ such that $(F\llp t \rrp,F,\ZZ)$ is a model of
$T$ and whose morphisms are field morphisms. Given a
$k\llp t \rrp$-variety $\cX$, a $k$-variety $X$, and $W=\cX\times
X\times \ZZ^r$, we can restrict the functor $h_W$ as defined above
to ${\rm Field}_k(T)$ and we also write $h_W$ to denote this
functor. We can speak of definable $T$-subassignments of $h_W$ in
exactly the same way as we did above for definable subassignments
of $h_W$. A definable $T$-morphism between $T$-subassignments $Z$
and $Z'$ is also defined accordingly.

We define the category $\GDef_k(\LPas, T)$ of definable
$T$-subassignments as the category whose objects are pairs $(Z,
h_W)$ with $Z$ a definable $T$-subassignment of $h_W$, where $W$ is
of the form $W=\cX\times X\times \ZZ^r$, and $T$-morphisms being
definable $T$-morphisms. One defines similarly $\Def_k(\LPas, T)$.

For $S$ a subring of $k \llp t \rrp$, if one restricts moreover the
coefficients in the valued field sort to $S$, one defines the
categories $\Def_k(\LPas (S),T)$ and ${\rm GDef}_k(\LPas (S),T)$
correspondingly.

\section{Dimension theory for definable
subassignments}\label{sec:dim}

\subsection{}
The \emph{Zariski closure} of a definable subassignment $Z$ of
$h_\cX$ with $\cX$ a variety over $k\llp t \rrp$ is the intersection $W$
of all subvarieties $Y$ of $\cX$ such that $Z\subset h_Y$. We
define the dimension of $Z$ as ${\Kdim} \, Z:=\dim W$ if $W$ is
not empty and as $-\infty$ if $W$ is empty. More generally, if $Z$
is a subassignment of $h_W$ with $W= \cX \times X \times \ZZ^r$,
$\cX$ a variety over $k \llp t \rrp$, and $X$ a variety over $k$, we
define ${\Kdim} \, Z$ as the dimension of the image of $Z$
under the projection $h_W\to h_\cX$.
We shall establish basic properties of this
dimension using work by van den Dries \cite{vdDries},
by Denef and Pas \cite{Pas}, and by Denef and van den Dries \cite{DvdD}.
A similar dimension theory in a setting of first order languages
with analytic functions has been developed by \c{C}elikler in
\cite{Celikler}.

\subsection{}\label{sec:man}
Since $K\llp t \rrp$ is a complete field for any field $K$, we can
use the theory of $K\llp t \rrp$-analytic manifolds as developed for
instance in \cite{Bour}, thus using the $t$-adic topology on $K\llp
t \rrp$.
 By a $K\llp t \rrp$-analytic manifold of dimension $n
\geq 0$, we mean a separated topological space endowed with an
analytic atlas of charts into $K\llp t \rrp^n$. Note that we do not
assume $K\llp t \rrp$-analytic manifolds to have a countable basis
for their topology. We shall consider the empty set as a $K\llp t
\rrp$-analytic manifold of dimension $- \infty$.

For any smooth equidimensional variety $\cX$ over $k\llp t \rrp$ and
for any $K$ in ${\rm Field}_k$ the set $\cX(K\llp t \rrp)$ has a
natural structure of $K\llp t \rrp$-analytic manifold. More
generally, if $\cX$ is a smooth equidimensional variety over $k\llp
t \rrp$, $X$ is a variety over $k$, $K$ is in ${\rm Field}_k$, and
$r$ is in $\NN$, the set $h_{\cX \times X \times \ZZ^r}(K)$ has a
natural structure of $K\llp t \rrp$-analytic manifold as product
manifold of $\cX(K\llp t \rrp)$ and $X(K)\times \ZZ^r$ (the latter
considered as a discrete set).

The following theorem asserts that definable subassignments are
closely related to analytic manifolds; its proof is given below. We
thank \c{C}elikler for his help with the proof.

\begin{theorem}\label{prop:man}
Let $Z$ be a nonempty definable subassignment of $h_W$ and $f:Z\to
h[1,0,0]$ a definable morphism, with $W= \cX \times X \times
\ZZ^r$, $\cX$ a variety over $k \llp t \rrp$, and $X$ a variety over
$k$. Let $\{\cX_i\}_i$ be a finite partition of $\cX$ into smooth
equidimensional varieties.
\begin{itemize}
\item[(i)] There exists a finite partition of $Z$ into definable
subassignments $Z_j$ such that, for each $K$ in ${\rm Field}_k$
and each $j$,   the set $Z_j(K)$ is a $K\llp t \rrp$-analytic
submanifold of $h_{\cX_i\times X\times \ZZ^r}(K)$ for some $i$
only depending on $j$ and such that the restriction
$f_{|Z_j}(K):Z_j(K)\to K\llp t \rrp$ is $K\llp t \rrp$-analytic.
 \item[(ii)] If $Z_i$ is a partition as in (i), then ${\Kdim} \, Z$ equals
 $$
 \max_{j,\, K\in{\rm Field}_k,} (\dim Z_j(K)),
 $$
where $\dim Z_j(K)$ denotes the dimension as a $K\llp t \rrp$-analytic
manifold.
 \item[(iii)]\label{lem:diff} There exists a definable
subassignment $Z'$ of $Z$ satisfying $\Kdim\, (Z\setminus
Z')<\Kdim\, Z$ and such that $Z'(K)$ is a $K\llp t \rrp$-analytic
manifold and $f_{|Z'}(K):Z'(K)\to K\llp t \rrp$ is $K\llp t \rrp$-analytic for
each $K$ in ${\rm Field}_k$.
\end{itemize}
\end{theorem}

\subsection{}\label{sec:top}
Given a subassignment $Z$ of $h_W$ with $W$ as in Theorem
\ref{prop:man} and $\cX\to\cX'$ an embedding of $\cX$ into a
smooth equidimensional  variety $\cX'$ over $k\llp t \rrp$ and $K$
in ${\rm Field}_k$, we endow $Z(K)$ with the induced topology
coming from the manifold structure on $h_{\cX'\times X\times
\ZZ^r}(K)$. This topology is independent of the embedding
$\cX\to\cX'$.  Many notions of general topology have a meaning in
$\GDef_k$, for example, if $f:Z\to Y$ is a definable morphism in
$\GDef_k$, we say $f$ is continuous if for each $K$ in ${\rm
Field}_k$ the map $f(K):Z(K)\to Y(K)$ is continuous. Similarly,
for $Z\subset Y$ in $\GDef_k$, one can construct definable
subassignments ${\rm int}(Z)$ and ${\rm cl}(Z)$ such that for each
$K$ in ${\rm Field}_k$ the set ${\rm int}(Z)(K)$, resp.~${\rm cl}
Z(K)$, is the interior, resp.~the closure, of $Z(K)$ in $Y(K)$.

\begin{theorem}\label{compdim}
Let $Z$ and $Y$ be in $\GDef_k$ and nonempty.
\begin{enumerate}
 \item[(i)] If $f:Z\to Y$ is a definable morphism in $\GDef_k$, then ${\Kdim} \, Z\geq {\Kdim} \, f(Z)$. If $f$ is a definable
isomorphism, then ${\Kdim} \, Z={\Kdim} \, Y$.
 \item[(ii)] The inequality $${\Kdim} \, (Z\times Y)\leq {\Kdim} \,
Z+{\Kdim} \, Y$$ holds.
 \item[(iii)] If $Z$ and $Y$ are definable subassignments of the same subassignment in $\GDef_k$, one has
\[
{\Kdim} \, (Z\cup Y)=\max({\Kdim} \, Z,\ {\Kdim} \, Y).
\]

 \item[(iv)] If $Z\subset Y$, let ${\rm cl}(Z)$ be the definable subassignment  which is the closure
of $Z$ in $Y$ as in \ref{sec:top}. Then
 \[
{\Kdim} \, ({\rm cl}(Z)\setminus Z) < {\Kdim} \, Z.
 \]
 \item[(v)] The integer ${\Kdim} \, Z$ is equal to the largest integer $d$ such that there exists a
 definable morphism $f:Z\to h[d,0,0]$ such that $f(Z)$ has nonempty
 interior in $h[d,0,0]$ for the topology of section \ref{sec:top},
 and with $f(Z)$ the image of $Z$ under $f$.
\end{enumerate}
\end{theorem}
\begin{example} The inequality in (ii) of Theorem
\ref{compdim} can be strict: Suppose that $-1$ does not have a
square root in $k$. Let $Z$, resp.~$Y$, be a definable
subassignment of $h[1,0,0]$ given by the formula $x=x\wedge
\exists y\in h[1,0,0]\ y^2=-1$, resp.~$x=x\wedge \forall y\in
h[1,0,0]\ y^2\not=-1$, where $x$ runs over $h[1,0,0]$. Then the
Zariski closure of both $Z$ and $Y$ is $\AA_{k\llp t \rrp}^1$, hence
they both have dimension $1$, but the definable subassignment
$Z\times Y$ is empty, hence its dimension is $-\infty$.
\end{example}

\begin{proof}[Proof of Theorem \ref{compdim}]
Let first $Z$ be a definable subassignment of $h[m,n,r]$ for some
$m,n,r$, and $Y$ a definable subassignment of $h[m',n',r']$ for
some $m',n',r'$.

Let $\cL$ be the language  $\cL_{\rm Val}$ together with the
following additional relation symbols:
\begin{itemize}
\item[(1)] for each $\bfL_{\rm Ord}$-formula $\varphi$ in $n$ free
variables the $n$-ary relation symbol $R_\varphi$ interpreted as
follows: $R_\varphi(a_1,...,a_n)$ if and only if
$\varphi(\ord(a_1),...,\ord(a_n))$ holds;

\item[(2)] for each $\bfL_{\rm Res}$-formula $\varphi$ in $n$ free
variables the $n$-ary relation symbol $R_\varphi$ interpreted as
follows: $R_\varphi(a_1,...,a_n)$ if and only if
$\varphi(\ac(a_1),...,\ac(a_n))$ holds.
\end{itemize}

For each $K$ in ${\rm Field}_k$ let $\cL_K$ be the language $\cL$
with additional constant symbols for all elements of $K\llp t \rrp$.
Then, for each $K$ in ${\rm Field}_k$, it follows from  Theorem
\ref{pqe} that the structure $(K\llp t \rrp,\cL_K)$ has elimination of
quantifiers; moreover, $(K\llp t \rrp,\cL_K)$ satisfies the conditions
of Proposition 2.15 of \cite{vdDries}, with a topology as in
\ref{sec:top}. Thus, using the terminology of \cite{vdDries},
there is a dimension function ${\rm algdim}$ - defined via the
Zariski closure of definable sets - on  $\cL_K$-definable
subsets of the structure $(K\llp t \rrp,\cL_K)$ for each $K$.

We claim that
 \begin{equation}\label{eq:dim}
 \max_{K\in {\rm Field}_k,\ (y,z)\in K^n\times \ZZ^r}
 {\rm algdim}\, (Z_{(y,z),K})= {\Kdim} \, Z,
 \end{equation}
with $Z_{(y,z),K}= Z(K)\cap (K\llp t \rrp^m\times \{(y,z)\})$, from
which the theorem will follow.

The subassignment $Z$ is given by a formula $\varphi(x,y,z)$,
where $x$ are the ${\rm Val}$-variables, $y$ the ${\rm
Res}$-variables and $z$ the ${\rm Ord}$-variables. By Denef-Pas
quantifier elimination Theorem \ref{pqe}, we can write $\varphi$ as a
disjunction over $j$ of formulas of the form
 \begin{equation}\label{eq:man}\begin{array}{l}
\psi_j(z,\ord\, f_{1j}(x),\ldots,\ord\, f_{rj}(x))\ \wedge_i\ f_{ij}(x)\not=0 \\
\wedge\ \vartheta_j(y,\ac\, g_{1j}(x),\ldots,\ac\, g_{sj}(x))\ \wedge_i\ g_{ij}(x)\not=0 \\
\wedge\ h_{1j}(x)=0\wedge\ldots \wedge h_{tj}(x)=0,\\
\end{array}
\end{equation}
with $f_{ij},g_{ij},h_{ij}$ polynomials
over $k\llp t \rrp$
 - strictly speaking, these polynomials are
defined over the constant symbols of $\cL_{\rm Val}$ - in the
variables $x$, $\psi_j$ $\bfL_{\rm Val}$-formulas and $\vartheta_j$
$\bfL_{\rm Res}$-formulas.
Note that the first two lines of (\ref{eq:man}) determine open
conditions. Let $V_j$ be the variety over $k\llp t \rrp$ associated
to the ideal $(h_{1j},\ldots,h_{tj})$. Just by rewriting the
disjunction (\ref{eq:man}), we may suppose that the $V_j$ are
irreducible over $k\llp t \rrp$.
 We prove
(\ref{eq:dim}) by induction on the maximum of the Zariski dimensions
of the $V_{j}$.
 For each $j$ and each point $w=(y_0,z_0,K)$ in $h[0,n,r]$ let
$U_{jw}$ be the definable subassignment in $\Def_K$ given by
 $$\begin{array}{l}
\psi_j(z_0,\ord(f_{1j}(x)),\ldots,\ord(f_{rj}(x)))\ \wedge_i\ f_{ij}(x)\not=0 \\
\wedge\ \vartheta_j(y_0,\ac(g_{1j}(x)),\ldots,\ac(g_{sj}(x)))\
\wedge_i\ g_{ij}(x)\not=0 \\
\wedge\ h_{1j}(x)=0\wedge\ldots \wedge h_{tj}(x)=0.\\
\end{array}
$$

Let $j_0$ be such that $V_{j_0}$ has maximal Zariski dimension, say
$\ell_0$, among the $V_j$. In the case where
 there exists
$w=(y_0,z_0,K)$ such that the set $U_{j_0w}(K)$ has a Zariski
closure over $K\llp t \rrp$ of dimension equal to $\ell_0$,
(\ref{eq:dim}) follows. In the case where
 for all $w=(y_0,z_0,K)$ the
set $U_{j_0w}(K)$ has a Zariski closure $V_w$ over $K\llp t \rrp$ of
dimension $<\ell_0$, then, %since $K\llp t \rrp[x]$ is faithfully
%flat over $k\llp t \rrp[x]$,
%by the form of the disjunction (\ref{eq:man}) and
since the field of definition of $V_w$ is contained in the algebraic
closure of $k\llp t \rrp$ intersected with $K\llp t \rrp$, for each
$w=(y_0,z_0,K)$, the set $U_{j_0w}(K)$ is contained in a Zariski
closed set defined over $k\llp t \rrp$ of dimension $<\ell_0$.
Hence, by compactness, $V_{j_0}$ can be replaced by a $k\llp t
\rrp$-variety $V'_{j_0}$ of dimension $<\ell_0$. This proves
(\ref{eq:dim}).

Theorem \ref{compdim} now follows for $Z\subset h[m,n,r]$,
$Y\subset h[m',n',r']$ by (\ref{eq:dim}) and by the equivalent
properties of $ {\rm algdim}$  stated in \cite{vdDries}. For
example, to prove (i) we consider the definable morphism
 $$f':=\pi\times f:Z\to h[0,n,r]\times Y:x \longmapsto (\pi(x),f(x))
 $$
with $\pi:Z\to h[0,n,r]$ the projection. We then compute
 \begin{eqnarray*}
{\Kdim}\, (Z) & = &  \max_{K\in {\rm Field}_k,\ x\in K^n\times
\ZZ^r} {\rm algdim}\, (Z_{x,K})\\
 & \geq & \max_{K\in {\rm Field}_k,\ y\in K^{n+n'}\times \ZZ^{r+r'}}
 ({\rm algdim}\, (f'(Z)_{y,K}) ) \\
 & = & {\Kdim} \, f'(Z)\\
 & = &  {\Kdim} \, f(Z),
 \end{eqnarray*}
with $f'(Z)_{y,K}:=Y(K)\cap (K\llp t \rrp^{m'}\times \{y\})$. Indeed,
the inequality holds by the equivalent property for $ {\rm
algdim}$ for fixed $x$, $K$, and $y$ above $x$, and the last
equality follows from the definition of $\Kdim$.

These results extend to any definable subassignments $Z$ and $Y$
of functors of the form $h_W$ with $W= \cX \times X \times \ZZ^r$,
$\cX$ a variety over $k \llp t \rrp$, $X$ a variety over $k$, by using
affine charts on $\cX$ and $X$.
\end{proof}

\begin{proof}[Proof of Theorem \ref{prop:man}]
Theorem \ref{prop:man} for $Z\subset h[m,n,r]$ follows by the same
proof as the proof of Lemma 3.12 of \cite{DvdD}, the first part of
the proof of Lemma 3.18 of \cite{DvdD}, and the proof of
Proposition 3.29 of \cite{DvdD}. For the convenience of the reader
we give an outline of this argument and refer to \cite{DvdD} for
the details.

We use the notation of the proof of Theorem \ref{compdim}. Let
$Z\subset h[m,n,r]$, $\varphi$, and $V_{j}$ be as in the proof of
Theorem \ref{compdim}. If one takes a partition  of each $V_{j}$
into smooth $k\llp t \rrp$-subvarieties which are irreducible over
$k\llp t \rrp$, noting that the first two lines of (\ref{eq:man})
describe an open set, one can easily partition $Z$ into definable
manifolds as in statement (i) by taking appropriate Boolean
combinations (see Lemma 3.12 in \cite{DvdD} for details).

For the part of statement (i) about $f:Z \subset h[m,n,r]\to
h[1,0,0]$ one uses induction on the dimension of $Z$ to obtain a
finite partition of $Z$ such that the restriction of $f$ to each
part is continuous (as in the proof of Proposition 3.29 in
\cite{DvdD}), hence one may suppose that $f$ is continuous. Then
one partitions the graph $\Gamma(f)$ of $f$ into manifolds as in
(i) and one refines the partition in such a way  that the tangent map of the
projection $\pi:\Gamma(f)\to h[m,n,r]$ has constant rank on each
part (as in the first part of the proof of Lemma 3.18 of
\cite{DvdD}). It then follows that on each part of this partition
the map $\pi$ is an analytic isomorphism between manifolds with
analytic inverse $f$.

Statement (ii) for $Z \subset h[m,n,r]$ follows from Theorem
\ref{compdim} (iii) and (\ref{eq:dim}) in its proof.

Statement (iii) for $Z \subset h[m,n,r]$ follows easily from (i)
and Theorem \ref{compdim}, by taking the parts $Z_j$ of maximal
dimension among the parts obtained in (i), and taking the union of
$Z_j\setminus \cup_{i\not=j}{\rm cl}(Z_i)$ for $Z'$.

Again, this extends to any definable subassignment $Z$ of $h_W$
with $W= \cX \times X \times \ZZ^r$, $\cX$ a variety over $k
\llp t \rrp$, $X$ a variety over $k$, by using affine charts on $\cX$
and $X$.
\end{proof}

\subsection{Relative dimension}\label{reld}

Let $Z$ and $Y$ be in $\GDef_k$ and let $f:Z\to Y$ be a definable
morphism. For every point $x$ on $Y$ let $Z_x$ be its fiber,
as defined in
section \ref{nnn}.

For $i$ in $\NN \cup\{- \infty \}$, we say $Z$ is of relative
dimension $\leq i$ rel.~$f$ if ${\Kdim} \, Z_x \leq i$ for every
point $x$ in $Y$. We say $Z$ is equidimensional of relative
dimension $i$ rel.~$f$ if ${\Kdim} \, Z_x = i$ for every point $x$
in $Y$.

By Proposition 1.4 of \cite{vdDries} and by using similar
arguments as the ones in the proof of Theorem \ref{compdim}, we
deduce the following proposition:
\begin{prop}\label{prop:reld} Let $Z$ and $Y$ be in $\GDef_k$ and let $f:Z\to Y$ be a definable
morphism. For every point $x$ on $Y$ let $Z_x$ be its fiber (as in
section \ref{nnn}). The morphism $H:Y\to h_{\ZZ}$ which sends $x$
to ${\Kdim} \, Z_x$ if $Z_x$ is nonempty and to $-1$ otherwise
is a definable morphism. For $i$ in $\NN$ let $Y_i$ be the
definable subassignment of $Y$ given by $H(x)=i$. Then, the
definable
subassignment $f^{-1}(Y_i)$ has dimension $i+{\Kdim} \, Y_i$.
\end{prop}

The next proposition is a relative version of Theorem
\ref{prop:man}.
\begin{prop}\label{prop:man:rel} Let $\Lambda$ be in $\Def_k$.
Let $Z\subset \Lambda[m,n,r]$ be a nonempty definable subassignment
over $\Lambda$ for some $m,n,r$ such that $Z\to \Lambda$ is
surjective. Let $f:Z\to h[1,0,0]$ be a definable morphism.
\begin{itemize}
\item[(i)] There exists a finite partition of $Z$ into definable
subassignments $Z_j$ such that, for each $\lambda$ in $ \Lambda$, for
each $K$ in ${\rm Field}_{k(\lambda)}$, and each $j$, the fiber
$i^*_\lambda(Z_{j})(K)$ is a $K\llp t \rrp$-analytic submanifold of
$K\llp t \rrp^m\times K^n\times \ZZ^r$ and such that the morphism
$$f_{|i^*_\lambda(Z_{j})}(K):i^*_\lambda(Z_{j})(K)\longrightarrow  K\llp t \rrp$$ is
$K\llp t \rrp$-analytic.
 \item[(ii)]\label{lem:diff:rel} There exists a definable
subassignment $Z'$ of $Z$ satisfying $$\Kdim\,
\Bigl(i^*_\lambda(Z)\setminus i^*_\lambda(Z')\Bigr)<\Kdim\,
i^*_\lambda(Z)$$ for each $\lambda$ in $\Lambda$ and such that, for
each $K$ in ${\rm Field}_{k(\lambda)}$, the fiber
$i^*_\lambda(Z')(K)$ is a $K\llp t \rrp$-analytic manifold on which
$f_{|i^*_\lambda(Z')}(K)$ %:i^*_\lambda(Z')(K)\to K\llp t \rrp$
is $K\llp t \rrp$-analytic.
\end{itemize}
\end{prop}
\begin{proof}
Apply Theorem \ref{prop:man} to $Z$, to $f$, and to the (${\rm
Val}$ coordinate functions of the) structure map $g:Z\to\Lambda$.
Partition further into finitely many definable subassignments such
that the restriction of $g$ to each of the parts has constant rank
$d$ with respect to the ${\rm Val}$-variables, and with a constant
submatrix of the Jacobian matrix of $g$ of size $d$ having nonzero
determinant. Now the proposition follows from the Implicit
Function Theorem.
\end{proof}

%\begin{prop}\label{prop:zar:rel}
%Let $\Lambda\subset h[m_0,n_0,r_0]$ be in $\Def_k$. Let $Z\subset
%h[m,n,r]$ be a nonempty definable subassignment over $\Lambda$ for
%some $m,n,r$. Let $Z_\Lambda$ be the graph of the structure
%morphism $Z\to \Lambda$. Let $\cY$ be the Zariski closure of the
%image of $Z_\Lambda$ under the projection to $h[m_0+m,0,0]$. Then
%$Z$ is of relative dimension $d$ over $\Lambda$ if and only if
%$h_\cY$ is of relative dimension $d$ over $h[m_0,0,0]$.
%\end{prop}
%The proof of this proposition will be given in section
%\ref{sec:an:cel}.
%\begin{proof}

%\end{proof}

\section{Summation over Presburger sets}\label{secpre}

\subsection{Presbuger sets}\label{llo}
Let $G$ denote  a $\ZZ$-group, that is, a group which is elementary
equivalent to the integers $\ZZ$ in the Presburger language
$\bfL_{\rm PR}$. We call $(G,\Lp)$ a Presburger structure. By a
Presburger set, function, etc., we mean a $\Lp$-definable set,
function. We recall that the theory ${\rm Th}(\ZZ,\Lp)$ has
quantifier elimination in $\Lp$ and is decidable \cite{Pr}. Let $S$
be a Presburger set. We call a function $$f: X\subset S \times
G^m\to G$$ $S$-linear (or linear for short)
 if there is a definable function
$\gamma$ from $S$ to $G$, and integers $0\leq c_i < n_i$ and $a_i$,
for $i=1,\ldots,m$, such that for every $x = (s, x_1,\ldots,x_m)$ in
$X$, $x_i-c_i\equiv 0 \pmod{n_i}$ and
\begin{equation}\label{linear}
f(x)=\sum_{i=1}^m a_i  \Bigl(\frac{x_i-c_i}{n_i}\Bigr)+ \gamma (s).
\end{equation}
We define similarly $S$-linear  maps $g:X\to G^n$.

\par From now on in this section we shall assume that $G = \ZZ$.

\subsection{Constructible Presburger functions}\label{cpf}
We consider a formal symbol $\LL$ and the ring
$$\AA := \ZZ \Bigl[\LL, \LL^{- 1},
\Bigl(\frac{1}{1 - \LL^{- i}}\Bigr)_{i >0}\Bigr].$$ Note that for
every  real number $q > 1$, there is a unique morphism of rings
$\vartheta_q : \AA \rightarrow \RR$ mapping $\LL$ to $q$ and that,
for $q$ transcendental,  $\vartheta_q$ is injective. We define a
partial ordering of $\AA$ by setting $a \geq b$ if, for every real
number $q > 1$, $\vartheta_q (a) \geq \vartheta_q (b)$. We denote by
$\AA_+$ the set $\{a \in \AA \, | \, a \geq 0\}$. Note that $\LL^i$,
for $i$ in $\ZZ$, $\LL^i - \LL^j$, for $i > j$,  and $\frac{1}{1 -
\LL^{- i}}$, for $i >0$, all lie in $\AA_+$, but, for instance, $\LL
- 2$ does not.
One has $a=b$ in $\AA$ if and only if
$\vartheta_q(a)=\vartheta_q(b)$ for all $q>1$. Indeed, considering a
single transcendental $q>1$ is enough.

Now if $S$ is a definable subset of $\ZZ^m$ we define the ring $\cP
(S) $ of constructible Presburger functions on $S$ as the subring of
the ring of functions $S \rightarrow \AA$ generated by all constant
functions into $\AA$, all definable functions $S \rightarrow \ZZ$
and all functions of the form $\LL^{\beta}$ with $\beta$ a
$\ZZ$-valued definable function on $S$. We denote by $\cP_+ (S) $
the semiring of functions in $\cP (S) $ with values in $\AA_+$ and
write $f \geq g$ if $f -  g$ is in $\cP_+ (S) $. This defines a
partial ordering on $\cP (S)$. When $S$ is one point we identify
$\cP (S)$ and $\AA$.

\subsection{Cell decomposition for Presburger sets}\label{cdps}
In this subsection we recall the cell decomposition for Presburger sets
as presented in
\cite{pres}.
Let $G$ denote  a $\ZZ$-group. Fix a Presburger set $S$.
We  define Presburger cells parametrized by $S$, or Presburger
$S$-cells.
 \begin{definition}\label{def cell}
An $S$-cell of type $(0)$ (also called a $(0)$-cell or cell for
short) is a subset of $S \times G$ which is the graph of a
$S$-linear function $S' \rightarrow G$,  with $S'$ a definable
subset of $S$. An $S$-cell of type $(1)$ (also called $(1)$-cell or
cell for short) is a subset $A$ of $S \times G$ of the form
 \begin{equation}\label{cell dim 1}
  \{(s, x) \in S' \times G\mid \alpha (s) \sq_1 x\sq_2 \beta (s),\ x\equiv c \pmod{n}\},
 \end{equation}
with $S'$ a definable subset of $S$, $\alpha$ and $\beta$ $S$-linear
functions $S' \rightarrow G$, $c$ and $n$ integers such that $0\leq
c<n$, and $\sq_i$ either $\leq$ or no condition, and such that the
cardinality of the fibers $A_s=\{x \in G\mid (s,x)\in A\}$ cannot be
bounded uniformly for $s$ in $S'$ by an integer.

Let us consider $i_j$ in $\{0,1\}$, for
$j=1,\ldots,m$, and $x=(x_1,\ldots,x_m)$. A
$(i_1,\ldots,i_{m},1)$-cell is a subset $A$ of $S \times G^{m + 1}$
of the form
 \begin{equation}\label{cell}
  A = \{(x,t)\in S \times G^{m+1} \mid x\in D,\ \alpha(x)\sq_1 t \sq_2
  \beta(x),\ t\equiv c \pmod{n}\},
 \end{equation}
with $D=\pi_{m}(A)$ a $(i_1,\ldots,i_m)$-cell in $S \times G^m$,
$\pi_m$ denoting the projection $S \times G^{m + 1} \rightarrow S \times G^m$,
 $\alpha,\ \beta:D\to G$ $S$-linear functions, $\sq_i$ either
$\leq$ or no condition and integers $0\leq c<n$ such that the
cardinality of the fibers $A_x=\{t\in G\mid (x,t)\in A\}$ cannot
be bounded uniformly for  $x$ in $D$ by an integer.
 \\
 A $(i_1,\ldots,i_{m},0)$-cell is a set of the form
 \[
\{(x,t)\in S \times G^{m+1} \mid x\in D,\ \alpha(x)=t\},
 \]
with $\alpha:D\to G$ a $S$-linear function and $D\subset S \times G^{m}$ a
$(i_1,\ldots,i_m)$-cell.
\end{definition}
A subset of $S \times G^m$ is called a $S$-cell if it is a
$(i_1,\ldots,i_{m})$-cell
for some $i_j$ in $\{0, 1\}$.

Now we can state the following:
\begin{theorem}[Presburger Cell Decomposition \cite{pres}]\label{cell decomp}
Let $S$ be a $\Lp$-definable set, let $X$ be a $\Lp$-definable
subset of $S \times G^m$ and $f:X\to G$ a $\Lp$-definable map.
Then there exists a finite partition $\Pm$ of $X$ into $S$-cells,
such that the restriction $f|_A:A\to G$ is $S$-linear for every
cell $A$ in $\Pm$.
\end{theorem}

\begin{remark}Of course, one could assume $S$ is the one point set
in the above statement, but it is more convenient to express it that way,
in view of further generalizations.
\end{remark}

\subsection{The basic rationality result}\label{brs}
Let $S$ be a definable Presburger set. We consider the ring $\cP
(S) \llb T_1, \cdots, T_r \rrb$ of formal series with coefficients in
the ring $\cP (S)$. If $\alpha$ is a definable function on $S$
with values in $\NN^r$, we write $T^{\alpha}$ for the series
$\sum_{j \in \NN^r} \11_{C_j} T^j$ in $\cP (S) \llb T_1, \cdots,
T_r \rrb$, where $\11_{C_j}$ is the characteristic function of the
subset $C_j$ of $S$ defined by the formula $\alpha(x)=j$. We
consider the subring $\cP (S) \{T_1, \cdots, T_r\}$ of power
series of the form $\sum_{i \in I} a_i T^{\alpha_i}$ with $I$
finite, $a_i$ in $\cP (S)$, and $\alpha_i$ a definable function on
$S$ with values in $\NN^r$. In other words $\cP (S) \{T_1, \cdots,
T_r\}$ is the $\cP (S)$-subalgebra of $\cP (S) \llb T_1, \cdots,
T_r \rrb$ generated by elements of the form $T^{\alpha}$ with $\alpha
: S \rightarrow \NN^r$ definable.

We denote by $\Gamma$ the multiplicative set of polynomials in $\cP
(S) [T_1, \cdots, T_r]$ generated by the polynomials $1 -
\LL^{\alpha} T^{\beta} :=1 - \LL^{\alpha} \prod_{1 \leq i \leq r}
T_i^{\beta_i}$, for $\alpha$ in $\ZZ$ and $\beta = (\beta_1, \cdots,
\beta_r)$ in $\NN^r \setminus \{(0, \cdots, 0)\}$. We denote by $\cP
(S) \{T_1, \cdots, T_r\}_{\Gamma}$ the localisation of $\cP (S)
\{T_1, \cdots, T_r\}$ with respect to $\Gamma$. Since the
polynomials $1 - \LL^{\alpha} T^{\beta}$ are invertible in $\cP (S)
\llb T_1, \cdots, T_r \rrb$, there exists a canonical morphism of
rings
$$
\cP (S) \{T_1, \cdots, T_r\}_{\Gamma}
\longrightarrow
\cP (S) \llb T_1, \cdots, T_r \rrb,
$$
which is injective. We denote by $\cP (S) \llb T_1, \cdots,
T_r \rrb_{\Gamma}$, or by $\cP (S) \llb T \rrb_{\Gamma}$ for short, the
image of this morphism, which we identify with $\cP (S) \{T_1,
\cdots, T_r\}_{\Gamma}$.

We shall consider the $\cP (S)$-module $\cP (S) \llb T_1,T_1^{-1}
\cdots, T_r, T_r^{-1} \rrb$, or $\cP (S) \llb T,T^{-1} \rrb$ for short.
Note that there is a natural product on $\cP (S) \llb T,T^{-1} \rrb$,
the Hadamard product, defined by
$$
f \ast g = \sum_{i \in \ZZ^r} f_i g_i T^i,
$$
for $f = \sum_{i \in \ZZ^r} f_i  T^i$ and $g = \sum_{i \in \ZZ^r}
g_i T^i$, that endows $\cP (S) \llb T,T^{-1} \rrb$ with a ring
structure.

For $\varepsilon$ in $ \{+1, -1\}^{r}$, we denote by $\varepsilon_*$
the $\cP (S)$-module automorphism of $\cP (S) \llb T, T^{-1} \rrb$ that
sends $\sum_{i\in \ZZ^r}a_iT^i$ to $\sum_{i\in \ZZ^r}a_i
T^{\varepsilon i}$ with $\varepsilon i=(\varepsilon_1
i_1,\ldots,\varepsilon_r i_r)$. We denote by $\cP (S)
\llb T,T^{-1} \rrb_{\Gamma}$ the $\cP (S)$-submodule of $\cP (S) \llb T,
T^{-1} \rrb$ generated by the submodules $\varepsilon_* (\cP (S)
\llb T \rrb_{\Gamma})$ for all $\varepsilon$ in $\{+1, -1\}^{r}$.

For $\varphi$ in $\cP (S \times \ZZ^r)$,  and $i$ in $\ZZ^r$, we
shall write $\varphi_i$ for the  restriction of $\varphi$ to $S
\times \{i \}$, viewed as an element of $\cP (S)$, and consider
the series
$$\cM  (\varphi) := \sum_{i \in \ZZ^r} \varphi_i T^i
$$
in $\cP (S) \llb T,T^{-1} \rrb$.

\begin{theorem}\label{presrat}Let $S$
be a definable set. For every $\varphi$ in $\cP (S \times \ZZ^r)$,
the series $\cM (\varphi)$ belongs to $\cP (S)
\llb T,T^{-1} \rrb_{\Gamma}$. Furthermore, the mapping $\varphi \mapsto
\cM (\varphi)$ induces $\cP (S)$-algebra isomorphisms
$$\cM : \cP (S \times \ZZ^r) \longrightarrow \cP (S) \llb T,T^{-1} \rrb_{\Gamma}$$
and
 $$ \cM:\cP (S \times \NN^r) \longrightarrow \cP (S)
\llb T \rrb_{\Gamma},$$ the product on the power series rings being
the Hadamard product.
\end{theorem}

\begin{remark}
Note  that $\cP (S) \llb T,T^{-1} \rrb_{\Gamma}$ and $\cP (S)
\llb T \rrb_{\Gamma}$ are stable by Hadamard product since $\cM$ is a
bijection.
\end{remark}

\begin{proof}The proof of the first statement is quite easy using the cell
decomposition for Presburger sets recalled in Theorem \ref{cell
decomp} and quite similar statements (compare with Lemma 3.2 of
\cite{D85}) may be found in the literature. Let us now prove that
the $\cM$ are isomorphisms. Take $\varphi$ in $\cP (S \times
\ZZ^r)$. We may first assume the support of $\varphi$ is contained
in $S \times \NN^r$. By the Cell Decomposition Theorem we may
furthermore assume that the support of $\varphi$ is contained in a
$S$-cell $A$ and that the restriction of $\varphi$ to $A$ is of the
form $\prod_{1 \leq k \leq d} \alpha_k \LL^{\beta}$ where $\alpha_k$
and $\beta$ are $S$-linear functions on $A$. Let us first consider
the case $r = 1$. When $A$ is a $(0)$-cell, there is nothing to
prove. Assume now $A$ is a $(1)$-cell. Consider first the case where
there is no  condition $\sq_2$ in (\ref{cell dim 1}). By Lemma
\ref{trl}, we can perform a direct computation of $\cM (\varphi)$
(which essentially amounts to summing (derivatives) of geometric
series of monomials in  $\LL$ and $T_1$ along an infinite arithmetic
progression) which yields that $\cM (\varphi)$ is a finite sum of
terms of the form $\psi (s) \frac{T_1^{\gamma (s)}}{(1 - \LL^a
T_1^b)^c}$, with $\psi  (s)$ in $\cP (S)$, $\gamma : S \rightarrow
\NN$ $S$-linear, $a$ in $ \ZZ$, $b > 0$ and $c$ in $\NN$.

\begin{lem}\label{trl} Let $R$ be a ring
and let $P$ be a degree $d$ polynomial in $R [X]$. The equality
\begin{equation}\label{delta}
\sum_{n \geq a} P (n) T^n
=
\sum_{i = 0}^d \frac{[\Delta^i P (a)] T^{ a + i}}{(1 - T)^{i + 1}}
\end{equation}
holds in $R\llb T \rrb$ for all $a$ in $ \NN$. Here $\Delta^i$ is the
$i$-th iterate of the difference operator $P \mapsto P (X + 1) - P
(X)$ with the convention $\Delta^0 P = P$. \qed
\end{lem}

When there is a condition $\sq_2$ in (\ref{cell dim 1}),
we may express $\cM (\varphi)$ as the difference of two series of the preceding type.

Consider now the case $r = 2$. Let us first sum with respect to
the variable $T_2$ in the series $\cM (\varphi)$. By what we know
about the case $r= 1$, relatively to $S \times \NN$, we get that
$\cM (\varphi)$ is a finite sum of terms of the form $$\sum_{i \in
\NN}  \psi (s, i) \, T_1^i \, \frac{T_2^{\gamma (s, i)}}{(1 - \LL^a
T_2^b)^c},$$ with $\psi  (s, i)$ in $\cP (S \times \NN)$, $\gamma :
S \times \NN \rightarrow \NN$ definable, $a$ in $ \ZZ$, $b > 0$
and $c$ in $\NN$. So we just need summing up series of the type
$$\sum_{i \in \NN} \psi (s, i) \, T_1^i  \, T_2^{\gamma (s, i)},$$
which can be done exactly in the same way as the case $r = 1$,
except that instead of dealing with geometric series in monomials
of $\LL$ and $T_1$, we have now to deal with geometric series in
monomials of $\LL$, $T_1$ and $T_2$ which will have the effect of
producing denominators of the form $1 - \LL^a T_1^b T_2^c$, with
$a$ in $\ZZ$, and $b$ and $c$ strictly positive integers.

 The case
where $\varphi$ belongs to  $\cP (S \times \ZZ^r)$ for general $r$
is completely similar. Next we show that every element $m$ in $\cP
(S) \llb T,T^{-1} \rrb_{\Gamma}$ is of the form $\cM (\varphi)$,
by using Theorem-Definition \ref{thm:kjh}. Note that this part of
the theorem is not used in the proof of Theorem-Definition
\ref{thm:kjh}.
 We may assume $m$ is
in $\cP (S) \llb T \rrb_{\Gamma}$, and furthermore, by linearity, that
it is of the form
$$\frac{T^{\alpha}}{\prod_{i=1}^\ell(1-\LL^{a_i}T^{b_i})}$$ with
$\alpha : S \rightarrow \NN^r$ definable, $a_i$ in $\ZZ$, and
$b_i$ in $\NN^r\setminus \{(0,\ldots,0)\}$. The case of general
$\alpha$ following easily from the case $\alpha=0$,  we may assume
$\alpha=0$. Let $H$ be the definable subassignment of $S\times
\NN^{r+\ell}$ given by the condition $$(s,c,d)\in
S\times\NN^{r+\ell}\wedge \sum_{i=1}^\ell d_i\cdot b_i=c.$$ Note
that for each $(s,c)$, the set of $d$'s in $\NN^{\ell}$ satisfying
this condition is finite. Let $\beta:H\to\ZZ$ be the definable
morphism $(s,c,d)\mapsto \sum_{i=1}^\ell d_i a_i$. Let $\psi$ be
$\LL^\beta$ in $\cP (S \times \NN^{r+\ell})$. By the finiteness of
the fibers, $\psi$ is $S\times\NN^r$-summable in the sense of
\ref{poi}. If one sets $\varphi$ to be $\mu_{S\times\NN^r}(\psi)$,
as given by Theorem-Definition \ref{thm:kjh}, then, by
construction, $\cM(\varphi)=m$,  which is what we had  to prove.
\end{proof}

\subsection{Summation of constructible Presburger functions}\label{poi}

Recall the notion of summable families in $\RR$ or $\CC$,
cf.~\cite{bbk} VII.16. In particular, a family $(z_i)_{i \in I}$ of
complex numbers is summable if and only if the family $(|z_i|)_{i
\in I}$ is summable in $\RR$. We shall say a  family $(a_i)_{i \in
I}$ in $\AA$ is summable if, for every $q > 1$, the family
$(\vartheta_q (a_i))_{i \in I}$ is summable in $\RR$. We shall say a
function $\varphi$ in $\cP (S \times \ZZ^r)$ is $S$-integrable if,
for every $s$ in $S$, the family $(\varphi (s, i))_{i \in \ZZ^r}$ is
summable. We shall denote by ${\rm I}_S\cP (S \times \ZZ^r)$ the
$\cP (S)$-module of $S$-integrable functions.

\begin{def-theorem}\label{thm:kjh}
For each $\varphi$ in ${\rm I}_S\cP (S \times \ZZ^r)$ there exists
a unique function $\mu_S(\varphi)$ in $\cP (S)$ such that for all
$q
> 1$ and all $s$ in $S$
\begin{equation}\label{kjh}
\vartheta_q (\mu_S(\varphi) (s)) = \sum_{i \in \ZZ^r} \vartheta_q
(\varphi (s, i)).
\end{equation}
Moreover, the mapping $\varphi\mapsto \mu_S(\varphi)$ yields a
morphism of $\cP (S)$-modules
$$\mu_S : {\rm I}_S \cP (S \times
\ZZ^r) \longrightarrow \cP (S).$$
\end{def-theorem}
\begin{proof}
By induction and  Fubini's Theorem, it is enough to consider the case when
$r=1$.
 Using Theorem \ref{cell decomp} as in the proof
of Theorem \ref{presrat}, we may assume that $\varphi$ is of the
form
\begin{equation}
\varphi (s, i) =
\begin{cases}
 \LL^{\frac{a(i-c)}{n}} \, \Bigl (\frac{i-c}{n}\Bigr)^{b} \, h(s) &
 \text{if $i$ belongs to $I(s)$}\\
 0 & \text{otherwise},
\end{cases}
\end{equation}
where $a$ lies in  $\ZZ$, $b$ in $\NN$, $s$ in $S$, $h$ in $\cP(S)$,
and $$I(s)= \{i\in\ZZ\mid \alpha(s) \sq_1 i \sq_2 \beta(s),\ i\equiv
c \bmod n\},$$ with $0\leq c < n $ integers, $\alpha,\beta:S\to\ZZ$
definable functions, $\sq_i$ either $<$ or no condition, and that
the sum
\begin{equation}\label{eq:cons}
\sum_{i\in I(s)} q^{\frac{a(i-c)}{n}} \, \Bigl
(\frac{i-c}{n}\Bigr)^{b} \, h(s)
\end{equation}
is summable for all $s$ in $ S$ and $q>1$. Now the theorem follows
from Lemma \ref{trl} and Lemma \ref{trl2}, which is a refinement of
the claim in the proof of Lemma 3.2 of \cite{D85}.
\end{proof}

\begin{lem}\label{trl2} Let $b$ and $0\leq c<n$ be integers.
There exist Presburger functions $\gamma_{\ell j}:\ZZ\to \ZZ$ for
$\ell$ in a finite set $L$ and $j=0,\ldots, b$,  such that
$$
\sum_{\sur{0\leq i \leq a} {i\equiv c\bmod n}} i^b=\sum_{\ell\in
L}\prod_{j=0,\ldots,b}\gamma_{\ell j}(a),
$$
for each $a\geq 0$. \qed
\end{lem}

\subsubsection{Characterization of ${\rm I}_S\cP $ in terms of power series}\label{jihy}
We denote by $\cP (S) \{\!\{T\}\!\}$ the subring of $\cP (S) \llb T
\rrb_{\Gamma}$ consisting of series with coefficients in $\cP (S)$
such that, for every $s$ in $S$, at most a finite number of
coefficients have non zero value at $s$. For instance, for $\gamma :
S \rightarrow \NN$, the series $\frac{1 - T^{\gamma}}{1 - T}$
belongs to $\cP (S) \{\!\{T\}\!\}$, say if $r = 1$. Let $\Sigma $ be
 the multiplicative set generated by the polynomials $1 -
\LL^{\alpha} T^{\beta} :=1 - \LL^{\alpha} \prod_{1 \leq i \leq r}
T_i^{\beta_i}$, for $\alpha$ in $\ZZ \setminus \NN$ and $\beta =
(\beta_1, \cdots, \beta_r)$ in $\NN^r\setminus\{(0,\ldots,0)\}$.
 We denote by $\cP (S) \llb T \rrb_{\Sigma}$ the subring of $\cP (S)
\llb T \rrb_{\Gamma}$
 (with the product being the usual product in power series rings)
  whose elements are of the form $\frac{P}{Q}$
with $P$ in $\cP (S) \{\!\{T\}\!\}$ and $Q$ in $\Sigma$. The ring
$\cP (S) \llb T \rrb_{\Sigma}$ captures the summable series among
the series in $\cP (S) \llb T \rrb_{\Gamma}$, as is shown by Theorem
\ref{sintmel}. We shall denote by $\cP (S) \llb T, T^{-1}
\rrb_{\Sigma}$ the $\cP (S)$-submodule of $\cP (S) \llb T, T^{-1}
\rrb_{\Gamma}$ generated by the submodules $\varepsilon_* (\cP (S)
\llb T \rrb_{\Sigma})$ for all $\varepsilon$ in $\{+1, -1\}^{r}$.

\begin{theorem}\label{sintmel}Let $S$ be a definable set.
The transformation $\cM$ induces isomorphisms of $\cP (S)$-modules
$$
\cM : {\rm I}_S\cP (S \times \ZZ^r) \longrightarrow \cP (S) \llb T,
T^{-1} \rrb_{\Sigma}$$
and $$\cM : {\rm I}_S\cP (S \times \NN^r) \longrightarrow \cP (S)
\llb T \rrb_{\Sigma}.$$
\end{theorem}

\begin{proof}Let $\varphi$ be in
$\cP (S \times \ZZ^r)$. We want to prove that $\varphi$ is $S$-integrable
if and only if $\cM (\varphi)$ lies in
$\cP (S) \llb T, T^{-1} \rrb_{\Sigma}$. We may assume the support of
$\varphi$ is contained in $S \times \NN^r$, so that $\cM (\varphi)$ belongs
to
$\cP (S) \llb T \rrb_{\Gamma}$.
It is quite clear that if
$\cM (\varphi)$ belongs
to
$\cP (S) \llb T \rrb_{\Sigma}$, then
$\varphi$ is $S$-integrable.
Assume now
$\cM (\varphi)$ is not in
$\cP (S) \llb T \rrb_{\Sigma}$.
Then, there exists $s_0$ in $S$ and $q > 1$ such that,
extending $\vartheta_q$ coefficientwise to series,
$$\vartheta_q (\cM (\varphi))_{|s = s_0} =
\frac{P_{s_0} (T_1, \dots, T_r)}{Q_{s_0} (T_1, \dots, T_r)},$$ with
$P_{s_0}$ and $Q_{s_0}$ in $\RR [T_1, \dots, T_r]$, such that
$P_{s_0}$ and $Q_{s_0}$ have no non constant common factor in $\RR
[T_1, \dots, T_r]$, and such that for some $\alpha \geq 0$ and
$\beta_i \geq 0$, $\beta \not= (0, \cdots, 0)$, the polynomials
$Q_{s_0}$ and $1 -q^{\alpha} T_1^{\beta_1} \dots T_r^{\beta_r}$ have
a non constant common factor in $\RR [T_1, \dots, T_r]$. Indeed,
otherwise, since one can take $q$ to be transcendental, for every
$s_0$ in $S$, one could write $\cM (\varphi)_{|s = s_0}$ as a
quotient $\frac{P_{s_0}}{Q_{s_0}}$ of polynomials in $\AA [T]$, with
$Q_{s_0}$ in $\Sigma$. Since the polynomials $Q_{s_0}$ all divide a
fixed non zero polynomial in $\AA [T]$, they have a common multiple
$Q$ in $\Sigma$, so we can assume $Q_{s_0}$ is independent of $s_0$,
hence  there exists $P$ in $\cP (S) \{\!\{T\}\!\}$ which gives
$P_{s_0}$ when evaluated at $s_0$ for every $s_0$ in $S$, and $\cM
(\varphi)$ would belong to $\cP (S) \llb T \rrb_{\Sigma}$.

It follows there exists $z_1$, \dots, $z_r$ in $\CC$, with
$|z_1| \leq 1$, \dots, $|z_r| \leq 1$,
such that
$P_{s_0} (z_1, \dots, z_r) \not=0$
and
$Q_{s_0} (z_1, \dots, z_r) = 0$.
In particular the family
$(\vartheta_q (\varphi (s_0, i)))_{i \in \NN^r}$
cannot be summable, since
the summability of a family
of real numbers $(a_i)_{i \in \NN^r}$
implies that the series $\sum_{i \in \NN^r} {a_i} z^i$ is convergent for
every $z = (z_1, \dots, z_r)$ in $\CC^r$ with
$|z_1| \leq 1$, \dots, $|z_r| \leq 1$.
\end{proof}

\subsubsection{}\label{Ldeg}
We may also characterise $S$-integrability in terms of the
$\LL$-degree as follows. We consider the unique extension $\dl : \AA
\rightarrow \ZZ \cup \{ - \infty\}$ of the function degree in $\LL$
from $\ZZ [\LL]$ to $\ZZ \cup \{ - \infty\}$ which satisfies $\dl
(ab) = \dl (a) + \dl (b)$. We have $\dl (a + b) \leq {\rm sup} (\dl
(a), \dl (b))$, with equality if $a$ and $b$ are both in $\AA_+$.
Now if $S  $ is a definable set and $\varphi$ is a function in $\cP
(S)$, we denote by $\dl (\varphi)$ the function $ S \rightarrow \ZZ
\cup \{ - \infty\}$ which sends $s$ to $\dl (\varphi (s))$.

\begin{prop}\label{ldegint}
The following conditions are equivalent for
a function  $\varphi$
in  $\cP (S \times \ZZ^r)$:
\begin{enumerate}
\item[(i)]$\varphi$ is $S$-integrable. \item[(ii)]For every $s$ in
$S$, $\lim_{\vert x \vert \mapsto \infty} \dl (\varphi ( s, x) ) =
- \infty$, where $\vert x \vert$ stands for $\vert x_1\vert +
\dots + \vert x_r\vert$. \item[(iii)] For every $q>1$,
$\lim_{\vert x \vert \mapsto
\infty} \vartheta_q(\varphi( s, x)) = 0$.
\end{enumerate}
\end{prop}

\begin{proof}
Take $\varphi$ in $\cP (S \times \ZZ^r)$. If $\varphi$ is
$S$-integrable,  we know by Theorem \ref{sintmel} that $\cM
(\varphi)$ is in $\cP (S) \llb T, T^{-1} \rrb_{\Sigma}$. But if $\cM
(\varphi)$ belongs to  $\cP (S) \llb T, T^{-1} \rrb_{\Sigma}$, then
the condition that for every $s$ in $S$, $\lim_{\vert x \vert
\mapsto \infty} \dl (\varphi)  ( s, x) = - \infty$, clearly holds.
For the reverse implication, we may by the Cell Decomposition
Theorem assume that the support of $\varphi$ is contained in a
$S$-cell $A$ and that the restriction of $\varphi$ to $A$ is of the
form $\prod_{1 \leq k \leq d} \alpha_k \LL^{\beta}$ where $\alpha_k$
and $\beta$ are $S$-linear functions on $A$ and that furthermore,
for fixed $s$, $\lim_{\vert x \vert \mapsto \infty} \beta (s, x) =
-\infty$. These conditions clearly imply the summability of the
corresponding series. This proves the equivalence of (i) and (ii).
The equivalence of (ii) and (iii) is clear.
\end{proof}

We have the following statement of Fubini type:
\begin{lem}\label{fubpres}Let $S$ be a definable set
and let $\varphi$ be in $\cP (S \times \ZZ^r)$.
Write $r = r_1 + r_2$ and identify $\ZZ^r$ with
$\ZZ^{r_1} \times \ZZ^{r_2}$.
\begin{enumerate}
\item[(1)]If $\varphi$ is $S$-integrable, then $\varphi$, as a
function in $\cP (S \times \ZZ^{r_1} \times \ZZ^{r_2})$, is $S
\times \ZZ^{r_1}$-integrable, $\mu_{S \times \ZZ^{r_1}} (\varphi)$
is $S$-integrable and
$$
\mu_S (\mu_{S \times \ZZ^{r_1}} (\varphi))
=
\mu_{S} (\varphi).
$$
\item[(2)]Assume $\varphi$ is in $\cP_+ (S \times \ZZ^r)$.
Then $\varphi$ is $S$-integrable if and only if it is $S \times
\ZZ^{r_1}$-integrable and $\mu_{S \times \ZZ^{r_1}} (\varphi)$ is
$S$-integrable.
\end{enumerate}
\end{lem}

\begin{proof}In view of Theorem-Definition \ref{thm:kjh}, the statement amounts
to the fact that if a family of real numbers $(a_{i, j})_{(i, j) \in \ZZ^{r_1} \times \ZZ^{r_2}}$ is summable then, for every $i$,
$(a_{i, j})_{j \in \ZZ^{r_2}}$ is summable, the family $(b_i =
\sum_{j \in \ZZ^{r_2}} a_{i, j})_{i \in \ZZ^{r_1}}$ is summable,
$\sum_{i \in \ZZ^{r_1}}b_j = \sum_{(i, j) \in \ZZ^{r_1} \times
\ZZ^{r_2}} a_{i, j}$, and that the reverse statement holds if the
$a_{i, j}$'s are all in $\RR_+$.
\end{proof}

Let $\lambda : S \times \ZZ^r \rightarrow S \times \ZZ^s$ be a
definable function commuting with the projections to $S$. Let $Z$
be a definable subset of $S \times \ZZ^r $ on which $\lambda$ is
injective. Let $\varphi$ be a function in $\cP (S \times \ZZ^r)$
which is zero outside $Z$. We define the function $\lambda_+
(\varphi)$ on $S \times \ZZ^s$ by $\lambda_+ (\varphi) (\lambda
(s, i)) = \varphi (s, i)$ and $\lambda_+ (\varphi) (s, j)= 0$ if
$(s, j)$ does not lie in the image of $\lambda$. Clearly
$\lambda_+ (\varphi)$ lies in $\cP (S \times \ZZ^{s})$.

The following statement will be useful in the proof of the change
of variable formula.

\begin{lem}\label{presinj}
Let $S$ be a definable set and let $\lambda : S \times \ZZ^r
\rightarrow S \times \ZZ^s$ be a definable function commuting with
the projections to $S$. Let $Z$ be a definable subset of $S \times
\ZZ^r $ on which $\lambda$ is injective. Let $\varphi$ be a
function in $\cP (S \times \ZZ^r)$ which is zero outside $Z$. Then
$\varphi$ is $S$-integrable if and only if $\lambda_+ (\varphi)$
is $S$-integrable. Furthermore, if these conditions hold, then
$$\mu_S (\varphi) = \mu_S (\lambda_+ (\varphi)).$$
\end{lem}

\begin{proof}The first statement follows directly
from the definition of $S$-integrability and the second from
Theorem-Definition \ref{thm:kjh}.
\end{proof}

\subsection{Generalization: From Presburger
sets to definable subassignments}\label{uio}
Note that any Presburger subset  $S$ of $\ZZ^m$
is  clearly also $\LPre$-definable.
Furthermore, it  follows from the Denef-Pas quantifier
elimination Theorem \ref{pqe} that
a function
$f : S \rightarrow \ZZ$ is
$\LPre$-definable if and only if it is
$\Lp$-definable.

Let us generalize what we did in
\ref{llo}-\ref{poi} for Presburger subsets
to definable subassignments in
$\GDef_k (\LPre)$.

Let $S$ be a definable subassignment in $\GDef_k (\LPre)$. We denote
by $|S|$ its set of points, defined in \ref{nnn} (this is indeed a
set by \ref{sec:def:sub}). Note that to any definable morphism
$\alpha:S\to h[0,0,1]$ corresponds a function $\tilde \alpha: |S|
\to \ZZ$. (Since $\alpha$ and $\tilde \alpha$ determine each other
we shall not distinguish their notation after \ref{uio}.)
 We define the ring $\cP (S) $ of constructible Presburger
functions on $S$  as the subring of the ring of functions $|S|
\rightarrow \AA$ generated by constant functions $|S| \rightarrow
\AA$, and by functions $\tilde \alpha :|S| \rightarrow\ZZ$ and
$\LL^{\tilde\beta}:|S| \rightarrow \AA$ for definable morphisms
$\alpha,\beta:S \rightarrow h[0,0,1]$.
  We also denote by $\cP_+
(S) $ the semiring of functions in $\cP (S) $ with values in
$\AA_+$. Everything  we did in \ref{llo}-\ref{poi}, including the
proof of Theorem \ref{cell decomp}, generalizes mutatis mutandis  to
that more general situation, up to minor changes like replacing $S
\times \ZZ^r$ by $S \times h_{\ZZ^r} = S [0,0,r]$ or $s$ in $S$ by
$s$ in $|S|$.
 To give an example of this adaptation, an element $\varphi$ of
$\cP(S [0,0,r])$ lies in ${\rm I}_S\cP (S[0,0,r])$ if and only if,
for every $s$ in $|S|$, the family $(\varphi (s, i))_{i \in \ZZ^r}$
is summable in the sense of \ref{poi}.
  Thus it is  allowed
  to use
constructions and results in \ref{llo}-\ref{poi} for definable
subassignments in $\GDef_k (\LPre)$ in the rest of the paper by
referring to the corresponding ones for Presburger sets.

Note
 that $h_{\ZZ^r} =h[0,0,r]$ has more points than $\ZZ^r$ since a
point on $h[0,0,r]$ consists of a tuple $(a,K)$ with $K$ in ${\rm
Field}_k$ and $a\in\ZZ^r$. Likewise,  the ring $\cP (h[0,0,r])$ is
larger than the ring $\cP (\ZZ^r)$ defined in \ref{cpf}, since it
contains $\cP (\ZZ^r)$ and is generated as a ring by $\cP (\ZZ^r)$
and characteristic functions of definable subassignments of the
final object $h_{\Spec k}$.

\section{Constructible motivic functions}\label{sec2}From now on and until the end of the paper, we shall work with
$\LPas=\LPre$ as Denef-Pas language.

\subsection{Grothendieck rings and semirings}\label{rd}
In previous publications on motivic integration, free abelian
groups of varieties (or Chow motives) over $k$ were used to build
up Grothendieck rings. Here we shall  consider a $\LPas$-variant
using definable subassignments, which has several advantages, in
particular, of being ``universal'' in our setting. In the absolute
case, one has natural ring morphisms to the previous Grothendieck
rings (cf.~section \ref{compa}). Since in integration theory
``positive'' functions play an important role, we shall also
consider Grothendieck semirings. We first recall some basics from
the theory of semirings.

\subsubsection{Semirings}
Let us recall that a (commutative) semiring $A$ is a set equipped
with two operations: addition and multiplication. With respect to
addition $A$ is a commutative semigroup (monoid) with  $0$ as unit
element. With respect to multiplication $A$ is a commutative
semigroup  with $1$ as a unit element. Furthermore the two
structures are connected by the axioms $x (y + z) = xy + xz$ and
$0x = 0$. A morphism of semirings is a mapping compatible with the
unit elements and the operations. A module (or semimodule) over a
semiring $A$ is a commutative semigroup $M$ with an operation
$\cdot : A \times M \rightarrow A$ satisfying the familiar axioms
$(ab) \cdot m = a \cdot (b \cdot m)$, $(a + b) \cdot m = a  \cdot
m + b \cdot m$, $a \cdot (m + n) = a \cdot m + a \cdot n$, $0
\cdot m = 0$, $a \cdot 0 = 0$ and $1 \cdot m = m$. One defines
morphisms of $A$-modules in the usual way. Also, if $M$ and $N$
are $A$-modules, one can define their tensor product $M \otimes_A
N$ in the usual way using generators and relations. It is an
$A$-module representing the functor of bilinear morphisms on $M
\times N$ and its existence also follows from classical
representability results. If $B$ is an $A$-algebra (that is, a
semiring together with a morphism $A \rightarrow B$), for every
$A$-module $M$ the module $B \otimes_A M$ has a natural $B$-module
structure compatible with the $A$-module structure. Also, if $B$
and $C$ are $A$-algebras, the formula $\sum_i (b_i \otimes c_i)
\sum_j (b'_j \otimes c'_j) =\sum_{i, j} b_i b'_j \otimes c_i c'_j$
endowes $B \otimes_A C$ with a structure of $A$-algebra.

\subsubsection{Definition and properties of $\RDef$}\label{rdbis}
Let $Z$ be a definable subassignment in $\GDef_k$. We shall use in a
essential way the full subcategory $\RDef_Z$ of $\GDef_Z$, whose
objects are definable subassignments $Y$ of $Z \times h_{\AA^n_k}$,
for some $n$, the morphism $Y \rightarrow Z$ being the one induced
by projection on the $Z$ factor. If $Y$ and $Y'$ are two objects of
$\RDef_Z$, their fiber product $Y \otimes_Z Y'$ together with the
canonical morphism $Y \otimes_Z Y' \rightarrow Z$ yields an object
of  $\RDef_Z$.

It is this category $\RDef_Z$, and not $\Def_Z$, that is used to
built relative Grothendieck rings over $Z$ with. This is so because
we want the Grothendieck rings to capture information over the
residue fields, not over the valued field nor the integers.

We define the Grothendieck semigroup
$SK_0 (\RDef_{Z,k} (\LPas))$ - or $SK_0 (\RDef_Z)$ for short -,
as the quotient of
the free abelian semigroup
over symbols $[Y \rightarrow Z]$ with $Y
\rightarrow Z$ in $\RDef_Z$ by relations
\begin{equation}\label{eq0}
[\emptyset \rightarrow Z] = 0,
\end{equation}
\begin{equation}\label{eq1}
[Y \rightarrow Z] = [Y' \rightarrow Z]
\end{equation}
if $Y \rightarrow Z$ is isomorphic to $Y' \rightarrow Z$ and
\begin{equation}\label{eq2}
[(Y \cup Y') \rightarrow Z] + [(Y \cap Y') \rightarrow Z]
= [Y \rightarrow Z] + [Y' \rightarrow
Z]
\end{equation}
for $Y$ and $Y'$ definable subassignments of some $Z[0,n,0]
\rightarrow Z$. Similarly one defines the Grothendieck group $K_0
(\RDef_{Z,k} (\LPas))$, or $K_0 (\RDef_Z)$ for short, as the
quotient of the free abelian group over symbols $[Y \rightarrow Z]$
with $Y \rightarrow Z$ in $\RDef_Z$ by relations (\ref{eq1}) and
(\ref{eq2}). Cartesian fiber product over $Z$ induces a natural
semiring, resp.~ring, structure on $SK_0 (\RDef_Z )$, resp.~$K_0
(\RDef_Z )$, by setting
$$
[Y \rightarrow Z] [Y' \rightarrow Z] = [Y \otimes_Z Y' \rightarrow Z].
$$
Let us remark that $[Z\rightarrow Z]$ is the multiplicative unit
and that $K_0 (\RDef_Z )$ is nothing but the ring obtained from
$SK_0 (\RDef_Z )$ by inverting additively every element. However
note that
there is no reason for the canonical
morphism $SK_0 (\RDef_Z ) \rightarrow K_0 (\RDef_Z )$ to be
injective in general.

We extend some operations from section \ref{sec1}. If $f$ is a
morphism $ Z  \rightarrow Z'$ in $\GDef_k$, then the functor $f^*$,
as defined in \ref{fstar},
 induces a semiring morphism $f^{*}: SK_0 (\RDef_{Z'} ) \rightarrow
SK_0 (\RDef_Z )$ and a ring morphism $f^{*}: K_0 (\RDef_{Z'} )
\rightarrow K_0 (\RDef_Z )$. Also, if $z$ is a point of $Z'$, then
the functor $i_z^{*}$,
 as defined in \ref{nnn},
 induces a semiring morphism $ i_z^{*} : SK_0 (\RDef_{Z'})
\rightarrow SK_0 (\RDef_{k (z)})$ and a ring morphism $ i_z^{*} :
K_0 (\RDef_{Z'}) \rightarrow K_0 (\RDef_{k (z)})$. Note that
 $f^*\circ g^*=(g \circ f )^*$
  for composable morphisms $f,g$ in
$\GDef_k$.

\par
 Constructing a direct image functor $f_!$ (on constructible
functions) for general morphisms $f:X\to Y$ in $\Def_Y$ is one of
the main purposes of the paper. However, on the level of $\RDef_Y$
and its Grothendieck rings, with $f$ a morphism in $\RDef_Y$, $f_!$
is easy to define and is merely a universal operator. Namely, if
$f:X\to Y$ is a morphism $\RDef_Y$, the functor $f_!$, as defined in
\ref{fstar},
 restricts to a functor $f_!:\RDef_X\to \RDef_Y$ and this induces a
semiring morphism $f_!: SK_0 (\RDef_X ) \rightarrow SK_0 (\RDef_Y)$
and a ring morphism $f_!: K_0 (\RDef_X) \rightarrow K_0 (\RDef_Y)$.
Note that $g_! \circ f_! =(g \circ f)_!$ for composable such
morphisms. It is also clear that the projection formula $f_! (x f^*
(y)) = f_! (x) y$ holds for $x$ in
 $SK_0 (\RDef_{X})$ (resp.~$K_0 (\RDef_{X})$) and $y$ in $SK_0
(\RDef_{Y})$ (resp.~$K_0 (\RDef_{Y})$).

Note that any $a$ in $SK_0 (\RDef_Z)$ is of the form $[\pi : Y
\rightarrow Z]$, with $Y$ in $\RDef_Z$ and $\pi$ the projection.
The definable subassignment $\pi (Y)$ depends only on $a$ and we
denote it $\Supp (a)$. Note that, by definition, a point $z$ of
$Z$ is a point of $\Supp (a)$ if and only if $i_z^*(a)$ is
different from the empty subassignment. Also, for $a$ and $b$ in
$SK_0 (\RDef_Z)$, we have $\Supp (a + b) = \Supp (a) \cup \Supp
(b)$ and $\Supp (a b) = \Supp (a) \cap \Supp (b)$.

\subsection{Constructible Presburger functions}\label{hhhp}
In \ref{uio}, we assigned to  every $Z$ in $\GDef_k$ the ring $\cP
(Z)$ of constructible Presburger functions on $Z$. If $f : Z
\rightarrow Y$ is a morphism in $\GDef_k$, composition with $f$
yields natural morphisms $f^* : \cP (Y) \rightarrow \cP (Z)$ and
$f^* : \cP_+ (Y) \rightarrow \cP_+ (Z)$, namely, by sending
$\varphi$ to $\varphi\circ f$. Similarly, if $z$ is point of $Z$
we have morphisms $i^*_z : \cP (Z) \rightarrow \cP (h_{\Spec {k
(z)}})$ and $i^*_z : \cP_+ (Z) \rightarrow \cP_+ (h_{\Spec {k
(z)}})$.

For $Y$ a definable subassignment of $Z$, we denote by ${\bf 1}_Y$
the function in $\cP (Z)$ with value 1 on $Y$ and zero on $Z
\setminus Y$. We shall denote by $\cP^0 (Z)$ (resp.~$\cP^0_+ (Z)$)
the subring (resp.~subsemiring) of $\cP (Z)$ (resp.~$\cP_+ (Z)$)
generated by the functions ${\bf 1}_Y$ for all definable
subassignments $Y$ of $Z$ and by the constant function $\LL - 1$.

Let us denote by $\LL_Z = \LL$ the class of $Z \times h_{\AA^1_k}$ in
$K_0 (\RDef_Z)$ and in
$SK_0 (\RDef_Z)$.
We also denote
by $\LL_Z - 1 = \LL - 1$ the class of $Z \times h_{\AA^1_k \setminus \{0\}}$ in
$SK_0 (\RDef_Z)$. Note that $(\LL - 1) + 1 = \LL$ in
$SK_0 (\RDef_Z)$.

We have a canonical ring, resp.~semiring, morphism $\cP^0 (Z)
\rightarrow K_0 (\RDef_Z)$, resp.~$\cP^0_+ (Z) \rightarrow SK_0
(\RDef_Z)$, sending ${\bf 1}_Y$ to $[i:Y \rightarrow Z]$, with $i$
the inclusion, and $\LL- 1$ to $\LL - 1$.

\begin{prop}\label{speprod}Let $S$ be in $\GDef_k$.
\begin{enumerate}
\item[(1)]Let $W$ be a definable subassignment of
$h [0, n, 0]$. The canonical morphisms
$$
\cP (S) \otimes_{\cP^0 (S)} \cP^0 (S \times W)
\longrightarrow
\cP (S \times W)
$$
and
$$
\cP_+ (S) \otimes_{\cP^0_+ (S)} \cP^0_+ (S \times W)
\longrightarrow
\cP_+ (S \times W)
$$
are isomorphisms.
\item[(2)]
Let $W$
be a definable subassignment of
$h [0, 0, r]$. The canonical morphisms
$$
K_0 (\RDef_S) \otimes_{\cP^0 (S)}  \cP^0 (S \times W)
\longrightarrow K_0 (\RDef_{S \times W})
$$
and
$$
SK_0 (\RDef_S) \otimes_{\cP^0_+ (S)}  \cP^0_+ (S \times W)
\longrightarrow SK_0 (\RDef_{S \times W})
$$
are isomorphisms.
\end{enumerate}
%In either case, we use the morphisms $\pi^*:\cP^0 (S)\to  \cP^0 (S
%\times W)$ and $\pi^*:\cP^0_+ (S)\to  \cP^0_+ (S \times W)$ with
%$\pi :S\times W\to S$ the projection.
%
%%the ones which send $\LL- 1$ to $\LL - 1$ and $\11_Y$ to
%%$\11_{Y\times W}$.
\end{prop}

\begin{proof}Follows directly from Theorem \ref{pqe}.
\end{proof}

\subsection{Constructible motivic functions}\label{cmf}
Let $Z$ be a definable subassignment in $\GDef_k$.
We define the semiring
$\cC_+ (Z)$ of positive constructible motivic functions on $Z$
as
$$\cC_+ (Z) :=  SK_0 (\RDef_Z) \otimes_{\cP^0_+ (Z)} \cP_+ (Z).$$
Similarly we define the ring $\cC (Z)$ of constructible motivic
functions on $Z$ as
$$\cC (Z) :=  K_0 (\RDef_Z) \otimes_{\cP^0 (Z)} \cP (Z).$$

Let us remark that $\cC (Z)$ is nothing but the ring obtained from
$\cC_+ (Z)$ by inverting additively every element and that in
general there is no reason for the canonical morphism $\cC_+ (Z)
\rightarrow \cC (Z)$ to be injective.

For $\varphi$ and $\varphi'$ in $\cC_+ (Z)$, we shall write
$\varphi \geq \varphi'$ if
$\varphi = \varphi' + \varphi''$
for some  $\varphi''$ in $\cC_+ (Z)$.

\begin{prop}\label{peprod}Let $S$ be in $\GDef_k$.
\begin{enumerate}
\item[(1)]Let $W$ be a definable subassignment of
$h [0, n, 0]$. The canonical morphisms
$$
\cP (S) \otimes_{\cP^0 (S)} K_0 (\RDef_{S \times W})
\longrightarrow
\cC (S \times W)
$$
and
$$
\cP_+ (S) \otimes_{\cP^0_+ (S)} SK_0 (\RDef_{S \times W})
\longrightarrow
\cC_+ (S \times W)
$$
are isomorphisms.
\item[(2)]
Let $W$
be a definable subassignment of
$h [0, 0, r]$. The canonical morphisms
$$
K_0 (\RDef_S) \otimes_{\cP^0 (S)}  \cP (S \times W)
\longrightarrow  \cC (S \times W)
$$
and
$$
SK_0 (\RDef_S) \otimes_{\cP^0_+ (S)}  \cP_+ (S \times W)
\longrightarrow \cC_+ (S \times W)
$$
are isomorphisms.
\end{enumerate}
\end{prop}

\begin{proof}Direct consequence of Proposition \ref{speprod}.
\end{proof}

Note that $\cC (h_{\Spec k})$ is canonically isomorphic to $K_0
(\RDef_k) \otimes_{\ZZ [\LL]} \AA$ and that $\cC_+ (h_{\Spec k})$ is
canonically isomorphic to $SK_0 (\RDef_k) \otimes_{\NN [\LL - 1]}
\AA_+$.

\subsection{Inverse image of constructible motivic functions}\label{invim}
Let $f : Z \rightarrow Y$ be a morphism in $\GDef_k$.
Since $f^*$  as defined on $\cP (Y)$ and
$K_0 (\RDef_Y)$
is compatible with the morphism
$\cP^0 (Y) \rightarrow K_0 (\RDef_Y)$, one gets by tensor product
an inverse image
morphism $f^* : \cC (Y) \rightarrow \cC(Z)$. One defines similarly
$f^* : \cC_+ (Y) \rightarrow \cC_+(Z)$.
Clearly $f^*  \circ g^* = (g \circ f)^*$ and
${\rm id}^* = {\rm id}$.

Similarly, if $z$ is a point of $Z$, there are  natural extensions
$i_z^*: \cC (Z) \rightarrow \cC (h_{\Spec k (z)})$ and
$i_z^*
: \cC_+ (Z) \rightarrow \cC_+ (h_{\Spec k (z)})$ of the restrictions
$i_z^*$ already defined.

If $Z_1$ and $Z_2$ are disjoint definable subassignments
of some $h_W$, then
\begin{equation}\label{ff}
\cC_+ (Z_1 \cup Z_2) \simeq \cC_+ (Z_1) \oplus \cC_+ (Z_2) \quad
\text{and} \quad
\cC (Z_1 \cup Z_2) \simeq \cC (Z_1) \oplus \cC(Z_2).
\end{equation}
If $Z_1$ and $Z_2$ are in $\GDef_S$,
then we have canonical morphisms
\begin{equation}\label{gg}\cC_+  (Z_1) \otimes _{\cC_+  (S)} \cC_+ (Z_2) \rightarrow
\cC_+  (Z_1 \times_S Z_2)
\quad
\text{and} \quad
\cC (Z_1) \otimes _{\cC (S)} \cC(Z_2) \rightarrow
\cC (Z_1 \times_S Z_2).
\end{equation}

\subsection{Push-forward for inclusions}\label{inc}
Let $i : Z \hookrightarrow Z'$ be an  inclusion between two
definable subassignments of $h_W$. Composition with $i$ yields
morphisms $i_! : K_0 (\RDef_Z) \rightarrow K_0 (\RDef_{Z'})$ and
$i_! : SK_0 (\RDef_Z) \rightarrow SK_0 (\RDef_{Z'})$. Extension by
zero induces morphisms $i_! : \cP (Z) \rightarrow \cP (Z')$ and
$i_! : \cP_+(Z) \rightarrow \cP_+ (Z')$. Since they are compatible
on $\cP^0 (Z)$ and $\cP^0_+ (Z)$, we get morphisms $i_! : \cC (Z)
\rightarrow \cC (Z')$ and $i_! : \cC_+(Z) \rightarrow \cC_+ (Z')$
by tensor  product.

\subsection{Push-forward for $k$-projections}\label{kproj}
Let $S$ be a definable subassignment in $\GDef_k$ and consider the
projection $f : S [0,n,0] \rightarrow S$ on the first factor. Recall
that, by Proposition \ref{peprod}, we have a canonical isomorphism
$\cC (S [0,n,0]) \simeq K_0 (\RDef_{S [0,n,0]}) \otimes_{\cP^0 (S)}
\cP (S)$, so that we may define a ring morphism $f_! : \cC (S
[0,n,0]) \rightarrow \cC (S)$ by sending $\sum_i a_i \otimes
\varphi_i$ to $\sum_i f_! (a_i) \otimes \varphi_i$, with $a_i$ in $
K_0 (\RDef_{S [0,n,0]})$, $\varphi_i$ in $\cP (S)$, and $f_! (a_i)$
as in section \ref{rdbis}; this is clearly independent  of the
choices. We define a semiring morphism $f_! : \cC_+ (S [0,n,0])
\rightarrow \cC_+ (S)$ in the same way. Clearly these morphisms
satisfy the projection formula
\begin{equation}\label{proform}
f_! (x f^* (y)) = f_! (x) y,
\end{equation}
 for $x$ in $\cC(S [0,n,0])$, resp. $\cC_+ (S [0,n,0])$, and $y$ in
$\cC(S)$, resp. $\cC_+ (S)$.

\subsection{Rational series and integrability}\label{piint}
Let $S$ be a definable subassignment in $\GDef_k$. As in
\ref{brs}, we consider the power series ring $\cC (S) \llb T \rrb = \cC
(S) \llb T_1, \dots, T_r \rrb$ and $\cC (S) \llb T, T^{-1} \rrb = \cC (S) \llb T_1, \dots,
T_r, T_1^{-1}, \dots, T_r^{-1} \rrb$. We shall set $\cC (S)
\llb T \rrb_{\Gamma} := \cC (S) \otimes_{\cP (S)} \cP (S)
\llb T \rrb_{\Gamma}$ and $\cC (S) \llb T, T^{-1} \rrb_{\Gamma}:= \cC (S)
\otimes_{\cP (S)} \cP (S) \llb T, T^{-1} \rrb_{\Gamma}$ and view them as
$\cC (S)$-submodules of $\cC (S) \llb T \rrb $ and $\cC (S) \llb T,
T^{-1} \rrb$, respectively.

Now, for $\varphi$ in $\cC (S [0, 0, r])$ and $i$ in $\ZZ^r$,
we denote by $\varphi_i$ the restriction of $\varphi$
to $S \times \{i\}$, viewed as an element in $\cC (S)$,
and, as in \ref{brs}, we set
$$
\cM (\varphi) := \sum_{i \in \ZZ^r} \varphi_i T^i
$$
in $\cC (S) \llb T, T^{-1} \rrb$.

By (2) of Proposition \ref{peprod}, we have a canonical isomorphism
\begin{equation}\label{cais}
\cC (S [0, 0, r])
\simeq
\cC (S) \otimes_{\cP (S)} \cP (S [0, 0, r]).
\end{equation}

Since, by the extension of Theorem \ref{presrat}
to
the  definable subassignment setting,
we have an isomorphism of $\cP (S)$-modules
$$
\cM : \cP (S [0,0, r])
\longrightarrow
\cP (S) \llb T, T^{-1} \rrb_{\Gamma},
$$
we get by tensoring with $\cC (S)$ the
following general rationality statement.

\begin{theorem}\label{rateau}Let $S$ be a definable subassignment in $\GDef_k$.
The mapping $\varphi \mapsto \cM (\varphi)$
induces a ring isomorphism
$$
\cM : \cC (S [0,0, r])
\longrightarrow
\cC (S) \llb T, T^{-1} \rrb_{\Gamma}.
$$
\end{theorem}

Similarly, we
define the $\cC(S)$-modules
$${\rm I}_S \cC (S [0, 0, r]) :=
\cC (S) \otimes_{\cP (S)} {\rm I}_S \cP (S [0, 0, r])$$
and
$$\cC (S) \llb T, T^{-1} \rrb_{\Sigma}:=
\cC (S) \otimes_{\cP (S)}
\cP (S) \llb T, T^{-1} \rrb_{\Sigma}.$$
We also define
$${\rm I}_S \cC_+ (S [0, 0, r]) :=
\cC_+ (S) \otimes_{\cP_+ (S)} {\rm I}_S \cP_+ (S [0, 0, r]),$$
where ${\rm I}_S \cP_+ (S [0, 0, r])$ stands for the $\cP_+
(S)$-module $\cP_+ (S [0, 0, r]) \cap {\rm I}_S \cC (S [0, 0,
r])$. A function in  $\cC (S [0, 0, r])$ (resp.~$\cC_+ (S [0, 0,
r])$) will be called $S$-integrable if it belongs to ${\rm I}_S
\cC (S [0, 0, r])$ (resp.~${\rm I}_S \cC_+(S [0, 0, r])$).

By Theorem \ref{sintmel} and tensoring with $\cC (S)$, the
isomorphism $\cM$ induces an isomorphism of $\cC (S)$-modules
 \begin{equation}\label{guio}
\cM : {\rm I}_S \cC (S [0,0, r]) \longrightarrow \cC (S) \llb T,
T^{-1} \rrb_{\Sigma}.
\end{equation}

By tensoring
the morphism of $\cP (S)$-modules
$\mu_S : {\rm I}_S \cP (S [0,0,r])
\rightarrow
\cP (S)$ with $\cC (S)$, we get
a morphism of
$\cC (S)$-modules
$$\mu_S : {\rm I}_S \cC (S [0,0,r])
\longrightarrow
\cC(S).$$
Similarly we have a morphism of $\cC_+ (S)$-modules
$$\mu_S : {\rm I}_S \cC_+ (S [0,0,r])
\longrightarrow
\cC_+ (S).$$

\begin{prop}\label{wfubpres}Let $S$ be a definable subassignment in $\GDef_k$
and let $\varphi$ be in $\cC (S \times h_{\ZZ^r})$.
Write $r = r_1 + r_2$ and identify $h_{\ZZ^r}$ with
$h_{\ZZ^{r_1}} \times h_{\ZZ^{r_2}}$.
If $\varphi$ is $S$-integrable,
then
$\varphi$, as a function in
$\cC (S \times h_{\ZZ^{r_1}} \times h_{\ZZ^{r_2}})$,
is $S \times h_{\ZZ^{r_1}}$-integrable,
$\mu_{S \times h_{\ZZ^{r_1}}} (\varphi)$ is $S$-integrable
and
$$
\mu_S (\mu_{S \times h_{\ZZ^{r_1}}} (\varphi))
=
\mu_{S} (\varphi).
$$
The statement with $\cC$ replaced by $\cC_+$ holds also.
\end{prop}

\begin{proof}Follows directly from Lemma \ref{fubpres} (1).
\end{proof}

Let $\lambda :  S [0,0,r] \rightarrow S [0,0,s]$ be a  morphism in
$\Def_S$. Let $\varphi$ be  a function in $\cC (S [0, 0, r] )$,
resp.~$\cC_+ (S [0, 0, r] )$. Assume $\varphi = \11_Z \varphi$ with
$Z$ a definable subassignment of $S [0, 0, r] $ on which $\lambda$
is injective.  Recall that this means that the function $\lambda(K)$
is injective on $Z(K)$ for each $K$ in ${\rm Field}_{k}$,
cf.~section \ref{sec:sub}.
 Thus $\lambda$ restricts to an isomorphism $\lambda'$ between $Z$
and $Z':= \lambda (Z)$. We set $\lambda_+ (\11_Z \varphi) := [i'_!
(\lambda'{}^{-1})^* i^* ](\varphi)$ in $\cC (S [0, 0, r] )$,
resp.~in $\cC_+ (S [0, 0, r] )$, where $i $ and $i'$ denote
respectively the inclusions of $Z$ and $Z'$ in $S [0,0,r]$ and $S
[0,0,s]$. Clearly this definition does not depend on the choice of
$Z$
\footnote{The $+$ in $\lambda_+$ does not refer to positivity but
denotes an alternative direct image, which shall only be later on,
Proposition \ref{same}, shown to be related with $!$ of $\lambda_!$.
Also it should not be confused with the functor $\lambda_+$ in ${\rm
A0(b) } $ of Theorem \ref{mt}.}.

The following statement follows directly from Lemma \ref{presinj}:

\begin{prop}\label{inj}
Let $\lambda :  S [0,0,r] \rightarrow S [0,0,s]$ be a morphism in
$\GDef_S$. Assume $\lambda$ is injective on a definable
subassignment $Z$
 of $S [0,0,r]$.
 Let $\varphi$ be a function in  $\cC (S [0,0,r])$ such that $\11_Z
\varphi = \varphi$. Then $\varphi$ is $S$-integrable if and only if
$\lambda_+(\varphi)$ is $S$-integrable. If these conditions are
satisfied then $\mu_S (\varphi) = \mu_S  (\lambda_+(\varphi))$. The
statement with $\cC$ replaced by $\cC_+$ holds also. \qed
\end{prop}

\subsection{Positivity and Fubini}
Let $S$ be a definable subassignment in $\GDef_k$.
It is quite clear that if $f$ and $g$ are
in
$\cP_+ (S [0,0, r])$, $f \geq g$ and
$f$ is $S$-integrable, then $g$ is $S$-integrable.
We shall now prove a  similar  statement for $\cC_+$.

For $a$ in $SK_0 (\RDef_S)$, we shall write $\11_a := \11_{\Supp
(a)}$,  with $\Supp (a)$ as defined in \ref{rdbis}.

\begin{prop}\label{recip}
Let $S$ be a definable subassignment in $\GDef_k$ and
let $f$ be a function in $\cC_+ (S [0, 0, r])$.
Write $f = \sum_i a_i \otimes \varphi_i$,
with $a_i$
in $SK_0 (\RDef_S)$ and $\varphi_i$
in $\cP_+ (S [0, 0, r])$.
Then $f$ is $S$-integrable if and only every function
$\11_{a_i} \varphi_i$ is $S$-integrable.
\end{prop}

\begin{proof}Let $f$ be a $S$-integrable function in
$\cC_+ (S [0, 0, r])$. Write $f = \sum_i a_i \otimes \varphi_i$,
with $a_i$ in $SK_0 (\RDef_S)$ and $\varphi_i$ in $\cP_+ (S [0, 0,
r])$. Since  $f$ is in ${\rm I}_S \cC_+ (S [0, 0, r])$, we may
also write $f = \sum_j b_j  \otimes \psi_j$, with $b_j$ in $SK_0
(\RDef_S)$ and $\psi_j$ in ${\rm I}_S \cP_+ (S [0, 0, r])$. We now
use the degree function $\dl$ defined in \ref{Ldeg}. Recall that,
for two functions $\varphi$ and $\varphi'$ in $\cP_+ (S [0, 0,
r])$, $\dl (\varphi + \varphi') = \sup (\dl (\varphi), \dl
(\varphi'))$. Let us also remark that if $a$ belongs to $\cP^0_+
(S)$ and $\varphi$ to $\cP_+ (S [0, 0, r])$, the difference $\dl
(a \varphi) - \dl (\11_a \varphi) $ may take only a finite number
of distinct values, uniformly for $s$ in $|S|$. It then follows from
the relations defining the tensor product $SK_0 (\RDef_S)
\otimes_{\cP^0_+ (S)} \cP_+ (S [0, 0, r])$ that there is a
constant $C$ such that $\dl (\11_{a_i} \varphi_i) \leq C + {\rm
sup}_j (\dl (\psi_j))$. From Proposition \ref{ldegint} we deduce
that every function  $\11_{a_i} \varphi_i$ is $S$-integrable. The
reverse implication being clear, this concludes the proof.
\end{proof}

\begin{cor}\label{repo}Consider a morphism
$g : S \rightarrow \Lambda$ in $\Def_k$. For every point $\lambda$
in $\Lambda$ consider the fiber $S_{\lambda}$ of $g$ at $\lambda$. A
function $f$ in $\cC_+ (S [0, 0, r])$ is $S$-integrable if and only
if, for every point $\lambda$ in $\Lambda$, the restriction
 $f_{\lambda}\in \cC_+ (S_{\lambda} [0, 0, r])$ of $f$ to
$S_{\lambda}[0,0,r]$ is $S$-integrable. Furthermore, if these
conditions are satisfied, $\mu_{S_{\lambda}} (f_{\lambda})$ is equal
to the restriction of $\mu_{S} (f)$ to $S_{\lambda}$.
\end{cor}

\begin{proof}The analogous
result
with $\cC_+$ replaced by $\cP_+$ being clear, the statement follows directly
from Proposition \ref{recip}.
\end{proof}

\begin{prop}\label{geq}Let $S$ be a definable subassignment in $\GDef_k$.
Let $f$ and $g$ be functions in
$\cC_+ (S [0, 0, r])$. If $f \geq g$ and $f$ is
$S$-integrable,
then $g$ is $S$-integrable.
\end{prop}

\begin{proof}We may write $f = g + h$ with $h$ in $\cC_+ (S [0, 0, r])$.
We write
$g = \sum_j b_j  \otimes \psi_j$,
with $b_j$
in $SK_0 (\RDef_S)$ and $\psi_j$ in $\cP_+ (S [0, 0, r])$,
and similarly $h = \sum_{j'} b_{j'}  \otimes \psi_{j'}$.
Since $f =  \sum_j b_j  \otimes \psi_j + \sum_{j'} b_{j'}  \otimes \psi_{j'}$,
it follows from Proposition \ref{recip}
that every function
$\11_{b_j} \psi_j$ is $S$-integrable, which concludes the proof.
\end{proof}

We now can state the analogue of Lemma \ref{fubpres} (2).
 Recall
that $h_{\ZZ^r}$ stands for $h[0,0,r]$.
 \begin{prop}\label{strongpresfub}
Let $S$ be a definable subassignment in $\GDef_k$ and let
$\varphi$ be in $\cC_+ (S \times h_{\ZZ^r})$. Write $r = r_1 +
r_2$ and identify $h_{\ZZ^r}$ with $h_{\ZZ^{r_1}} \times
h_{\ZZ^{r_2}}$. The function  $\varphi$ is $S$-integrable if and
only if, as a function in $\cC_+ (S \times h_{\ZZ^{r_1}} \times
h_{\ZZ^{r_2}})$, it is $S \times h_{\ZZ^{r_1}}$-integrable and
$\mu_{S \times h_{\ZZ^{r_1}}} (\varphi)$ is $S$-integrable.
\end{prop}

\begin{proof}If $\varphi$ is
$S \times h_{\ZZ^{r_1}}$-integrable, we may write $\varphi =
\sum_i a_i \otimes \varphi_i$ with $a_i$ in $SK_0 (\RDef_S)$ and
$\varphi_i$ in ${\rm I}_{S \times h_{\ZZ^{r_1}}}\cP_+ (S \times
h_{\ZZ^r})$. Replacing $\varphi_i$ by $\11_{a_i} \varphi_i$ we may
even assume $\varphi_i = \11_{a_i} \varphi_i$. Hence $\mu_{S
\times h_{\ZZ^{r_1}}} (\varphi) = \sum_i a_i \otimes \mu_{S \times
h_{\ZZ^{r_1}}} (\varphi_i)$. If $\mu_{S \times h_{\ZZ^{r_1}}}
(\varphi)$ is $S$-integrable, it follows from Proposition
\ref{recip} that the functions $\11_{a_i} \mu_{S \times
h_{\ZZ^{r_1}}} (\varphi_i) = \mu_{S \times h_{\ZZ^{r_1}}}
(\varphi_i)$ are all $S$-integrable. One then deduces from Lemma
\ref{fubpres} (2) that the functions $\varphi_i$ are all $
S$-integrable, hence $\varphi$ is $S$-integrable. The reverse
implication is already known (Proposition \ref{wfubpres}).
\end{proof}

\begin{prop}\label{mixing}Let $S$ be a definable subassignment in $\GDef_k$
and consider the projections $\pi^1: S [0, n, r] \rightarrow S [0,
0, r]$ and $\pi^2: S [0, n, 0] \rightarrow S [0, 0, 0]$. Let
$\varphi$ be a function in $\cC_+ (S [0, n, r])$. Then $\varphi$ is
$S [0, n, 0]$-integrable if and only if the function $\pi^1_!
(\varphi)$ in $\cC_+ (S [0, 0, r])$ is $S$-integrable. If these
conditions hold, then
$$
\pi^2_! (\mu_{S [0, n, 0]} (\varphi)) = \mu_S (\pi^1_! (\varphi)).
$$
\end{prop}

\begin{proof}
Let $\varphi$ be a function in $\cC_+ (S [0, n, r])$. We may write
$\varphi = \sum_i a_i \otimes \varphi_i$ with $a_i$ in $SK_0
(\RDef_{S [0,n,0]})$ and $\varphi_i$ in $\cP_+ (S [0,0,r])$. Indeed,
it follows from the second part of Theorem \ref{pqe} that the
canonical morphism
$$\cP^0_+ (S [0,n,0]) \otimes_{\cP^0_+ (S)} \cP^0_+ (S [0,0,r])
\longrightarrow \cP^0_+ (S [0,n,r])$$ is an isomorphism, from
which we deduce a canonical isomorphism
$$\cC_+ (S [0,n,r]) \simeq
SK_0 (\RDef_{S [0,n,0]}) \otimes_{\cP^0_+ (S)} \cP_+ (S [0,0,r])$$
using Propositions \ref{speprod} and \ref{peprod}. We have $\pi^1_!
(\varphi) = \sum_i \pi^2_! (a_i) \otimes \varphi_i$. The key remark
is that  $\11_{a_i} \pi^1{}^* (\varphi_i)$ is $S [0,n,0]$-integrable
if and only if $\11_{\pi^2_! a_i} \varphi_i$ is $S$-integrable,
which holds since integrability is defined by a pointwise condition.
Hence, it follows from Proposition \ref{recip} that $\varphi$ is $S
[0,n,0]$-integrable if and only if $\pi^1_! (\varphi)$ is
$S$-integrable. Let us assume that these conditions hold, so that,
by Proposition \ref{recip}, each $\11_{a_i} \varphi_i$ is
$S$-integrable. We assume $\varphi_i = \11_{a_i} \varphi_i$ for
every $i$. Since $\mu_{S [0,n,0]} (\pi^1{}^* (\varphi_i)) =
\pi^2{}^* (\mu_S (\varphi_i))$, we get $\mu_{S [0,n,0]} (\varphi) =
\sum_i a_i \otimes \mu_S (\varphi_i)$, and hence, we deduce $\pi^2_!
\mu_{S [0,n,0]} (\varphi) = \sum_i \pi^2_! (a_i) \otimes \mu_S
(\varphi_i) = \mu_S (\pi^1_! (\varphi))$.
\end{proof}

Let $\lambda :  S [0,n,r] \rightarrow S [0,n',r']$ be a  morphism
in $\Def_S$. Let $\varphi$ be a function in $\cC (S [0, n, r] )$,
resp.~in $\cC_+ (S [0, n, r] )$. Assume $\varphi = \11_Z \varphi$
with $Z$ a definable subassignment of $S [0, n, r]$ on which
$\lambda$ is injective. Thus $\lambda$ restricts to an isomorphism
$\lambda'$ between $Z$ and $Z':= \lambda (Z)$. We define
$\lambda_+(\varphi)$ in $\cC(S [0,n',r'])$, resp.~in $\cC_+ (S [0,
n', r'] )$ as $[i'_! (\lambda'{}^{-1})^* i^* ](\varphi)$, where $i
$ and $i'$ denote respectively the inclusions of $Z$ and $Z'$ in
$S [0,n,r]$ and $S [0,n',r']$. Clearly this definition does not
depend on the choice of $Z$ and is  compatible with the definition
of $\lambda_+$ in section \ref{piint} when $n= n'= 0$.

\begin{prop}\label{stronginj}
Let $\lambda : S [0, n, r] \rightarrow S [0, n' , r']$
be a morphism in $\Def_S$. Let $\varphi$ be a
function in $\cC_+ (S [0, n, r])$ such that  $\varphi = \11_Z \varphi$
with
$Z$ a definable subassignment of
$S [0, n, r]$ on which $\lambda$ is injective.
Then $\varphi$ is
$S[0,n,0]$-integrable if and only if $\lambda_+(\varphi)$ is
$S[0,n',0]$-integrable and if this is the case then
\begin{equation}\label{opu}
p_! (\mu_{S[0,n,0]}(\varphi)) = p'_! (\mu_{S[0,n',0]}( \lambda_+
(\varphi))),
\end{equation}
with $p:S[0,n,0]\to S$ and $p':S[0,n',0]\to S$ the projections.
\end{prop}

\begin{proof}In the case where
$n = n'$ and $\lambda$ is the identity on the $\AA^n_k$-factor,
that is, $\lambda$ is of the form $(s, \xi, \alpha) \mapsto (s,
\xi, g (s, \alpha))$, with $s$ variable on the $S$-factor, $\xi$
on the $\AA^n_k$-factor and $\alpha$ on the $\ZZ^r$-factor, the
statement follows directly from Proposition \ref{inj}. Assume now
that $r = r'$ and that $\lambda$ is the identity on the
$\ZZ^r$-factor, that is, $\lambda$ is of the form $(s, \xi,
\alpha) \mapsto (s, f (s, \xi), \alpha)$. Write $g:S[0,n,0]\to
S[0,n',0]:(s, \xi) \mapsto (s, f (s, \xi))$. Let $\varphi$ be a
function in $\cC_+ (S [0, n, r])$ such that $\varphi = \11_Z
\varphi$ with $Z$ a definable subassignment of $S [0, n, r]$ on
which $\lambda$ is injective. As in the proof of Proposition
\ref{mixing} we may write
\begin{equation}\label{huy}
\varphi = \sum_i a_i \otimes \varphi_i
\end{equation}
with $a_i$ in $SK_0 (\RDef_{S [0,n,0]})$, $\varphi_i$ in $\cP_+ (S
[0,0,r])$. By the second part of Theorem \ref{pqe}, we may assume
that $g$ is injective on $\Supp(a_i)$ for each $i$. Using the
injectivity of $g$ on $\Supp(a_i)$, we can write
\begin{equation}\label{huybis}
\lambda_+ (\varphi) = \sum_i g_!(a_i) \otimes \varphi_i,
\end{equation}
where $g_!(a_i)$ in $SK_0(\RDef_{S[0,n',0]})$ is defined as
$[X_i\otimes_{S[0,n',0]}S[0,n',0]\to S[0,n',0]]$ with $a_i=[X_i\to
S[0,n,0]]$\footnote{This definition of $g_!$ is consistent with
the definition of the push-forward in section \ref{kproj}. Note
that $a_i$ can always be written in the form $[X_i\to S[0,n,0]]$.}.
If $\varphi$ is $S [0, n, 0]$-integrable, we may assume all the
functions $\varphi_i$ are $S$-integrable, hence $\lambda_+
(\varphi)$ is $S [0, n', 0]$-integrable. For the reverse
implication, note that $\varphi = {\tilde \lambda}_+ (\lambda_+
(\varphi))$ for any morphism $\tilde \lambda : S [0, n', r']
\rightarrow S [0, n, r]$ which restricts to the inverse of
$\lambda'$ on $Z'$. Relation (\ref{opu}) then follows from
(\ref{huy}), (\ref{huybis}), and from the obvious relation
$p_!=p'_!\circ g_!$.

Note that if the statement of Proposition \ref{stronginj} holds
for two composable morphisms $\lambda_1$ and $\lambda_2$, it still holds
for $\lambda_2 \circ \lambda_1$. In particular it follows from the previous
discussion that the statement we want to prove holds for
$\lambda$ of the form
$(s, \xi, \alpha) \mapsto (s, f (s, \xi), g (s, \alpha))$.
Now consider the case of a general
morphism $\lambda : S [0, n, r] \rightarrow S [0, n' , r']$.
Let $\varphi$ be a
function in $\cC_+ (S [0, n, r])$ such that  $\varphi = \11_Z \varphi$
with
$Z$ a definable subassignment of
$S [0, n, r]$ on which $\lambda$ is injective.
By the
second part of Theorem \ref{pqe}
there is a finite partition of $Z$ into definable subassignments
$Z_i$, such that the restriction
of
$\lambda$ to each $Z_i$ is of the form
$(s, \xi, \alpha) \mapsto (s, f_i (s, \xi), g_i (s, \alpha))$.
Since $\lambda_+ (\11_{Z_i} \varphi)$ only depends on the restriction of
$\lambda$ to $Z_i$, it follows that the statement we want to prove
holds
for $\11_{Z_i} \varphi$, hence also for $\varphi = \sum_i \11_{Z_i} \varphi$.
\end{proof}

\section{Constructible motivic Functions}\label{sec4}

\subsection{Dimension and relative dimension}
Let $Z$ be in $\Def_k$ and let $\varphi$ be in $\cC(Z)$,
resp.~$\cC_+(Z)$. We say $\varphi$ is of $K$-dimension $ \leq d$
if it may be written as a finite sum $\varphi = \sum \lambda_i
\11_{Z_i}$ in $\cC (Z)$, resp.~in $\cC_+ (Z)$, with ${\Kdim} \,
Z_i \leq d$. We say   $\varphi$ is of $K$-dimension $d$ if it is
of $K$-dimension $ \leq d$ and not of $K$-dimension $ \leq d - 1$.

More generally, if $Z\to S$ is in $\Def_S$ for some $S$ and
$\varphi$ in $\cC(Z)$, resp.~in $\cC_+(Z)$, we say that $\varphi$
is of dimension $\leq d$ rel.~the projection $Z\to S$ if it may be
written as a finite sum $\varphi = \sum \lambda_i \11_{Z_i}$ in
$\cC (Z)$, resp.~in $\cC_+ (Z)$, with $Z_i$ of relative dimension
$\leq d$ rel.~the projection $Z\to S$ (as in section \ref{reld}).
We also use the notion of dimension $d$ rel.~the projection $Z \to
S$, for $\varphi$ in $\cC (Z)$ or in $\cC_+ (Z)$ if it is of
relative dimension $\leq d$ but not of relative dimension $\leq
d-1$.

\subsection{Constructible motivic Functions}\label{cmc} Let $Z$ be
a definable subassignment in $\GDef_k$. We denote by $\cC^{\leq d}
(Z)$, resp.~$\cC_+^{\leq d} (Z)$, the subgroup,
resp.~subsemigroup, of elements of $\cC (Z)$, resp.~$\cC_+ (Z)$,
of $K$-dimension $\leq d$. We denote by $C^{d} (Z)$ the quotient
$$C^{d} (Z) : = \cC^{\leq d} (Z) / \cC^{\leq d - 1} (Z)$$
and we set
$$
C (Z) =  \bigoplus_{d \geq 0} C^{d} (Z),
$$
which is a graded abelian group. Similarly, we denote by $C_+^{d}
(Z)$ the quotient
$$C^{d}_+ (Z) : = \cC_+^{\leq d} (Z) / \cC_+^{\leq d - 1} (Z)$$
and we consider the graded abelian semigroup
$$
C_+ (Z) =  \bigoplus_{d \geq 0} C^{d}_+ (Z).
$$

An element in $C (Z)$, resp.~in $C_+ (Z)$, will be called a
constructible motivic Function, resp.~a positive constructible
motivic Function
(note the capital $F$).
It is an equivalence class of constructible motivic functions,
resp.~of positive constructible motivic functions.

If $\varphi$ in $\cC(Z)$ or in $\cC_+(Z)$ is of $K$-dimension $d$,
or if $\varphi=0$ in  $\cC(Z)$,
we denote by $[\varphi]$ the class of $\varphi$ in $C^d (Z)$ and
in $C^d_+(Z)$.

Let us remark that, since $ \cC^{\leq d} (Z)$ is an ideal in
$ \cC (Z)$, the product on $ \cC (Z)$
induces a $\cC (Z)$-module structure on $C (Z)$
and on each $C^d(Z)$.
Similarly the product on $ \cC_+ (Z)$ induces a $\cC_+ (Z)$-module
structure on $C_+ (Z)$.

\subsection{}\label{gcmc}More generally, let us fix a definable
subassignment $S$ in $\GDef_k$ and consider the category
$\GDef_S$. For $Z$ in $\GDef_S$, we define $\cC^{\leq d} (Z
\rightarrow S)$ and $\cC_+^{\leq d} (Z \rightarrow S)$ as in
\ref{cmc}, but by replacing $K$-dimension by relative
$K$-dimension relative to $S$.

We set
$C^d (Z \rightarrow S) := \cC^{\leq d} (Z\rightarrow S ) /
\cC^{\leq d - 1} (Z\rightarrow S)$ and
$$
C (Z\rightarrow S) =  \bigoplus_{d \geq 0} C^{d} (Z\rightarrow S).
$$
One defines similarly $C^d_+ (Z \rightarrow S)$ and $C_+
(Z\rightarrow S)$. Also, if $\varphi$ is of relative dimension
$d$,
or if $\varphi=0$ in $\cC(Z\rightarrow S)$,
we denote by $[\varphi]$ the class of $\varphi$ in $C^d (Z
\rightarrow S)$ and in $C^d_+ (Z \rightarrow S)$. Let us remark
that $C ({\rm id} : S \rightarrow S) = \cC (S)$ and that $C (Z
\rightarrow h_{\Spec k}) = C (Z)$ and similarly for $C_+$.

Let $f$ be a Function in $C (Z\rightarrow S)$.
For every point $s$ in $S$, $f$ naturally restricts
to a Function $f_s$ in $C (Z_s)$, where
$Z_s$ denotes the fiber of $Z$ at $s$.

If $Z_1$ and $Z_2$ are disjoint definable subassignment of
some
definable subassignment in $\GDef_S$,
then (\ref{ff}) induces isomorphisms
$$
C (Z_1 \cup Z_2
\rightarrow S) \simeq C (Z_1 \rightarrow S) \oplus C(Z_2 \rightarrow S)
$$
and
$$
C_+ (Z_1 \cup Z_2
\rightarrow S) \simeq C_+ (Z_1 \rightarrow S) \oplus C_+(Z_2 \rightarrow S).
$$
Also,  (\ref{gg}) induces morphisms
$$C (Z_1 \rightarrow S) \otimes_{\cC (S)} C (Z_2 \rightarrow S)
\rightarrow
C (Z_1 \times_S Z_2 \rightarrow S)$$
and
$$C_+ (Z_1 \rightarrow S) \otimes_{\cC_+ (S)} C_+ (Z_2 \rightarrow S)
\rightarrow
C_+ (Z_1 \times_S Z_2 \rightarrow S).$$

\section{Cell decomposition}\label{sec5}
In this section we shall state some variants and mild
generalizations of the Cell Decomposition Theorem of \cite{Pas} in
a form suitable for our needs. Our terminology concerning cells
differs slightly from that used in \cite{Pas}.

\subsection{Cells}\label{ce}
Let $S$ be in $\Def_k$. We will define the notion of a cell
$Z\subset S[1,0,0]$, at first with a base $C\subset S$, and secondly
with a more general base $C\subset S[0,r,s]$.

First let $C$ be a definable subassigment of $S$. Let $\alpha$,
$\xi$, and $c$ be definable morphisms $\alpha : C \rightarrow \ZZ$,
$\xi : C \rightarrow h_{\GG_{m, k}}$, and $c : C \rightarrow
h[1,0,0]$. The $1$-cell $Z_{C, \alpha, \xi, c}$ with basis $C$,
order $\alpha$, center $c$, and angular component $\xi$ is the
definable subassignment of $S[1,0,0]$ defined by $y$ in $C$, $\ord
(z - c (y)) = \alpha (y)$, and $\ac (z - c(y)) = \xi (y)$, where $y$
lies in $S$ and $z$ in $h[1,0,0]$. Similarly, if $c$ is a definable
morphism $c : C \rightarrow h[1,0,0]$, we define the $0$-cell $Z_{C,
c}$ with center $c$ and basis $C$ as the definable subassignment of
$S[1,0,0]$ defined by $y \in C$ and $z = c(y)$.

Secondly and finally, a definable subassignment $Z$ of $S[1,0,0]$
will be called a 1-cell, resp.~a 0-cell, if there exists a definable
isomorphism
$$\lambda : Z
\rightarrow Z_C = Z_{C, \alpha, \xi, c} \subset S[1,s,r],
$$
resp.~a definable isomorphism
$$\lambda : Z
\rightarrow Z_C = Z_{C, c} \subset S [1,s,0],
$$
for some $s,r\geq 0$, some basis $C\subset S[0,s,r]$,
resp.~$C\subset S[0,s,0]$,
 and some $1$-cell $Z_{C, \alpha, \xi, c}$,
resp.~$0$-cell $Z_{C, c}$, such that the morphism $\pi\circ\lambda$,
with $\pi$ the projection on the $S[1,0,0]$-factor, is the identity
on $Z$.

We shall call the data $(\lambda, Z_{C, \alpha, \xi, c})$,
resp.~$(\lambda, Z_{C, c})$, sometimes written for short
$(\lambda, Z_C)$, a presentation of the cell $Z$.

One should note that $\lambda^{*}$ induces a canonical bijection
between $\cC (Z_C)$ and $\cC (Z)$.

\begin{remark}\label{rem:pcell}
Cells as defined in \cite{Pas} fall within our definition, but not
the other way around. Also, any presentation as in \cite{Pas} of a
cell as in \cite{Pas} can be adapted to a presentation in the
above sense; for a $0$-cell this is trivial, and for a 1-cell this
can be done by adding one more ${\rm Ord}$-variable.
\end{remark}

\subsection{Cell decomposition}\label{cd}The following
variant of the Denef-Pas Cell Decomposition Theorem \cite{Pas}
will play a fundamental role in the present paper:

\begin{theorem}\label{np}
Let $X$ be a definable subassignment  of $S[1,0,0]$ with $S$ in
$\Def_k$.
\begin{enumerate}
 \item[(1)]The subassigment $X$ is a finite disjoint union of cells.
 \item[(2)] For every $\varphi$ in $\cC (X)$ there exists a finite partition of $X$ into
cells $Z_i$ with presentation $(\lambda_i, Z_{C_i})$, such that
$\varphi_{| Z_i} = \lambda_i^{*} p_i^{*} (\psi_i)$, with $\psi_i$
in $\cC (C_i)$ and $p_i : Z_{C_i} \rightarrow C_i$ the projection.
Similar statements hold for $\varphi$ in $\cC_+ (X)$, in $\cP(X)$,
in $\cP_+(X)$, in $K_0(\RDef_X)$, and in $SK_0(\RDef_X)$.
\end{enumerate}
\end{theorem}

We shall call a finite partition of $X$ into cells $Z_i$ as in
Theorem \ref{np}(1), resp.~\ref{np}(2) for a function $\varphi$, a
cell decomposition of $X$, resp.~a cell decomposition of $X$
adapted to $\varphi$.

\begin{proof}[Proof of Theorem \ref{np}]
Clearly (2) implies (1). We show how (2) follows from the cell
decomposition Theorem.~3.2 of \cite{Pas}. To fix notation, let $S$
be a definable subassignment of $h[m,n,r]$.

First we let $\varphi$ be in $\cC (X)$. Write $\varphi$ as $\sum_i
a_i\otimes \varphi_i$ with $a_i$ in $ K_0(\RDef_X)$ and $\varphi_i$
in $\cP(X)$. Let $f_1,\ldots,f_t$ be all the polynomials in the
${\rm Val}$-variables occurring in the formulas\footnote{By this we
mean that we take the defining formulas for $X$ and, for each $i$,
of subassignments in $\RDef_X$ representing $a_i$ and the defining
formulas for all definable morphisms $\alpha:X\to\ZZ$ occurring in
$\varphi_i$, written as a sum of products of constants in $\AA$ and
functions of the forms $\LL^\alpha$ and $\alpha$.} describing the
data $X$, $a_i$, and $\varphi_i$, where we may suppose that these
formulas do not contain quantifiers over the valued field sort.
Apply the cell decomposition Theorem 3.2 of \cite{Pas} to the
polynomials $f_i$. Using Remark \ref{rem:pcell}, we see that this
yields a partition of $h[m+1,0,0]$ into cells $\tilde Z_i$ with
presentations $\lambda_i:\tilde Z_i\to \tilde Z_{C_i}$ and with some
center $c_i$. Write $x=(x_1,\ldots,x_{m+1})$ for the ${\rm
Val}$-variables, $\xi=(\xi_i)$ for the ${\rm Res}$-variables and
$z=(z_i)$ for the ${\rm Ord}$-variables on $\tilde Z_{C_i}$.

If $\tilde Z_i$ is a $1$-cell, we may suppose that for $(x,\xi,z)$
in $\tilde Z_{C_i}$ we have $\ord(x_{m+1}-c_i)=z_1$ and
$\ac(x_{m+1}-c_i)=\xi_1$, by changing the presentation of $\tilde
Z_i$ if necessary (that is, by adding more ${\rm Ord}$-variables and
${\rm Res}$-variables). By the application of Theorem 3.2 of
\cite{Pas} and by changing the presentation as before if necessary,
we may also assume that
\[
\begin{array}{l}
\ord\, f_j(x)= z_{k_j},\\
\ac\, f_j(x) =  \xi_{l_j},
\end{array}\]
for $(x,\xi,z)$ in a $1$-cell $\tilde Z_{C_i}$, where the indices
$k_j$ and $l_j$ only depend on $j$ and $i$.

Since the condition $f(x)=0$ is equivalent to $\ac(f(x))=0$, we may
suppose that, in the formulas describing $X$, $a_i$, and
$\varphi_i$, the only terms involving ${\rm Val}$-variables are of
the forms $\ord\, f_j(x)$ and $\ac\, f_j(x)$. Combining this with
the above description of $\ord\, f_j(x)$ and $\ac\, f_j(x)$ one can
then easily construct a partition of $X$ into cells $Z_i$ and for
each such cell a constructible functions $\psi$ which satisfies the
requirements of the theorem. If $\varphi$ is in $\cC_+ (X)$,
resp.~in $\cP(X)$, $\cP_+(X)$, $K_0(\RDef_X)$, or in
$SK_0(\RDef_X)$, the same argument works.
\end{proof}

\begin{example}\label{trivexa}A cell decomposition for
$h_{\AA^1_{k \llp t \rrp}}$:  take  the disjoint union of the 0-cell
$\{0\}$ and the 1-cell $h_{\AA^1_{k \llp t \rrp}} \setminus \{0\}$ with
presentation $(\lambda, Z')$ with $\lambda (z) = (z, \ac (z), \ord
(z))$ and $Z'$ defined by $\ac (z) = \xi, \ord (z) = i,\
\xi\not=0$.
\end{example}

\begin{cor}\label{finitechoice}
Let $Y$ and $Z\subset Y[m',0,0]$ be definable subassignments and let
$f:Z\to Y$ be the projection. Suppose that for each $y=(y_0,K)$ in
$Y$, the set $Z_y(K)$ is finite, with $Z_y$ the fiber above $y$.
Then the cardinality of $Z_y(K)$ is bounded uniformly in
$y=(y_0,K)$, and, there exists a definable isomorphism $g:Z\to
Y'\subset Y[0,n,r]$ over $Y$ such that $p\circ g=f$, with
$p:Y[0,n,r]\to Y$ the projection.
\end{cor}
\begin{proof}
Replacing $Z$ by the graph of $f$, the essential case to prove is
when $Z\subset h[m+m',n,r]$, $Y=h[m,n,r]$, and  $f$ is the
projection. First suppose $m'=1$. Applying Theorem \ref{np} to
$Z$, the proposition follows immediately piecewise, and hence
globally, since a finite partition can be replaced by one part by
allowing for extra parameters. The case of $m'>1$ is treated by an
inductive application of Theorem \ref{np}.
\end{proof}

\begin{cor}\label{inj-constant}
Let $S$ be in $\Def_k$ and $f:S[1,0,0]\to h[1,0,0]$ be a definable
morphism. Then there exists a definable isomorphism $g:S[1,0,0]\to
Z\subset S[1,n,r]$ over $S[1,0,0]$ and a partition of $Z$ into
parts $X$, $Y$ such that for each $s=(s_0,K)$ in $S[0,n,r]$ the map
$X_{s}(K)\to K\llp t \rrp:u\mapsto f\circ g^{-1}(s,u)$ is injective and
the map $Y_{s}(K)\to K\llp t \rrp:u\mapsto f\circ g^{-1}(s,u)$ is
constant, with $X_{s}$ and $Y_{s}$ the fiber of $X$, resp.~of $Y$,
above $s$ under the projection $X\to S[0,n,r]$, resp.~$Y\to
S[0,n,r]$.
\end{cor}
\begin{proof}
It follows from Theorem \ref{prop:man}, the definition of the
dimension of definable subassignments, and the Implicit Function
Theorem for $K\llp t \rrp$-analytic maps, that there is a partition
$\{X_0,Y_0\}$ of $S[1,0,0]$ such that for each $s=(s_0,K)$ in $S$
the map $X_{0s}(K)\to K\llp t \rrp:u\mapsto f(s,u)$ is finite to one and
the map $Y_{0s}(K)\to K\llp t \rrp:u\mapsto f(s,u)$ has finite image,
with $X_{0s}$ and $Y_{0s}$ the fibers. To obtain $X$, apply
Corollary \ref{finitechoice} to the map $X_0\to
S[1,0,0]:(s,u)\mapsto (s,f(s,u))$. To obtain $Y$, apply Corollary
\ref{finitechoice} to the projection map $\pi:p({\rm
Graph}(f_{|Y_0})) \to S$ where $p$ is the projection $S[2,0,0]\to
S[1,0,0]$ onto the $(s,f)$ variables.
\end{proof}

\begin{lem}\label{rem:presen}
Let $\varphi$ be in $\cC_+(S[1,0,0])$ for some $S$ in $\Def_k$ and
let $Z\subset S[1,0,0]$ be a $1$-cell which is adapted to
$\varphi|_Z$ and which has a presentation $\lambda:Z\to Z_C=Z_{C,
\alpha, \xi, c}$. There is some $\psi$ in $\cC_+(C)$ with
$\lambda^*p^*(\psi)=\varphi|_Z$, where $p:Z_C\to C$ is the
projection. Let $C_G$ be any nonempty definable subassignment of
$C$ and let $\alpha_G,\xi_G,c_G$ be the restrictions of $\alpha,
\xi, c$ to $C_G$. Then, the subassignment $Z_G$ of $Z$ given by
$p\circ \lambda(x)\in C_G$ is a $1$-cell with presentation
$\lambda_G:Z_G\to Z_{C_G, \alpha_G, \xi_G, c_G}$, where
$\lambda_G$ is the restriction of $\lambda$ to $Z_G$. Moreover,
$Z_G$ is adapted to $\varphi|_{Z_G}$ and
$\varphi|_{Z_G}=\lambda_G^*p_G^*(\psi_G)$, where $\psi_G$ is the
restriction of $\psi$ to $C_G$ and $p_G$ the restriction of $p$ to
$Z_G$. A similar statement holds for $0$-cells.
\end{lem}
\begin{proof}
Clear.
\end{proof}

\subsection{Refinements}
Let $\cP = (Z_i)_{i \in I}$ and $\cP' = (Z'_j)_{j \in J}$ be two
cell decompositions of a definable subassignment $X\subset
S[1,0,0]$ for some $S$ in $\Def_k$. We say $\cP$ is a refinement
of $\cP'$ and write $\cP \prec \cP'$ if for every $j$ in $J$ there
exists $i (j)$ in $I$ such that $Z_{i (j)} \subset Z'_j$.

\begin{lem}\label{refin:g}
Let $Z$ and $Z'$ be cells in $S[1,0,0]$ for some $S$ in $\Def_k$.
Let $(\lambda', Z'_{C'})$ be a presentation of $Z'$. If $Z \subset
Z'$ then there exists a presentation $(\lambda, Z_{C})$ of $Z$ and a
(necessarily unique) definable morphism $g:C\to C'$ such that
$g\circ p \circ \lambda = p'\circ \lambda'_{|Z}$, where $p:Z_C\to C$
and $p':Z'_{C'}\to C'$ are the projections.
\end{lem}
\begin{proof}
If $Z$ is a $0$-cell the statement is clear. Now suppose that $Z$ is
a $1$-cell with some presentation $(\lambda_0, Z_{C_0})$. Then $Z'$
is also a $1$-cell with some presentation $(\lambda', Z'_{C'})$. Let
$\lambda$ be the presentation of $Z$ given by
$$
Z\to Z_{C_0\times C'}:z\mapsto (\lambda_0(z),p' \circ \lambda'(z)).
$$
By the non archimedean property this is indeed a presentation of the
$1$-cell $Z$, and the uniqueness and existence of $g$ for this
$\lambda$ is clear.
\end{proof}

\begin{prop}\label{refprop} Let $X$ be a definable subassignment of
$S[1,0,0]$ with $S$ in $\Def_k$ and let $\varphi_i$ be a function
in $\cC (X)$ for $i=1,2$. Let $\cP_i$  be a cell decomposition of
$X$ adapted to $\varphi_i$ for $i=1,2$. Then there exists a cell
decomposition $\cP$ of $X$ such that $\cP \prec \cP_i$ for
$i=1,2$. Such $\cP$ is automatically adapted to both $\varphi_1$
and $\varphi_2$. Similar statements hold for $\varphi_i$ in $\cC_+
(X)$, in $\cP(X)$, in $\cP_+(X)$, in $K_0(\RDef_X)$, and in
$SK_0(\RDef_X)$.
\end{prop}
\begin{proof}It follows from Lemma \ref{refin:g} that if $\cP$ is a
refinement of $\cP_1$ then $\cP$ is automatically adapted to
$\varphi_1$. Thus we only have to show that there is a common
refinement $\cP$ of $\cP_1$ and $\cP_2$, but this follows at once
from the definitions and Theorem \ref{np}.
\end{proof}

\subsection{Bicells}\label{bice}We will have, for technical
reasons, to consider bicells, that is, cells with 2 special
variables.
More precisely, they will be needed for a basic version of
Fubini's Theorem.

Fix $S$ in $\Def_k$. As for the definition of cells in \ref{ce}, we
shall first define bicells $Z\subset S[2,0,0]$ with base $C\subset
S$, and then with more general base $C\subset S[0,s,i]$.
 Let first $C$ be a definable subassigment of $S$. Let  $\alpha$ and
$\beta$ be definable morphisms $C \rightarrow \ZZ$, $\xi$ and $\eta$
definable morphisms $C \rightarrow h_{\GG_{m, k}}$, $c$ a definable
morphism $C \rightarrow h [1,0,0]$, and let $d$ be a definable
morphism $C[1,0,0]\rightarrow h [1,0,0]$. We further assume that
either, for every point $y = (y_0, K)$ in $C$, the function $u
\mapsto d (y, u)$ is constant on $K \llp t \rrp$, or, for every
point $y = (y_0, K)$ in $C$, it is injective on $K \llp t \rrp$.

The bicell $Z_{C, \alpha, \beta, \xi, \eta, c, d}$ with basis $C$
is the definable subassignment of $S [2,0,0]$ defined by
\begin{gather*}
y \in C\\
\ord (z - d(y,u)) = \alpha (y)\\
\ac (z - d (y,u)) = \xi (y)\\
\ord (u - c (y)) = \beta (y)\\
\ac (u - c(y)) = \eta (y),
\end{gather*}
where $y$ denotes the $S$-variable, $z$ the first $\AA^1_{k
\llp t \rrp}$-variable and $u$ the second $\AA^1_{k  \llp t \rrp}$-variable.

Similarly, we define the bicell $Z'_{C, \beta, \eta, c, d}$ as the
definable subassignment of $S [2,0,0]$ defined by
\begin{gather*}
y \in C\\
z = d(y,u) \\
\ord (u - c (y)) = \beta (y)\\
\ac (u - c(y)) = \eta (y),
\end{gather*}
the bicell $Z''_{C, \alpha, \xi, c, d}$ as the definable
subassignment defined by
\begin{gather*}
y \in C\\
\ord (z - d (y,u)) = \alpha (y)\\
\ac (z - d (y,u)) = \xi (y)\\
u = c (y),
\end{gather*}
and the bicell $Z'''_{C, c, d}$ as the definable subassignment
defined by
\begin{gather*}
y \in C\\
z = d (y,u) \\
u = c (y).
\end{gather*}

Now we can define bicells similarly to \ref{ce} with general base
$C\subset S [0, s, i]$,
 except that we have 4 types of bicells instead of 2 types of cells.
 A definable subassignment $Z$ of $S[2,0,0]$ will be called a
$(1,1)$-bicell (resp.~a $(1,0)$-bicell, a $(0,1)$-bicell, a
$(0,0)$-bicell), if there exists a definable isomorphism
$$\lambda : Z \longrightarrow Z_C\subset S [2, s, i],$$ and a bicell $Z_C=Z_{C,
\alpha, \beta, \xi, \eta, c, d}$ (resp.~$Z_C=Z'_{C, \beta, \eta,
c, d}$, $Z_C=Z''_{C, \alpha, \xi, c, d}$, $Z_C=Z'''_{C, c, d}$)
with basis $C\subset S [0, s, i]$, such that the morphism
$\pi\circ\lambda$, with $\pi$ the projection on the $S [2, 0,
0]$-factor, restricts to the identity on $Z$.

We shall call the data $(\lambda, Z_C)$ with $Z_C$ of one of the
above forms a presentation of the cell $Z$.

We define similarly to \ref{cd} bicell decompositions of a
definable subassignment $Z$ of $S [2, 0, 0]$ and bicell
decompositions of $Z$ adapted to a given function $\varphi$ in
$\cC (Z)$.

The following statement is an easy consequence of Theorem \ref{np}
and its Corollary  \ref{inj-constant}:

\begin{prop}\label{binp}
\begin{enumerate}
\item[(1)]Every  definable subassignment $Z$ of $S [2, 0, 0]$
admits
 a bicell decomposition.
\item[(2)]For every $\varphi$ in $\cC (Z)$ there exists a bicell
decomposition of $Z$ adapted to the function $\varphi$, namely,
there exists a finite partition of $Z$ into bicells $Z_i$ with
presentation $(\lambda_i,Z_{C_i})$, such that $\varphi_{| Z_i} =
\lambda_i^{*} p_i^{*} (\psi_i)$, with $\psi_i$ in $\cC (C_i)$ and
$p_i : Z_{C_i} \rightarrow C_i$ the projection. Similar statements
hold for $\varphi$ in $\cC_+ (Z)$, in $\cP(Z)$, in $\cP_+(Z)$, in
$K_0(\RDef_Z)$, and in $SK_0(\RDef_Z)$.
\end{enumerate}
\end{prop}
\begin{proof}
First apply Theorem \ref{np} to obtain a partition of $Z$ into
cells, adapted to $\varphi$. We  apply now Corollary
\ref{inj-constant} to each center to partition each basis. By
Lemma \ref{rem:presen} this yields a partition of $Z$ into cells,
refining the previous partition. We finish the proof by applying
Theorem \ref{np} to each basis $C$ of the occurring cells and the
functions $\psi$ in $\cC(C)$ corresponding to $\varphi$ as in
Theorem \ref{np} (2).
\end{proof}

\subsection{Analyticity and cell decomposition}\label{sec:an:cel}

We consider the expansion $\LPas^*$ of $\LPas$ which is obtained
by adding the following function symbols for each integer $n>0$:
\begin{enumerate}
 \item[(1)]The symbol  $^{-1}:{\rm Val}\to {\rm Val}$ for the field inverse
on ${\rm Val}\setminus\{0\}$ extended by  $0^{-1}=0$.
 \item[(2)] The symbol $(\cdot,\cdot)^{1/n}: {\rm Val}\times{\rm Res}\to {\rm Val}$
 for the function sending $(x,\xi)$ to the (unique) element $y$
 with $y^n=x$ and $\ac(y)=\xi$ whenever $\xi^n=\ac(x)$ and $\ord(x)\equiv 0\bmod n$, and to $0$
otherwise.
 \item[(3)] The symbol $h_n:{\rm Val}^{n+1}\times{\rm Res}\to {\rm Val}$ for the function
 sending $(a_0,\ldots,a_{n},\xi)$ to the (unique) element $y$
satisfying $\ord(y)=0$, $\ac(y)=\xi$ and $\sum_{i=0}^{n} a_{i}
y^i=0$ whenever $\xi\not=0$, $\ord(a_i)\geq 0$, $\sum_{i\in
I}\ac(a_{i}) \xi^i=0$, and $\sum_{i\in I}i\ac(a_{i})
\xi^{i-1}\not=0$ for $I=\{i\mid \ord (a_i)=0\}$, and to $0$
otherwise.
\end{enumerate}

The following is a fundamental structure result for definable
functions with values in the valued field. Its proof uses an
analogue of Lemma 3.7 of \cite{Pas} which will be contained in
\cite{CLR}.
\begin{theorem}%[\cite{CLR}]
\label{normal} Let  $f:X\to h[1,0,0]$ be a morphism in $\Def_k$.
Then there exists a definable isomorphism $g:X\to X'\subset
X[0,n,r]$ over $X$ and a $\LPas^*$-term $\tilde f$ in variables
running over $X'$ such that $f= \tilde f\circ g$.
\end{theorem}
\begin{proof} Let $\varphi$ be a formula describing the graph of $f$
and suppose that $\varphi$ is of the form (\ref{eq:man}) as in the
proof of Theorem \ref{compdim}. Let $p_j$ be all polynomials in the
${\rm Val}$-sort which appear in $\varphi$. With exactly the same
proof as the proof of the Denef-Pas cell decomposition in \cite{Pas}
where one replaces the words strongly definable function by
$\LPas^*$-term (also in Lemma 3.7 of \cite{Pas}, cf.~\cite{CLR}) and
assuming quantifier elimination, one shows that there exists a cell
decomposition of $X[1,0,0]$ adapted to $\ord(p_j)$ for each $j$ such
that the centers $c_i$ of the occurring cells $Z_i$ are
$\LPas^*$-terms. For each cell $Z_i$, let $Z'_i$ be
$\lambda_i^{-1}({\rm Graph}(c_i))\cap {\rm Graph}(f)$, where
$\lambda_i$ is the representation of $Z_i$. Clearly each $Z_i'$ is a
$0$-cell with presentation the restriction of $\lambda_i$ to $Z_i'$.
It follows from the description (\ref{eq:man}) of $\varphi$ that for
each point $x$ in $X$ the point  $f(x)$ is a zero of at least one of
the polynomials $p_i(x,\cdot)$, and by cell decomposition that at
least one of the centers gives this zero. Hence, the cells $Z_i'$
form a partition of the graph of $f$ and one concludes that the
restriction of $f$ to each of finitely many pieces in a partition of
$X$ satisfies the statement. One can glue $s$ pieces together using
extra parameters contained in the definable subassignment
$A:=\{\xi\in h[0,1,0]^s\mid \sum_i \xi_i=1\, \wedge\, (\xi_i=0\,
\vee\, \xi_i=1)\}$ to index the pieces, by noting that for each
element $a$ in $A$ there exists a definable morphism $A\to
h[1,0,0]$, given by a $\LPas^*$-term, which is the characteristic
function of $\{a\}$.
\end{proof}

Let $K$ be a field. A subset $B$ of $K \llp t \rrp$ of the form $c
+ t^{\alpha} K\llb t \rrb$ is called a ball of volume $\LL^{-
\alpha}$. A function $f : B \rightarrow K\llp t \rrp$ is called
strictly analytic if there exists a power series $\varphi
:=\sum_{i \in \NN} a_i x^i$ in $K \llp t \rrp \llb x \rrb$
converging on $t^{\alpha} K\llb t \rrb$, equivalently, $\lim_{i
\mapsto \infty} (\ord\, a_i + i \alpha) = \infty$, such that $f (c
+ y) = \varphi (y)$ for every $y$ in $t^{\alpha} K\llb t \rrb$.
Note that this definition is independent of the choice of the
center $c$ and that if $f$ is strictly analytic on $B$, its
restriction to any ball contained in $B$ is also strictly
analytic.

\begin{lem}\label{lem:tate}
Let $K$ be a field of characteristic zero. Let $f$ be strictly
analytic on a ball $B\subset K\llp t \rrp$ of volume $\LL^{-\beta}$
for some $\beta$ in $\ZZ$, and write $f'$ for its derivative.
Suppose that there exists $\alpha$ in $\ZZ$ such that $\ord\, f'(x)=
\alpha$ for every $x$ in $B$. Then the image of $f$ is contained in
a ball of volume $\LL^{-\alpha-\beta}$ and cannot be contained in a
ball of volume $\LL^{-\alpha-\beta-1}$. For every $x_0$ in $tB$, the
restriction of $f$ to $x_0+tB$ is a bianalytic bijection onto a ball
of volume $\LL^{-\alpha-\beta-1}$ with strictly analytic inverse.
Also, $\ord (f(z)-f(y))=\alpha+\ord(z-y)$ for all $y, z$ in $x_0+tB$
and $x_0$ in $B$.
\end{lem}

\begin{proof}
We may assume that $B=t^{\beta} K\llb t \rrb$. Write $f =\sum_i a_i
x^i$ in $K\llp t \rrp\llb x \rrb$. By replacing $x$ by $x/t^\beta$ and $f(x)$ by
$(f(x)-a_0)/t^\alpha$ we may suppose that $a_0=\alpha=\beta=0$.
 First we prove that $\ord\,
a_1=0$ and $\ord\, a_i\geq 0$ for all $i$. Let $I$ be $\{i\mid
\ord\, a_i=\min (\ord\, a_j)_j\}$, let $p$ be the polynomial
$\sum_{i\in I} a_ix^i$ and let $p'$ be its derivative. If $p'=a_1$
this is trivial. Suppose that the degree of $p'$ is $>0$. Since
$K$ is infinite there exists for any $b$ in $K$ an element $y_0$ in
$K$ such that $\sum_{i\in I, i>0} i\ac(a_i)y_0^{i-1}\not=b$. Taking $b=0$, it follows that $\ord\, a_i=0$ for $i$ in $I$. Thus,
$1$ belongs to $I$ since otherwise $\ord\, f'(0)>0$. Fix $c$ in $K\llp t \rrp$.
Taking $b=\ac(c)$ if $\ord\, c=0$ and $b=0$ otherwise, it follows
that the image of $f$ cannot be contained in the set $c+tK\llb t \rrb$.

Now fix $x_0$ in $K$. It is clear that $f$ maps $x_0+t K\llb t
\rrb$ into $f(x_0)+tK\llb t \rrb$. The statement about
bianaliticity is well known and follows from the Inverse Function
Theorem for complete fields, stated in \cite{Igusa}, Corollary
2.2.1(ii). The statement about the orders follows easily by
developing $f$ into power series around $z$.
\end{proof}

\begin{theorem}\label{thm:ancell}Let $X$ be in $\Def_k$, $Z$ be a definable subassignment
of $X[1,0,0]$, and let $f:Z\to h[1,0,0]$ be a definable morphism.
Then there exists a cell decomposition of $Z$ into cells $Z_i$
such that the following conditions hold for every  $\xi$ in $C_i$, for
every  $K$ in $ {\rm Field}_{k(\xi)}$, and for every $1$-cell $Z_i$ with
presentation $\lambda_i:Z_i\to Z_{C_i}=Z_{C_i,\alpha_i, \xi_i,
c_i}$ and with projections $p_i:Z_{C_i}\to C_i$, $\pi_i:Z_{C_i}\to
h[1,0,0]$:
\begin{enumerate}
\item[(1)] The set $\pi_i(p_i^{-1}(\xi))(K)$ is either empty or a
ball of volume $\LL^{- \alpha_i(\xi)-1}$.  %(Note that this volume is independent of $K$).
 \item[(2)] When $\pi_i(p_i^{-1}(\xi))(K)$ is nonempty, the function
$$g_{\xi,K}:\pi_i(p_i^{-1}(\xi))(K)\to
K\llp t \rrp:x\longmapsto f\circ\lambda^{-1}(\xi,x)$$
%$f_{\vert p_i^{-1}(\xi)(K)} \circ \pi_{i, \xi}^{-1}$
is strictly analytic. % on the ball $\pi_i(p_i^{-1}(\xi))(K)$.

\end{enumerate}
We can furthermore ensure that either $g_{\xi,K}$ is constant or
(3), (4), and (5) hold.
\begin{enumerate}
 \item[(3)] There exists a definable morphism $\beta_i:C_i\to
h[0,0,1]$ such that
$$
\ord\, \frac{\partial }{\partial x}g_{\xi,K}(x)=\beta_i(\xi)
$$
for all $x$ in $\pi_i(p_i^{-1}(\xi))(K)$.  %(Note that $\beta_i(\xi)$ is independent of $K$ and $x$).
 \item[(4)] When
$\pi_i(p_i^{-1}(\xi))(K)$ is nonempty, the map $g_{\xi,K}$ is a
bijection onto a ball of volume $\LL^{-\alpha_i(\xi) - 1 -
\beta_i(\xi)}$.

\item[(5)] For every
$x,y$ in $\pi_i(p_i^{-1}(\xi))(K)$,
$\ord
(g_{\xi,K}(x)-g_{\xi,K}(y))=\beta_i(\xi)+\ord(x-y)$.

\end{enumerate}
\end{theorem}
Note that, in this theorem, $\alpha_i(\xi)$ and $\beta_i(\xi)$ are
independent of $K$.

\begin{proof}[Proof of Theorem \ref{thm:ancell}]
Statement (1) holds
automatically if $\alpha_i$ is a coordinate function in the image of the presentation $\lambda_i$.
\par
Take a $1$-cell $Z$ with presentation $\lambda :Z\to Z_C$ and
projection $p : Z_C \rightarrow C$. By Theorem \ref{normal} we may
suppose that $f\circ \lambda^{-1}$ is given by a $\LPas^*$-term.
We shall now prove statement (2) by induction on the complexity of
the term $f\circ \lambda^{-1}$. Fix $\xi_0$ in $C$. The case where
$\pi_i(p_i^{-1}(\xi_0))(K)$ is empty being clear, we may assume
after translation and homothety that $\pi_i(p_i^{-1}(\xi_0))(K)$
is in fact the ball $K \llb t \rrb$. Consider the term $a^{-1}$ with $a$
a term for which the statement already holds. By cell
decomposition we may assume that $\ord(a)$ and $\ac(a)$ only
depend on $\xi_0$. We may also assume $a$ is non zero. Let us
denote by $\tilde a$ the function induced by $a$ on the ball
$K\llb t \rrb$. We may write
\begin{equation}\label{nfo}
\tilde a = t^{\alpha} \eta (1 + \sum_{j \geq 1} P_j (x) t^j),
\end{equation}
with $\alpha$ in $\ZZ$, $\eta$ non zero in $K$, and $P_j$
polynomials in $K [x]$. We use here the fact that $K$ is infinite.
Note also that a series in $K \llb t \rrb \llb x \rrb$ converges on the ball
$K \llb t \rrb$ if and only if it lies in $K [x] \llb t \rrb$. Since $1 +
\sum_{j \geq 1} P_j (x) t^j$ is a unit in the ring $K [x] \llb t \rrb$,
the result follows in this case. Similarly, consider a term $(a,
\xi)^{1/n}$, with $a$ a term for which the statement already
holds. As before we may assume that the function $\tilde a$
induced by $a$ on the ball $\pi_i(p_i^{-1}(\xi_0))(K)$ is of the
form (\ref{nfo}). Furthermore we may assume that $\alpha$ lies in
$n \ZZ$ and $\eta = \xi^n$. The result follows since the series $1
+ \sum_{j \geq 1} P_j (x) t^j$ has a unique $n$-th root of the
form $1 + \sum_{j \geq 1} Q_j (x) t^j$ in the ring $K [x] \llb t \rrb$.
Now assume the term is $h_n(a_0,\ldots,a_{n},\xi)$ where the $a_i$
are terms for which the statement holds and $\xi_0\in X$. We may
assume by cell decomposition that $\xi$, $\ord(a_i)$ and
$\ac(a_i)$ only depend on $\xi_0$. Denoting by $\tilde a_i$ the
function induced by $a_i$ on the ball $\pi_i(p_i^{-1}(\xi_0))(K)$,
we may write
\begin{equation}\label{nfoee}
\tilde a_i = t^{\alpha_i} \eta_i (1 + \sum_{j \geq 1} P_{i,j} (x)
t^j),
\end{equation}
with $P_{i, j}$ in $K [x]$. We may assume that $\alpha_i \geq 0$
for all $i$, $\sum_{i \in I} \eta_i \xi^i = 0$, and $\sum_{i \in
I} i\eta_i \xi^{i -1} \not= 0$,
 where $I$ denotes the
set of $i$'s with $\alpha_i = 0$. By the usual proof of Hensel's
Lemma by successive approximations modulo higher powers of $t$,
one gets that there exists universal polynomials $Q_j$ in $K
[x_{i, \ell}]_{\sur{0 \leq i \leq n}{1 \leq \ell \leq j}}$ such
that $h_n(\tilde a_0,\ldots,\tilde a_{n},\xi)$ is equal to
\begin{equation}
h_n(\tilde a_0,\ldots,\tilde a_{n},\xi) = \xi + \sum_{j \geq 1}
Q_j (P_{i, \ell} (x)) t^j,
\end{equation}
from which the assertion follows. This  concludes the proof of
(2), the result being clear for the remaining types of terms of
the forms $+,-,\cdot$, and constants.

\par

Statement (3) follows easily: using Theorem \ref{prop:man} (iii)
and general model theory, there is a definable morphism $g$ which
is, almost everywhere, equal to the derivative of $f$ with respect
to the $h[1,0,0]$-variable; then take a refinement adapted to
$\ord\, g$.

\par
Clearly (4) is a definable condition on $\xi$ in  $C_i$. Let $C_{i1}$
be the definable subassignment of $C_i$ given by condition (4) and
$C_{i2}$ its complement in $C_i$. Set $X_j:= p_i^{-1}(C_{ij})$ for
$j=1,2$; these are cells by Lemma \ref{rem:presen}.
 It is enough to prove statement (4) for the $X_j$ and the restrictions $f_{|X_j}$.
For $X_1$ and $f_{|X_1}$, statement (4) is clear and (5) is
automatically true by (2), (4), and Lemma \ref{lem:tate}.
By Lemma \ref{lem:tate} and the construction of $X_2$, for $\xi$ in
$C_{i2}$, $p_i^{-1}(\xi)(K)$, if nonempty, is mapped under $f$ into
a proper subset of a ball of volume $\LL^{-
\alpha_i(\xi)-1-\beta_i(\xi)}$ for every $K$ in ${\rm
Field}_{k(\xi)}$.
 Let  $Y$ be the
definable subassignment of $C_{i2}[1,0,0]$
 determined by $(\exists x\in X_2) y=(p_i\circ\lambda_i(x),f(x))$. If we now apply cell
decomposition to $Y$, then the fibers of these cells will be
strictly smaller balls than those of volume $\LL^{-
\alpha_i(\xi)-1-\beta_i(\xi)}$, by the definition of cells and the
construction. By a fiber product argument, we may assume that the
new  parameters we just obtained, as well as those for $C_i$, are
already present as coordinate functions for $X_2$. If we now apply
again  cell decomposition to $X_2$, the fibers of the cells in
$X_2$ are strictly smaller balls than those of volume $\LL^{-
\alpha_i(\xi)-1}$ by Lemma \ref{lem:tate}. An application of Lemma
\ref{lem:tate} shows that (4) and (5) for $X_2$ and $f_{|X_2}$
hold on this cell decomposition of $X_2$.
\end{proof}

\begin{cor}\label{cor:reciprokecel}
Let $f:S[1,0,0]\to S[1,0,0]$ be a definable isomorphism over $S$
and let $\pi:S[1,0,0]\to h[1,0,0]$ be the projection. Then there
exists a finite partition of $S[1,0,0]$ into cells $Z_i$ with
presentation $\lambda_i:Z_i\to Z_{C_i}$ and projection
$p_i:Z_{C_i}\to C_i$ such that the $f(Z_i)$ are cells with
presentation $(p_i\circ \lambda_i\circ f^{-1},\pi):x\mapsto
(p_i\circ \lambda_i\circ f^{-1}(x),\pi(x))$ and such that (1) up
to (5) of Theorem \ref{thm:ancell} are fulfilled for each $1$-cell
$Z_i$ and the map $\pi\circ f$.

Moreover, one can take the $Z_i$ and $f(Z_i)$ adapted to
$f^*(\varphi)$ and $\varphi$ for given $\varphi$ in
$\cC_+(S[1,0,0])$.
\end{cor}
\begin{proof}
First apply Theorem \ref{thm:ancell} to $S[1,0,0]$ and the
function $\pi\circ f$. Then apply cell decomposition to
$f^*(\varphi)$ to refine the obtained cells. The corollary
follows.
\end{proof}

\begin{remark}
Alternatively to the proofs of section \ref{sec:dim}, based on the
work by van den Dries, one could proof many of the results of
section \ref{sec:dim} using cell decomposition and Theorems
\ref{normal} and \ref{thm:ancell}.
\end{remark}

\section{Volume forms and Jacobians}\label{df1}

\subsection{Differential forms on definable subassignments}
Let $W$ be of the form $W = \cX \times X \times \ZZ^r$ with $\cX$
a $k\llp t \rrp$-variety and $X$ a $k$-variety. Let $h$ be a definable
subassignment of $h_W$.
Denote
by $\cA (h)$ the ring of definable functions $h \rightarrow
h_{\AA^1_{k \llp t \rrp}}$ on $h$. We want to define, for every integer
$i$ in $\NN$, an $\cA (h)$-module $\Omega^i (h)$ of definable
$i$-forms on $h$,
which we do in (\ref{defdiffi}).

First consider a $k\llp t \rrp$-variety  $\cY$ and the sheaf
$\Omega^i_{\cY}$ of degree $i$ algebraic differential forms on
$\cY$, namely, the $i$-th exterior product of the sheaf
$\Omega^1_{\cY}$ of K\"ahler differentials.
Denote by $\cA_{\cY}$ the Zariski
sheaf
 associated to the
presheaf\footnote{This presheaf is actually
 already a sheaf.}
$U \mapsto \cA (h_U)$ on $\cY$.
Both $\Omega^i_{\cY}$ and $\cA_{\cY}$ are sheaves of $\cO_{\cY} (=
\Omega^0_{\cY})$-modules, so we can consider the sheaf
$$
\Omega^i_{h_{\cY}} := \cA_{\cY} \otimes_{\cO_{\cY}} \Omega^i_{\cY}
$$
of definable degree $i$ differential forms on $\cY$.
Note that, in general, the module of global sections
$\Omega^i_{h_{\cY}} (\cY)$ is much bigger than $\cA_{\cY} (\cY)
\otimes_{\cO_{\cY} (\cY)} \Omega^i_{\cY} (\cY)$.

\par
Now let $\cY$ be the subvariety of $\cX$ which is the Zariski
closure of the image of $h$ under the projection $\pi : h_W
\rightarrow h_{\cX}$. Using the ring morphism $\cA (h_{\cY})\to
\cA (h):f\mapsto f\circ\pi$, we define the $\cA (h)$-module
$\Omega^i (h)$ of definable $i$-forms on $h$ as
\begin{equation}\label{defdiffi}
\Omega^i (h) := \cA (h) \otimes_{\cA (h_{\cY})} \Omega^i_{h_{\cY}}
(\cY).
\end{equation}
Note that $\Omega^0 (h) = \cA (h)$.

Let $d$ be the $K$-dimension of $h$. We denote by $\cA^< (h)$ the
ideal of $\cA (h)$ consisting of definable functions on $h$ vanishing on the
complement of a definable subassignment of $K$-dimension $< d$,
and we set
$$\tilde \Omega^d (h):=\Omega^d (h)/\bigl(\cA^< (h)\Omega^d (h)\bigr).
 $$
It is a free $\cA(h)$-module of rank $1$ since, for $\cY$ a $k
\llp t \rrp$-variety of dimension $d$, the sheaf $\Omega^d_{\cY}$ is
locally free of rank one away from the singular locus.

\par
Let $f: h' \to h $ be a definable morphism between two definable
subassignments. Assume $h$ and $h'$ are both of $K$-dimension $d$
and that the fibers of $f$ all have $K$-dimension $0$. Under these
conditions, we define a natural pullback morphism
\begin{equation}\label{ouy0}
f^* : \tilde \Omega (h) \longrightarrow \tilde \Omega (h')
\end{equation}
as follows. Let $\omega$ be in $\tilde \Omega (h)$. By Theorem
\ref{prop:man} (iii), there exist definable subassignments
$Z\subset h$, $Z'\subset h'$ such that, for each $K$ in ${\rm
Field}_{k}$, $Z(K)$ and $Z'(K)$ are $K\llp t \rrp$-analytic manifolds,
$\Kdim (h\setminus Z)<\Kdim h$, $\Kdim (h'\setminus Z')<\Kdim h'$,
$f_K:=f(K)_{|Z'(K)}$ is $K\llp t \rrp$-analytic on $Z'(K)$, and such
that $\omega$ induces a $K\llp t \rrp$-analytic $d$-form $\omega_K$ on
$Z(K)$. Using partial differentials with respect to local
coordinates\footnote{On an affine piece, these local coordinates
can be taken out of the coordinate functions on an embedding
affine space.} on $Z(K)$, it is clear that there exists a unique
definable $d$-form $\omega'$ in $\tilde \Omega (h')$ which induces
a $K\llp t \rrp$-analytic $d$-form $\omega_K'$ on $Z'(K)$ with
$\omega_K'=f_{K}^*(\omega_K)$ for each $K$. Define $f^*(\omega)$
as the class of $\omega'$ in $\tilde \Omega (h')$.

\subsection{Volume forms on definable subassignments}\label{vol}
Let $h$ be a definable subassignment of $h_W$, $W = \cX \times X
\times \ZZ^r$. Assume $h$ is of $K$-dimension $d$. There is a
canonical morphism of
commutative semigroups
$$
\lambda : \cA (h) \longrightarrow C_+^d (h)
$$
sending a function $f$ to the class of  $\LL^{- \ord\, f}$, with
the convention $\LL^{ - \infty} = 0$.
 Define the space $\vert \tilde
\Omega \vert_+(h)$ of definable positive volume forms on $h$ as
the quotient of the free abelian semi-group on symbols $(\omega,
g)$, with $\omega$ in $\tilde \Omega^d (h)$ and $g$ in $C^d_+
(h)$, by the relations:
\begin{eqnarray*}
 (f \omega, g)&=&(\omega, \lambda (f) g) \\
(\omega, g + g') &=&(\omega, g) + (\omega, g')\\
(\omega, 0) &=&0,
\end{eqnarray*}
for all $f$ in $ \cA (h)$ and $g'$ in $C^d_+ (h)$. We shall  write
$g |\omega|$ for the class of $(\omega, g)$, so that $g \vert f
\omega \vert = g \LL^{- \ord\, f} \vert \omega \vert$. In
particular, if $\omega$ is a differential form in $\tilde \Omega^d
(h)$ (or in $\Omega^d (h)$), we shall denote by $|\omega|$ the
class of $(\omega, 1)$ in $\vert \tilde \Omega \vert_+ (h)$. If
$h'$ is a definable subassignment of $h$, there is a natural
restriction morphism $\vert \tilde \Omega \vert_+ (h) \rightarrow
\vert \tilde \Omega \vert_+ (h')$. When $h'$ is of $K$-dimension
$d$, it is induced by restriction of differential forms and
Functions. When $h'$ is of $K$-dimension $< d$, we define it to be
the zero morphism.

Note that $\vert \tilde \Omega \vert_+ (h)$ has a natural
structure of $C^d_+ (h)$-module. We shall say an element
$|\omega|$ with $\omega$  in $\tilde \Omega^d (h)$ is a gauge
form, if it is a generator of this $C^d_+ (h)$-module.
Gauge forms always exist, since $\tilde \Omega^d (h)$ is a free
$\cA(h)$-module of rank $1$.

Replacing $C_+^d (h)$ by $C^d (h)$, one defines similarly the $C^d
(h)$-module $\vert \tilde \Omega \vert (h)$.

Let $f: h' \to h $ be a definable morphism between two definable
subassignments. Assume $h$ and $h'$ are both of $K$-dimension $d$
and that the fibers of $f$ all have $K$-dimension $0$. Under these
assumptions, pullback of functions induces a morphism $f^* :
C^{d}_+ (h) \rightarrow C^{d}_+ (h')$. Still under these
conditions, the pullbacks $f^* : \tilde \Omega (h) \longrightarrow
\tilde \Omega (h')$ and $f^* : C^{d}_+ (h) \rightarrow C^{d}_+
(h')$ induce a natural pullback morphism
\begin{equation}\label{ouy}
f^* : \vert \tilde \Omega \vert_+ (h) \longrightarrow \vert \tilde
\Omega \vert_+ (h'),
\end{equation}
defined by sending the class of $(\omega, g)$ to the class of
$(f^* (\omega), f^* (g))$.

\subsection{Canonical volume forms}\label{canvf}
Let $h$ be a definable subassignment of $h [m, n, r]$ of
$K$-dimension $d$. We denote by $x_1$, \dots, $x_m$ the
coordinates on $\AA^m_{k\llp t \rrp}$ and we consider the $d$-forms
$\omega_I := dx_{i_1} \wedge \dots \wedge dx_{i_d}$ for $I =
\{i_1, \dots, i_d\} \subset \{1, \dots m\}$, $i_1 < \dots < i_d$.
We denote by $\vert \omega_I \vert_h$ the image of $\omega_I$ in
$\vert \tilde \Omega \vert_+ (h)$.

\begin{def-lem}\label{ccvv}There is a unique element
$\vert \omega_0\vert_h$ in $\vert \tilde \Omega \vert_+ (h)$, the
canonical volume form, such that, for every $I$, there exists
$\ZZ$-valued definable functions $\alpha_I$ and $\beta_I$ on $h$,
with $\beta_I$ only taking as values $1$ and $0$, such that
$\alpha_I+\beta_I>0$ on $h$, $\vert \omega_I \vert_h = \beta_I
\LL^{- \alpha_I} \vert \omega_0\vert_h$ in $\vert \tilde \Omega
\vert_+ (h)$, and such that $\inf_I \alpha_I = 0$.
\end{def-lem}

\begin{proof}Uniqueness is clear.
Fix a gauge form $\vert \omega \vert$ on $h$. We may write $\vert
\omega_I \vert_h = \beta_I
\LL^{ - \gamma_I}
 \vert \omega \vert$, with $\gamma_I$ and
$\beta_I$ $\ZZ$-valued definable functions on $h$, and $\beta_I$
only taking  $0$ and $1$ as values.
Clearly we may suppose that $\beta_I+\gamma_I>\inf_I \gamma_I$.
If one sets $\alpha = \inf_I \gamma_I$, then $\vert
\omega_0\vert_h :=
\LL^{ - \alpha}
 \vert \omega \vert$ satisfies
the required property.
\end{proof}

We call $\vert \omega_0\vert_h$ the canonical volume form on $h$.
It is a gauge form on $h$. It is an analogue of the canonical
volume form defined by Serre in \cite{serre} in the $p$-adic case.

\subsection{Order of jacobian}\label{jac}
Let $f:X\to Y $ be a definable morphism between two definable
subassignments of $h [m, n, r]$ and $h [m', n', r']$,
respectively. Assume $X$ and $Y$ are both of $K$-dimension $d$ and
that the fibers of $f$ all have $K$-dimension $0$. By (\ref{ouy})
we may consider $f^* \vert \omega_0\vert_{Y}$ and we may write
\begin{equation}\label{osty}
f^* \vert \omega_0\vert_{Y} = \LL^{ - \ordjac f}\vert
\omega_0\vert_{X},
\end{equation}
with $\ordjac f$ a $\ZZ$-valued function on $X$ defined outside a
definable subassignment of $K$-dimension $< d$. Since, basically,
 $\ordjac$ comes from calculating (the order of)
partial derivatives in the valued field, the restriction of $\ordjac
f$ to a definable subassignment of $K$-dimension $<d$ is a definable
morphism. Thus, ${\ordjac f}$ and $\LL^{\ordjac f}$ make sense as
Functions in $C_+^d (X)$.

\begin{prop}[Chain rule for  $\ordjac$]\label{transjac}
Let $f : X \rightarrow Y$ and $g : Y \rightarrow Z$ be definable
functions between definable subassignments of $K$-dimension $d$.
Assume the fibers of $f$ and $g$ all have $K$-dimension $0$. Then
$$\ordjac   (g \circ f)
= (\ordjac f) + ((\ordjac g) \circ f)$$ outside a definable
subassignment of $K$-dimension $< d$.
\end{prop}

\begin{proof}Follows directly from the chain rule
for the pullback of
usual differential
forms.
\end{proof}

\subsection{Relative variants}\label{relvar}
Let $h\to\Lambda$ and $h'\to\Lambda$ be morphisms in $\Def_k$.
Assume that $h \rightarrow \Lambda$ and $h' \rightarrow \Lambda$
are equidimensional of relative $K$-dimension $d$. Let $f:h'\to h$
be a definable morphism whose fibers have dimension $0$ and which
commutes with the projections to $\Lambda$.
 On each fiber $h_\lambda$ and $h'_\lambda$ for
$\lambda$ in $\Lambda$, let $|\omega_{0}|_{h_\lambda}$,
resp.~$|\omega_{0}|_{h'_\lambda}$, be the canonical volume form in
$\vert \tilde \Omega \vert_+ (h_\lambda)$ and $\vert \tilde \Omega
\vert_+ (h'_\lambda)$.
 For every $\lambda$ in $\Lambda$, let $f_\lambda:h'_\lambda\to
h_\lambda$ be the map induced by $f$ and write $\ordjac f_\lambda$
for the Function in $C_+^d (h'_\lambda)$ defined in \ref{jac}. As in
the non-relative setting \ref{jac}, behind this are partial
derivatives with respect to valued field variables, which are
compatible with definability by the $\varepsilon$, $\delta$
definitions of partial derivatives. Thus, by construction there
exists a unique Function $\ordjac_{\Lambda}f$ in $C_+^d
(h'\to\Lambda)$ which is the class of a definable morphism $h'\to
\ZZ$ such that the fiber of $\ordjac_{\Lambda}f$ at $\lambda$ equals
$\ordjac f_\lambda$ for every $\lambda$ in $\Lambda$.

 \par
%Writing $p : h \rightarrow \Lambda$ for the structure morphism,
We have a commutative diagram:
\begin{equation*}\xymatrix{
h \ar[rr]^{i} \ar[dr]_{}& &\Gamma \subset h \times \Lambda \ar[dl] \\
&\Lambda,& }
\end{equation*}
where $i$ denotes the isomorphism to the graph $\Gamma$ of
$h\rightarrow \Lambda$. By construction, $\ordjac_{\Lambda} i=0$
holds.

\begin{remark}\label{rem:rel}
It is also possible to define relative analogues of  $\Omega^i
(h)$, $\tilde \Omega^d (h)$, $|\tilde \Omega |_+ (h ) $, $|\tilde
\Omega | (h) $, pullbacks, and $\vert \omega_0\vert_h$. We won't
pursue this.
\end{remark}

\subsection{Models and volume forms}\label{mod/vol}The
following construction will not be needed until \S \kern .15em
\ref{compa}. Let $\cX^0$ be an algebraic variety over $\Spec
k\llb t \rrb$, say flat over $\Spec k\llb t \rrb$. Set $\cX := \cX^0
\otimes_{\Spec k\llb t \rrb}\Spec k \llp t \rrp$. In other words $\cX$ is the
generic fiber of $\cX^0$ and $\cX^0$ is a model of $\cX$. Assume
$\cX$ is of dimension $d$. Let us denote by $U^0$ the largest open
subset of $\cX^0$ on which the sheaf $\Omega^d_{\cX^0 | k\llb t \rrb}$
is locally free of rank 1 over $k\llb t \rrb$. Its generic fiber $U :=
U^0 \otimes_{\Spec k\llb t \rrb}\Spec k \llp t \rrp$ may be identified with
the smooth locus of $\cX$
 when $\cX$ is of pure dimension $d$.
Let us choose a
finite cover of $U^0$ by open subsets $U_i^0$ on which the sheaf
$\Omega^d_{\cX^0 | k\llb t \rrb}$ is generated by a non zero form
$\omega_i$ in $\Omega^d_{U_i^0 | k\llb t \rrb} (U_i^0)$. Each form
$\omega_i$ gives rise to a volume form $|\omega_i|$ in $\vert
\Omega \vert_+ (h_{U_i})$, where $U_i$ denotes the generic fiber
of $U_i$. The subsets $U_i$ form an open  cover of $U$. Clearly
there exists a unique element $|\omega_0|$ in $\vert \tilde \Omega
\vert_+ (h_{\cX})$ such that $|\omega_0|_{| h_{U_i}} = |\omega_i|$
in $\vert \tilde \Omega \vert_+ (h_{U_i})$. Furthermore,
$|\omega_0|$ only depends on the model $\cX^0$, not on the choice
of the cover by open subsets $U_i^0$.

\section{Integrals in dimension one}\label{dim1}
This section is only needed to show that integrals in dimension
$1$, as axiomatized by Theorems \ref{mt} and \ref{mtr}, are well
defined, and satisfy a basic change of variable formula. These
results will be used in the proofs of Theorems \ref{mt} and
\ref{mtr} and the general change of variables formulas.

\subsection{Relative integrals relative to the projection
$S[1,0,0]\to S$}\label{di1} Let $S$ be in $\Def_k$ and let $\varphi$
be in $C_+(S[1,0,0]\to S)$. Since $C_+(S[1,0,0]\to S)=\oplus_{i=0}^1
C^i_+(S[1,0,0]\to S)$, we can write
$\varphi=[\varphi_0]+[\varphi_1]$ with $\varphi_i$ in
$\cC_+(S[1,0,0])$ of relative dimension $i$ relative to the
projection $S[1,0,0]\to S$ for $i=0,1$. Let $\cP_i$ for $i=0,1$ be a
cell decomposition of $S[1,0,0]$ adapted to $\varphi_i$ as in
Theorem \ref{np} and set $\cP_{ii}=\{Z\in \cP_i\mid \mbox{ $Z$ is a
$i$-cell}\}$. Fix $Z_i$ in $\cP_{ii}$ for $i=0,1$. The cell $Z_i$
has a presentation $\lambda_i : Z_i \rightarrow Z_{C_i}$ with
$Z_{C_1} =Z_{C_1,\alpha_1, \xi_1, c_1} \subset S[1,s_1,r_1]$, and
$Z_{C_0}=Z_{C_0, c_0} \subset S[1,s_0,r_0]$, for some $r_i,s_i\geq
0$, $r_0=0$, some $C_i\subset S[0,s_i,r_i]$,
  and some definable morphisms $c_i$, $\alpha_1$, and
$\xi_1$, for $i=1,2$. By Theorem \ref{np}, there is $\psi_i$ in
$\cC_+(C_i)$ such that
$$
\varphi_{i|Z_i} = \lambda_i^{*} p_i^{*} (\psi_i),
$$
where $p_i$ is the projection $Z_{C_i} \rightarrow C_i$, $i=0,1$.
Note that $\psi_i$ is unique for fixed $\varphi_i$ since
$\lambda_i$ is an isomorphism and $p_i$ is surjective. For $i=0,1$
we write $j_i$ for the inclusion
$$j_i:C_i\to S[0,s_i,r_i],$$
and $\pi_i$ for the projection
 $$
 S[0,s_i,r_i] \rightarrow S[0,0,r_i].
 $$

\begin{def-lem}\label{int-1-cell} The following definitions do not  depend on
the choice of $\lambda_i$, $i=0,1$, where we use the above notation.
We set
$$
\mu_{S,Z_0}(\varphi_0\11_{Z_0}) := \pi_{0!}(j_{0!}(\psi_0))
$$
in $\cC_+(S)$. Also, we say $\varphi_1\11_{Z_1}$ is $S$-integrable
along $Z_1$ if $\pi_{1!}j_{1!}(\LL^{-\alpha_1-1}\psi_1)$ is
$S$-integrable and if this is the case we set
$$
\mu_{S,Z_1}(\varphi_1\11_{Z_1}) :=
\mu_{S}(\pi_{1!}(j_{1!}(\LL^{-\alpha_1-1}\psi_1)))
$$
in $\cC_+(S)$. Here, $j_{i!}$, $\pi_{i!}$, and $\mu_{S}$ are as in
sections \ref{inc}, \ref{kproj}, and \ref{piint}, respectively.
\end{def-lem}

\begin{def-lem}\label{dl:int} The following definitions do not  depend on
the choice of $\varphi_i$ and $\cP_i$, $i=0,1$, where we use the
above notation. We say $\varphi$ is $S$-integrable if $\varphi_1
\11_{Z}$ is $S$-integrable along $Z$ for each $Z$ in $\cP_{11}$. If
this is the case we define $ \mu_{S}(\varphi)$ in $\cC_+(S)$ as
$$ \mu_{S}(\varphi) :=
\sum_{i=0,1}\sum_{Z\in\cP_{ii}}\mu_{S,Z}(\varphi_i\11_{Z}).
$$
\end{def-lem}

\begin{proof}[Proof of  \ref{int-1-cell}] We first prove
independence from the choice of $\lambda_0$. Suppose there is another presentation
$\lambda_0':Z_0\to Z_{C_0',c_0'}$ and $\psi_0'$ in $\cC_+(C_0')$
with $\lambda_0'{}^*p_0'{}^*(\psi_0')=\varphi_{0|Z_0}$ with
$p'_0:Z_{C_0',c_0'}\to C_0'$ the projection. Then, by the
definition of $0$-cells and by  functoriality of the pullback,
there is a definable isomorphism $f_0:C_0\to C_0'$ over $S$ with
$f_0^*(\psi_0')=\psi_0$. Now independence from the choice of $\lambda_0$
follows from Proposition \ref{stronginj}.

Let us  prove now independence from the choice of $\lambda_1$. Let
$\lambda_1' : Z_1 \rightarrow Z_{C_1'}=Z'_{C_1', \alpha_1', \xi_1',
c_1'} \subset S[1,s_1',r_1']$ be another presentation, for some
$s_1',r_1'\geq 0$ and some definable morphisms $\alpha_1'$,
$\xi_1'$, and $c'_1$, and $\psi'_1$ in $\cC_+(C_1')$ such that
$$
\varphi_{1|Z_1} = \lambda_1'^{*} p_1'^{*} (\psi_1'),
$$
where $p_1'$ is the projection $Z_{C_1'} \rightarrow C_1'$. Write
$\pi_1'$ for the projection $S[0,s_1',r_1'] \rightarrow
S[0,0,r_1']$ and $j_1'$ for the inclusion $C_1'\to
S[0,s_1',r_1']$.
By Lemma \ref{refin:g} (with $Z$ and $Z'$ of Lemma \ref{refin:g}
both equal to $Z_1$ here), we may suppose that there is a (unique)
definable morphism $g:C_1\to C_1'$ over $S$ such that $g\circ p_1
\circ \lambda_1 = p_1'\circ \lambda_1'$. Indeed, it is enough to
compare both $\lambda_1$ and $\lambda_1'$ with the presentation
$$
\lambda_1'':Z_1\to Z_{C_1\times C_1'}:z\mapsto (p_1\circ
\lambda_1(z),\lambda_1'(z)),
$$
and by symmetry it is enough to compare $\lambda'_1$ with
$\lambda_1''$, hence we may suppose that $\lambda_1''=\lambda_1$ and
$g:C_1\to C_1'$ exists, cf.~the proof of Lemma \ref{refin:g}. By
functoriality of the pullback we have $g^*(\psi_1')=\psi_1$. The
result now follows from Lemma \ref{ballintoballs} and Proposition
\ref{mixing}.
\end{proof}
\begin{lem}\label{ballintoballs}Let $S$ be in $\Def_k$ and
let $X=Z_{C_1,\alpha_1, \xi_1, c_1}\subset S[1,0,0]$ be a $1$-cell
with basis $C_1\subset S$.
Let $\lambda
 : X \rightarrow Z_{C, \alpha, \xi, c}$
be another presentation of the $1$-cell $X$,
 with  $C\subset S[0,n,r]$
and   $\pi:C\to S[0,0,r]$ be the projection. Then
%br
\begin{equation}\label{formball}
%er
\mu_{S}(\pi_!(\LL^{-\alpha-1}))=\LL^{-\alpha_1-1} \11_{C_1}
\end{equation}
in
$\cC_+(S)$.
\end{lem}
\begin{example}\label{balinbal}Let us consider
the simple case where
$X$ is the $1$-cell
$Z_{C_1,\alpha_1, \xi_1, c_1}\subset h[1,0,0]$
with  $C_1 = h[0,0,0]$,
$\alpha_1=0$, $c_1=1$, $\xi_1=-1$.
So $X$ is the
definable subassignment of $h[1,0,0]$ given by
$\ord(x)>  0 \vee x = 0$ (a ball).
 Another presentation for $X$ is
$\lambda
 : X \rightarrow Z_{C, \alpha, \xi, c}$
 with
$C = h[0,1,0]$, $\alpha(\eta)=1$, $c(\eta)=0$ when $\eta\not = 0$,
$c(0)=t$, $\xi(\eta)=\eta$ when $\eta\not = 0$, $\xi(0)=-1$, and
$\lambda(x)=(\ac(x),x)$ for $x$ with $\ord(x)=1$ and
$\lambda(x)=(0,x)$ for $x$ with $\ord(x)>1$. Hence, the ball $X(K)$
is partitioned into smaller balls, and there are ``residue field
many'' of these smaller balls, that is, the smaller balls are
parameterized by $C$. In this example,  the formula
%br
(\ref{formball})
%er
in  Lemma \ref{ballintoballs} holds, since
\begin{eqnarray*}
\mu_{h[0,0,0]}(\pi_!(\LL^{-\alpha-1})) & = &  \mu_{h[0,0,0]}(
[C]\cdot \LL^{-2})
=\mu_{h[0,0,0]}( \LL\cdot \LL^{-2})\\
&=& \mu_{h[0,0,0]}(\LL^{-1})
= \LL^{-1}
= \LL^{-\alpha_1-1},
\end{eqnarray*}
in $\cC_+(h[0,0,0])$,
 $\mu_{h[0,0,0]}$ being trivial on $\cC_+(h[0,0,0])$.
\end{example}

\begin{example}\label{balinbal3}Let $X$ be as in Example \ref{balinbal}.
Fix $\gamma \geq 1$. Another presentation for $X$ is $\lambda
 : X \rightarrow Z_{C, \alpha, \xi, c}$
 with
$C = h_{\GG^m_k}\times \{i\in\ZZ\mid 1\leq i \leq \gamma \}\cup
(h[0,1,0]\times \{\gamma + 1\})$, $\alpha(\eta,i)=i$, $c(\eta,i)=0$
when $\eta\not=0$, $c(0,\gamma + 1)=t^{\gamma + 1}$,
$\xi(\eta,i)=\eta$ when $\eta\not = 0$, $\xi(0,\gamma + 1)=-1$, and
$\lambda(x)=(\ac(x),\ord(x),x)$ for $x$ with $\ord(x)\leq \gamma +
1$ and $\lambda(x)=(0,\gamma + 1,x)$ for $x$ with $\ord(x)>\gamma +
1$ or $x = 0$. The formula
%br
(\ref{formball})
%er
in Lemma \ref{ballintoballs} holds, since
\begin{eqnarray*}
\mu_{h[0,0,0]}(\pi_!(\LL^{-\alpha-1})) & = & \mu_{h[0,0,0]}
([\GG^m_k](\LL^{-2}+\ldots+\LL^{-\gamma - 1})+ \LL\LL^{-\gamma - 2}) \\
& = & \mu_{h[0,0,0]}
((\LL-1)(\LL^{-2}+\ldots+\LL^{-\gamma - 1})+ \LL^{- \gamma - 1}) \\
% &=& \mu_{h[0,0,0]}((\LL-1)\LL^{-2}\frac{1 - \LL^{-100}}{1-\LL^{-1}} + \LL^{-101}  )\\
%&=&\mu_{h[0,0,0]}( \LL^{-1}(1 - \LL^{-100}) + \LL^{-101}    )\\
&=& \mu_{h[0,0,0]}(\LL^{-1}) = \LL^{-1} = \LL^{-\alpha_1-1},
\end{eqnarray*}
in $\cC_+(h[0,0,0])$. Note that it may happen that $t^\gamma$ is
 not uniformly definable in a family but
this does not pose any problems by Remark \ref{simplerestr}. %and
%\ref{balinbal5}.
This example shows how one can reduce the case
where
 the center takes two values to the case
 where
 the center takes
only one value (corresponding to the first presentation of $X$ in
 Example
\ref{balinbal}).
\end{example}

\begin{remark}\label{balinbal2}
Let $X$ be as in Example \ref{balinbal}. Note that, if $X (K)$
contains a subset $Z$ of the form
$$
\ord(x-c)=\alpha,\quad \ac(x-c)=\xi,
$$
then, necessarily $\alpha\geq 0$ and $\ord (c) \geq 0$. Indeed, if
$\alpha < 0$, $Z$ is too large
 to be contained in $X (K)$, and if
 $\alpha \geq  0$ and $\ord (c) < 0$, $Z$ is disjoint from $X (K)$.
Furthermore, if $\ord (c) =  0$ then $\alpha = 0$, and  $Z = X (K)$,
 since $Z$ is not contained in $X (K)$ when
$\ord (c) =  0$ and $\alpha > 0$.
It follows that if
$\lambda
 : X \rightarrow Z_{C, \alpha, \xi, c}$ is a presentation of $X$ such that $\ord (c) = 0$ on $C$,
necessarily $C$ should be $h [0, 0, 0]$, and the presentation is
similar to the presentation $\lambda_1$ with $c_1 = 1$ replaced by
$c$ and $\xi_1 = - 1$ replaced by $- \ac (c)$.
\end{remark}

\begin{remark}\label{simplerestr}
A function $\varphi$ in $\cC_+(Z\times Z')$ with $Z$ and $Z'$ some
definable subassignments is said not to depend on the $Z'$-variables
when $\varphi$ can be build up (in finitely many steps) using
formulas which do not involve variables from $Z'$. Clearly $\varphi$
then defines a unique function in $\cC_+(Z)$ by restriction in the
obvious way.
\end{remark}

We shall prove  Lemma \ref{ballintoballs} by essentially reducing to
the previous examples, using the previous remarks.

\begin{proof}[Proof of Lemma \ref{ballintoballs}]
We start by assuming
%bf
$S = h[0, 0, 0]$
%ef
 for simplicity.
Let $K$ be a field containing $k$ such that $C_1 (K)$ is nonempty, i.e. a point
$\eta_1 (K)$. Note that $c (K) : C(K) \rightarrow K\llp t \rrp$ has finite image:
this can be seen by using a valued field quantifier free formula
defining the graph of $c$.
Now we show that  $\alpha(K) : C(K) \rightarrow\ZZ$ also has finite image.
 Namely, if  $\alpha(K)$ takes infinitely
many values, it must take arbitrarily large values, so there exists
an infinite sequence of points
$(\eta_n)_{n\in\NN}$ in
$C(K)$ such that
$\alpha(K)(\eta_n)$ is strictly increasing.
Since $c(K)$ can take only finitely many values, we may
assume  that $c(\eta_n)$ takes a constant value $\tilde c$  for $n\in\NN$.
But then the balls  $\lambda^{-1}p_C^{-1}(K)(\eta_n)$,
with $p_C:Z_{C, \alpha, \xi, c}\to C$ for
the projection, would have their $t$-adic distance with $\tilde c$ go to zero
as $n$ increases, which forces $\tilde c$ to lie in the ball $X(K)$.
Similarly, since all
balls $\lambda^{-1}p_C^{-1}(K)(\eta)$ are disjoint for different
$\eta\in C(K)$, no ball around $\tilde c$ can be of the form
$\lambda^{-1}p_C^{-1}(K)(\eta)$ for $\eta\in C(K)$ and hence the
domain of $\lambda$ can not contain $\tilde c$ (every ball around
$\tilde c$ intersects  $\lambda^{-1}p_C^{-1}(K)(\eta_n)$ for
sufficiently large $n$).

By the discussion in Remark
 \ref{balinbal2},
if for some point
$\eta$ in $C (K)$, $c (\eta)$ lies outside
$X (K)$, then $C (K)$ is a point and
$\alpha (\eta) = \alpha_1 (\eta_1 (K))$.
Reciprocally, if $\alpha (\eta) = \alpha_1 (\eta_1 (K))$,
then $c (\eta)$ lies outside
$X (K)$.

We now consider the case of a general $S$. Note that the morphism
${\rm Im} c \rightarrow C_1$ has globally finite fibers, meaning
that the number of fibers is finite and bounded uniformly in $K$,
 by quantifier
elimination of valued field quantifiers.

Denote by $q$ the projection $C \rightarrow C_1$. Let us consider
the  definable subassignment $C'_1$ of $C_1$ consisting of those
points $\eta_1 (K)$ for which there exists a point $\eta$ in $C (K)$
such that $q (\eta) = \eta_1$ and $c (\eta)$ lies outside $X (K)$.
Denote by $C''_1$ the complement  of $C'_1$ in $C_1$. By additivity
we may assume $C_1$ is either $C'_1$ or $C''_1$. By additivity and a
similar finite partitioning argument we may suppose that the number
of points in the fibers of ${\rm Im} c \rightarrow C_1$ is constant
and equal to an integer $\delta>0$. We will perform an induction
argument on $\delta>0$
%br
where we show that for $\delta=1$ the formula (\ref{formball})
holds, and that for any presentation $\lambda$ with some $\delta>1$
we can find a presentation $\lambda'$ with strictly smaller $\delta$
and such that $\lambda$ and $\lambda'$ 
%bbf clearly 
%eef
yield the same result
for the left hand side of (\ref{formball}).
%er

By Remark \ref{balinbal2}, one has that $\delta=1$ if and only if
$C_1 = C'_1$, and then the projection $C \rightarrow h [0, 0, 0]$
induces an isomorphism $q: C \rightarrow C_1$ and $\alpha = \alpha_1
\circ q$. The statement is clear in this case. Hence we may assume
the image of $c$ lies in $X$.

Now suppose that $\delta \geq 2$. For notational simplicity we shall
assume again $S = h [0, 0, 0]$, the general construction being
completely similar.

The cell $Z_{C, \alpha, \xi, c} (K)$ induces a partition of $X$ into
balls $B_{\eta}$ given by conditions
$$
\ord(x-c(\eta ))=\alpha(\eta) \quad \text{and} \quad \ac(x-c(\eta
))=\xi(\eta ),
$$
where $\eta$ runs over $C(K)$.

For every $K$ such that $C_1 (K)$ is nonempty, consider the
different centers $c_j (K)$, $j \in J (K)$ with $J (K)$ finite non
empty, where ${\rm Im}c(K)=\{c_j(K)\}$.
%bf
Let us denote by $\beta (K)$ the maximum of $\ord (c_i (K) - c_j
(K))$ for $i \not= j$ and consider the subset $B (K)$ of  ${\rm
Im}c(K) $ consisting of points $c_i (K)$ such that there exists $c_j
(K)$ with $\ord (c_i (K) - c_j (K)) = \beta (K)$. We write
 $B (K)$
%br
 (uniquely) as the
%er
  disjoint union of subsets $G_i$, $i=1,\ldots,d$,
 such that
 two points  $c_1$ and $c_2$  of
 $B (K)$ belong to the same
$G_i$ if and only if
$\ord(c_1-c_2) = \beta (K)$.
For each $G_i$ we denote by $B_i$ the
%bbf
largest ball containing $G_i$ and no point in
 ${\rm Im}c(K) \setminus G_i$
 %eef
 and by $|G_i|$ the number of elements in
$G_i$.
For each $i$ consider the barycenter $\tilde c_i (K) : =
\vert G_i \vert^{- 1} \sum_{c_j\in G_i } c_j (K)$. Clearly $\tilde
c_i (K)$ belongs to $B_i$ for each $i$. %So we have a definable
%morphism $\tilde c : C \rightarrow X$ sending $x$ to $\tilde c_i$
%with $i$ such that $c(x)$ lies in $B_i$. Fix $K$ with $C_1 (K)$ non
%empty.
%ef

Each point  $\tilde c_i (K)$ belongs to a unique ball $B_{\eta_i }$
for a unique $\eta_i \in C(K)$.
If $\eta \not= \eta_i$ and the ball $B_{\eta } \subset B_i$ does occur in the
partition we may rewrite it as
$$
\ord(x- \tilde c_i (K))=\alpha(\eta) \quad \text{and} \quad \ac(x-
\tilde c_i (K))=\xi(\eta ) + \ac (c (\eta)- \tilde c_i (K)).
$$
Hence we may assume all balls $B_{\eta} \subset B_i$  occurring  in
the partition have center
 $\tilde c_i (K)$ except for the
 ball
 containing $\tilde c_i
(K)$ which has center $c (\eta_i)$.
%bbf
Now if a ball $B_{\eta}$  occurring  in the partition but not
contained in  $B_i$
 has a  center in $B_i$, we can replace that center by
$\tilde c_i (K)$. This shows that we may suppose that each $G_i$
consists  exactly of the two points $\tilde c_i (K)$ and $c
(\eta_i)$ and that $c (\eta_i)$ is the center of exactly one ball in
the decomposition, the one containing $\tilde c_i (K)$.
%eef

When $\delta = 2$, one falls back to the computation done in Example
\ref{balinbal3} and we can reduce to the case where $\delta=1$ by
the computation done in that example.
 Indeed, with the above notation,
the two centers are $\tilde c_1 (K)$ and $c (\eta_1)$
for a unique $\eta_1\in C(K)$
 and $\ord(\tilde c_1 (K) - c (\eta_1))$ is equal to
 $\max\limits_{\eta\in C(K)}(\alpha(K)(\eta))$.

%bbf
When $\delta > 2$, one reduces to smaller $\delta$ as follows. For
every $i$, $1 \leq i \leq d$, we denote by $\gamma_i (K)$ the
supremum of $\ord (y - \tilde c_i (K))$, with $y$ running over ${\rm
Im}c(K) \setminus B_i$ and by $D_i$ the smallest ball containing
$B_i$ and $\tilde c_i (K) + t^{\gamma_i (K)}$, that is, the smallest
ball strictly containing $B_i$. Of course, one may possibly have
$D_i = D_j$ for some $i \not= j$. We denote by $\gamma (K)$ the
supremum of all $\gamma_i (K)$, and we denote by $J' $ the set of
$i$'s with $\gamma_i (K) = \gamma (K)$.

Fix $i$ in $J'$. Note that if $B_j$ is contained in $D_i$, then $j$
lies in $J'$. We denote by ${\tilde c_i}' (K)$ the barycenter of all
points in ${\rm Im}c(K) \cap  D_i$. Note that all balls in our cell
decomposition that have a center in $D_i$ but  are not contained in
$D_i$ may be rewritten so to have center ${\tilde c_i}' (K)$.
%br J to J'.
%er

We replace the presentation $\lambda$ by a presentation $\lambda'$
obtained in the following way. One keeps all balls
%br
$B_\eta$ of
%er
 $\lambda$ not contained in some $D_i$ for $i$ in $J'$,
replacing centers lying in $D_i$ by ${\tilde c_i}' (K)$. Now let us
explain how we change the presentation inside a ball $D_i$, $i \in
J'$. For each $x$ in ${\rm Im}c(K) \cap  D_i$, let  $\Gamma_x$ be
the maximal ball strictly contained in $D_i$ and containing $x$.
%br
There are finitely many such balls, and we name them $\Gamma_{j,
i}$.
%er
Note that, by construction, there exist at least two such balls.
  So,
each of the balls $\Gamma_{j, i}$ contains strictly less than
$\delta$ points in ${\rm Im}c(K)$, so we may apply
%br
the induction hypothesis to each of
%er
them. That is, we can remove all balls in the presentation $\lambda$
lying in some $\Gamma_{j, i}$, add as new balls the balls
$\Gamma_{j, i}$, taking as center ${\tilde c_i}' (K)$, except if
${\tilde c_i}' (K)$ lies in $\Gamma_{j, i}$, in which case one
 takes ${\tilde c_i}' (K) + t^{\gamma (K)}$
%br
as center.
%er
 One keeps the  balls
$B_{\eta}$ in the presentation $\lambda$ which are contained in $D_i$ but
not contained in some $\Gamma_{j, i}$, replacing their center by
${\tilde c_i}' (K)$, except when ${\tilde c_i}' (K)$ lies in $B_{\eta}$, in
which case one takes ${\tilde c_i}' (K) + t^{\gamma (K)}$
%br
as center.
%er
In this way, one gets a new presentation $\lambda'$ with strictly
smaller $\delta$, since for each $i$ in $J'$, the number of centers
lying in $D_i$ is at least $3$ for $\lambda$ and
%br
is equal to 
%er
 $2$ for $\lambda'$, and
the other centers did not change.

  Note that the writing of $t^{\gamma (K)}$
involves additional  parameters
which are harmless: one can
always allow parameters of order $\gamma$ to parameterize the center,
work relatively over these parameters,
and use Remark \ref{simplerestr} to get rid of them after
integrating.
%br 
By the previous application of the induction hypothesis, $\lambda$
and $\lambda'$ yield the same result for the left hand side of
(\ref{formball}).
%er
\end{proof}

The next Lemma is essential for the proof of Lemma-Definition
\ref{dl:int}.

\begin{lem}\label{ballinpresen}
Let $X$ be as in Lemma \ref{ballintoballs}. Let $\lambda$ be any
presentation of $X\setminus Y$ onto a $1$-cell $Z_{C, \alpha, \xi,
c}$  with basis $C$, order $\alpha$, center $c$, and angular
component $\xi$, with $Y\subset X$ a $0$-cell. Write $C\subset
S[0,n,r]$ and let $\pi:C\to S[0,0,r]$ be the projection. Then, in
$\cC_+(S)$,
\begin{equation}\label{eqlembinpres}
\mu_{S}(\pi_!(\LL^{-\alpha-1}))=\LL^{-\alpha_1-1} \11_{C_i}.
\end{equation}
\end{lem}

\begin{example}\label{balpresen}
Let $X$ be as in example \ref{balinbal}. A simple presentation
$\lambda$ for $X\setminus \{0\}$ is given by $C=h_{\GG_{m, k}}\times
\NN_0$, $c(\eta)=0$, with $\NN_0=\{a\in\ZZ\mid a>0\}$,
$$
\ord(x)=\alpha(\eta), \quad \ac(x)=\xi(\eta),
$$
where $\eta$ runs over $C$, and $\alpha$ is the projection on
$\NN_0$, and $\xi$ the projection on the multiplicative group of the
residue field $h_{\GG_{m, k}}$.
 For this example, one computes that the formula at the end of
Lemma \ref{ballinpresen} holds. Namely, in $\cC_+(h[0,0,0])$,
\begin{eqnarray*}
\mu_{h[0,0,0]}(\pi_!(\LL^{-\alpha-1})) & = &  \mu_{h[0,0,0]}(
[\GG_{m,k}]\cdot \LL^{-{\rm id}-1})
=  [\GG_{m,k}] \sum_{i>1}\LL^{-i}\\
 & = &  [\GG_{m,k}]\frac{\LL^{-2}}{(1-\LL^{-1})}  \
=  \frac{[\GG_{m,k}]\LL^{-2}}{\LL^{-1}(\LL-1)} =  \LL^{-1},
\end{eqnarray*}
since $[\GG_{m,k}]/(\LL-1)=1$, where the infinite sum is understood
as in section \ref{poi} and ${\rm id}$ is the identity function on
$\NN_0$.
\end{example}

\begin{proof}[Proof of Lemma \ref{ballinpresen}]
If $Y$ is the empty subassignment, we are in the situation of Lemma
\ref{ballintoballs} and we are done. So, by a partitioning argument
as in the proof of Lemma \ref{ballintoballs}, we may assume that
the preimage of every point $\eta_1$ of $C_1$ has a nonempty intersection with $Y$.
By the discussion in the proof of  Lemma
\ref{ballintoballs} this forces the image of $c$ to be contained in $X$.

As in the proof of  Lemma \ref{ballintoballs}, we shall assume $S =
h [0, 0, 0]$ for notational simplicity, and the constructions being
canonical, they will carry over directly to the general relative
case by working  fiberwise.
Fix $K$ with $C_1 (K)$ nonempty. For every point $x$ in $X (K)$ we
denote by $\gamma (x)$ the supremum
 of $\ord (c' - x )$ for $c'$ running
over all points in $c (C (K))$ different from $x$. Now consider a
point $y$ in $Y (K)$ and the ball $B_{y}$ defined by
$$
\ord(x- y + t^{\gamma (y)})= \gamma (y)
\quad
\text{and}
\quad
\ac(x-y + t^{\gamma (y)})= 1.$$
Note that
 no ball with center $c' \not=y$  occurring in the presentation $\lambda$
 can be contained in $B_{y}$, since it would necessarily be equal to
  $B_{y}$, which is impossible.
Note that the writing of $t^{\gamma (y)}$ is
%bf
again
%ef
harmless: one can
always allow parameters of order $\gamma$ to parameterize the center
and use Remark \ref{simplerestr} to get rid of them after
integrating. Hence $y$ belongs to
 $c (C (K))$ and  all balls occurring
in the presentation $\lambda$
that are contained in $B_{y}$ have center $y$.
It follows that the restriction of
the presentation $\lambda$
to $B_{y}$ is covered by
Example
\ref{balpresen}.
Let us remark that there exists a
definable subassignment $Y'$ such that, for every $K$,
 $Y' (K)$ is the union of the balls
$B_{y}$ when $y$ runs over $Y (K)$.
Hence,  if one considers the presentation $\lambda'$ of $X$
obtained from the presentation $\lambda$ of $X \setminus Y$
 by keeping the balls in $\lambda$ not contained
in $Y'$, removing the other ones and adding the balls $B_j$ as new
cells, the statement we have to prove follows from Lemma
\ref{ballintoballs} applied to $\lambda'$. Note that the
presentation  $\lambda'$ exists since one may view the balls $B_j$
as parameterized by $Y$.
\end{proof}

\begin{proof}[Proof of \ref{dl:int}] First we prove
independence of $\mu_S(\varphi)$ from  the choice of $\varphi_i$,
$i=0,1$. Actually, $\varphi_0$ is uniquely defined. For
$\varphi_1$, we suppose that there is
$\varphi_1'$ in $\cC_+(S[1,0,0])$ with $[\varphi_1]=[\varphi'_1]$ and
we suppose that the cell decomposition $\cP_1$ is adapted to both
$\varphi_1$ and $\varphi_1'$ (see Proposition \ref{refprop}). For
a $1$-cell $Z$ in $\cP_1$ with basis $C$, representation
$\lambda:Z\to Z_C$, and projection $p:Z_C\to C$, one has
$\psi,\psi'$ in $\cC_+(C)$ satisfying
$\lambda^*p^*(\psi)=\varphi_{1|Z}$ and
$\lambda^*p^*(\psi')=\varphi'_{1|Z}$. Since
$[\varphi_1]=[\varphi'_1]$ we must have $\psi=\psi'$ and hence
$\varphi_1\11_Z=\varphi_1'\11_Z$,
by the definitions of $1$-cells and adapted cell decompositions.
This shows that there is no dependence on the choice of
$\varphi_1$ either.

We now prove that $\mu_S(\varphi)$ does not depend on  the choice
of $\cP_0$. By Proposition \ref{refprop} it is enough to consider
a refinement $\cP_0'$ of $\cP_0$ adapted to $\varphi_0$ and to
compare $\sum_{Z\in\cP_{00}}\mu_{S,Z}(\varphi_0\11_{Z})$ with
$\sum_{Z\in\cP'_{00}}\mu_{S,Z}(\varphi_0\11_{Z})$, where
$\cP_{00}'$ is the collection of $0$-cells in $\cP_0'$.
 Clearly, for each $1$-cell $Z$ in $\cP_0$ (resp.~in $\cP_0'$) we
have $\varphi_0 \11_Z=0$, because $\varphi_0$ is of relative
dimension $0$ and $\cP_0$ (resp.~$\cP_0'$) is adapted to
$\varphi_0$.
Note that the union of two $0$-cells is a $0$-cell, and that for
two different $0$-cells $Z_1,Z_2$ in $\cP_{00}$ one has
$$
\mu_{S,Z_1\cup Z_2}(\varphi_0\11_{Z_1\cup Z_2})=
\mu_{S,Z_1}(\varphi_0\11_{Z_1})+\mu_{S,Z_2}(\varphi_0\11_{Z_2}).
$$
Let $Z$ be the union of all $0$-cells in $\cP_{00}$ and $Z'$ be  the
union of all $0$-cells in $\cP_{00}'$. Then $Z$ and $Z'$ are
$0$-cells.
Since $\cP_0'$ is a refinement of $\cP_0$, it follows that
$Z\subset Z'$ and that $Z'\setminus Z$ is also a $0$-cell. We also
have $\varphi_0\11_{Z'\setminus Z}=0$, since $\cP_0$ and $\cP'_0$
are adapted to $\varphi_0$. One computes
\begin{eqnarray*}
\mu_{S,Z'}(\varphi_0\11_{Z'}) & = & \mu_{S,Z}(\varphi_0\11_{Z}) +
\mu_{S,Z'\setminus Z}(\varphi_0\11_{Z'\setminus Z})\\
 & = & \mu_{S,Z}(\varphi_0\11_{Z})
\end{eqnarray*}
which proves that $\mu_S(\varphi)$ is independent of the choice of
$\cP_0$.\\

Now let us prove that $\mu_S(\varphi)$ is independent of the choice
of $\cP_1$. By Proposition \ref{refprop} it is enough to compare
two cell decompositions $\cP_1$ and $\cP_1'$ of $S[1,0,0]$ adapted
to $\varphi_1$ such that $\cP_1'$ is a refinement of $\cP_1$.
 Fix a $1$-cell $Z$ in $\cP_1$. Note that the union of two disjoint
$1$-cells is a single $1$-cell. Similarly, the union of two
disjoint $0$-cells is a single $0$-cell. For the disjoint union of
two $1$-cells $Z_1$, $Z_2$, adapted to $\varphi$, one has clearly
$$
\mu_{S,Z_1\cup Z_2}(\varphi_1\11_{Z_1\cup Z_2})=
\mu_{S,Z_1}(\varphi_1\11_{Z_1})+\mu_{S,Z_2}(\varphi_1\11_{Z_2}).
$$
 Also, a
$0$-cell cannot contain a $1$-cell.
 Hence, we may suppose that $Z'\subset Z$ and that $Z\setminus Z'$
is a $0$-cell, with $Z'$ in $\cP_1'$. Let $\lambda:Z\to
Z_C=Z_{C,\alpha,\xi,c}$ be a presentation of $Z$, and let
$\lambda':Z'\to Z'_{C'}=Z'_{C',\alpha',\xi',c'}$ be a presentation
of $Z'$. Write $p:Z_C\to C$ and $p':Z'_{C'}\to C'$ for the
projections.
  By a fiber product argument we may suppose that $\lambda={\rm id}:Z\to Z=Z_C$.
 Let $Z'_1$ be the definable subassignment of $Z'$ determined by
$$
x\in Z'\ \wedge\ \alpha \circ p(x) = \alpha'\circ p'\circ
\lambda'(x),
$$
and set $Z_1:=Z_1'$, $Z_2:=Z\setminus Z_1$, $Z'_2:=Z'\setminus
Z'_1$.

Then, by Proposition \ref{rem:presen}, $Z_i$ and $Z_i'$, $i=1,2$,
are either empty or $1$-cells. Also, for $x$ in $Z'_2$,
$$
\alpha \circ p(x) < \alpha'\circ p'\circ \lambda'(x),
$$
by the non archimedean property, since  $Z'\subset Z$.

Since the equalities $ \mu_{S,Z_i}(\varphi\11_{Z_i}) =
\mu_{S,Z_i'}(\varphi\11_{Z_i'})$ for $i=1,2$ imply that $
\mu_{S,Z}(\varphi\11_{Z}) = \mu_{S,Z'}(\varphi\11_{Z'})$ by the
above discussion, it is enough to prove the following claim.
\begin{claim}\label{lem:well} One has $
\mu_{S,Z_2}(\varphi\11_{Z_2}) = \mu_{S,Z_2'}(\varphi\11_{Z_2'})$.
\end{claim}
For the proof of the claim, we may suppose that $Z=Z_2$ and
$Z'=Z'_2$. It is enough to show that
$$
\mu_{S,Z'}(\11_{Z'})=\LL^{\alpha-1}\11_{C}
$$
holds in $\cC_+(S)$, which follows from Lemma \ref{ballinpresen}.
\end{proof}

\subsection{Direct image under the projection
$S[1,0,0]\to S$}\label{di2}

Let $S$ be in $\Def_k$ and write $\pi:S[1,0,0]\to S$ for the
projection. Let $\varphi$ be in $C_+(S[1,0,0])$.

We first suppose that $\varphi$ is in $C_+^d(S[1,0,0])$ for some
$d$.
Fix $\hat \varphi$ in $\cC_+^{\leq d}(S[1,0,0])$ such that $\varphi$
is the class of $\hat \varphi$.
Let $\cP$ be a cell decomposition of $S[1,0,0]$ adapted to $\hat
\varphi$ as in Theorem \ref{np} and set $\cP_{i}=\{Z\in \cP\mid
\mbox{ $Z$ is a $i$-cell}\}$ for $i=0,1$.

Fix $Z_i$ in $\cP_i$, $i=0,1$. The cell $Z_0$ has a presentation
$$\lambda_0 : Z_0
\rightarrow Z_{C_0}=Z_{C_0, c_0} \subset S[0,s_0,0],
$$
for some $s_0\geq 0$ and some definable morphism $c_0$. There is a
unique $\psi_0$ in $\cC_+(C_0)$ such that
$$
\hat \varphi_{|Z_0} = \lambda_0^{*} p_0^{*} (\psi_0),
$$
where $p_0$ is the projection $ p_0 : Z_{C_0} \rightarrow C_0$. We
write $j_0$ for the inclusion $ j_0: C_0 \to S [0,s_0, 0 ]$, and
$\pi_0$ for the projection $S[0,s_0,0] \rightarrow S$. Denote by
$\gamma: C_0 \rightarrow \ZZ$ the definable morphism $ y \mapsto
(\ordjac p_0) \circ p_0^{-1}$, where $\ordjac$ is
defined as in section \ref{jac}.

\begin{def-lem}\label{pi-1-cell} The following definitions
are  independent of the choice of $\lambda_0$, where we use the
above notation. We define
 $$
 \pi_{!Z_0,d}(\hat \varphi\11_{Z_0})
 $$
as the image of the constructible
function
 $
 \pi_{0!}(j_{0!}(\psi_0\LL^{\gamma} ))
 $
in $\cC_+^{\leq d}(S)$ under the natural morphism $\cC_+^{\leq
d}(S)\to C_+^{d}(S)$. Here, $j_{0!}$ and $\pi_{0!}$ are as in
sections \ref{inc} and \ref{kproj}, respectively.

Also, we say $\hat \varphi\11_{Z_1}$ is
$S$-integrable along $Z_1$ if $\hat \varphi\11_{Z_1}$ is
$S$-integrable along $Z_1$ as in \ref{int-1-cell}. If this is the
case,
 $$
 \mu_{S,Z_1}(\hat \varphi \11_{Z_1})
 $$
as defined in Lemma-Definition \ref{int-1-cell} lies in $\cC_+^{\leq d-1}(S)$
and we define
 $$
 \pi_{!Z_1,d}(\hat \varphi \11_{Z_1})
 $$
as the image of $\mu_{S,Z_1}(\hat \varphi\11_{Z_1})$ under the
natural morphism $\cC_+^{\leq d-1}(S)\to C_+^{d-1}(S)$.
\end{def-lem}

\begin{def-lem}\label{dl:int2} The following definitions are independent of
the choice of $\cP$ and $\hat \varphi$, where we use the above
notation. We say $\varphi$ is $S$-integrable if $\hat \varphi
\11_{Z}$ is $S$-integrable along $Z$ for each $Z$ in $\cP_{1}$. If
this is the case we define $\pi_!(\varphi)$ in $\cC_+(S)$ as
$$ \pi_!(\varphi) :=
\sum_{Z\in\cP}\pi_{!Z,d}(\hat \varphi\11_{Z}),
$$
where $\pi_{!Z,d}(\hat \varphi\11_{Z})$ is defined as in
\ref{pi-1-cell}.
\end{def-lem}

Finally we take a  general $\varphi$ in $ C_+(S[1,0,0])$ and we write
$\varphi=\sum_i \varphi_i$ with $\varphi_i$ in
$C^{i}_+(S[1,0,0])$. We set
 $$
 \pi_!(\varphi):=\sum_i \pi_!(\varphi_i),
 $$
where each $\pi_!(\varphi_i)$
is  defined as in \ref{dl:int2}. By the above discussion this is
independent of the choices.

\begin{proof}[Proof of  \ref{pi-1-cell}]
Let $\lambda_0':Z_0\to Z'_{C_0',c_0'}$ be a different presentation
of $Z_0$. Since $Z_0$ is a $0$-cell, clearly there is a definable
isomorphism $g:C_0\to C_0'$ compatible with the maps $Z_0\to C_0$
and $Z_0\to C_0'$.
By Proposition \ref{transjac} and the definition of $\ordjac$,
$\gamma= (\ordjac{p'_0}) \circ p'_0{}^{-1} \circ g$,
 where $p'_0$ denotes the projection $ p'_0 : Z_{C'_0} \rightarrow
 C'_0$.
Hence,
independence from the choice of $\lambda_0$ follows by functorial
properties of the pullback. Since $Z_0$ is a $0$-cell, one has by
Proposition \ref{prop:reld} that $\Kdim\, X=\Kdim\, \pi(X)$ for
each $X\subset Z_0$, where still $\pi:S[1,0,0]\to S$ is the
projection. Hence, it follows that $
\pi_{0!}(j_{0!}(\psi_0\LL^{\gamma} ))$ is in $\cC_+^{\leq d}(S)$.
Similarly it follows that  $\mu_{S,Z_1}(\11_{Z_1}\hat \varphi)$ is
in $\cC_+^{\leq d-1}(S)$ by Proposition \ref{prop:reld}.
\end{proof}

\begin{proof}[Proof of \ref{dl:int2}]
First we prove independence of $\pi_!(\varphi)$ from  the choice
of $\hat \varphi$. Suppose that there is $\hat
\varphi'$ in $\cC^{\leq d}_+(S[1,0,0])$
whose class in $C_+^d(S[1,0,0])$ is $\varphi$.
Then there exist $\varepsilon$ and $\varepsilon'$ in $\cC^{\leq
d-1}_+(S[1,0,0])$ such that $\hat \varphi+\varepsilon=\hat
\varphi'+\varepsilon'$. We may suppose that the cell decomposition
$\cP$ is adapted to $\hat \varphi$, $\hat \varphi'$,
$\varepsilon$, and $\varepsilon'$. By dimensional considerations
similar to the ones used in  the proof of Lemma-Definition \ref{pi-1-cell}, we get that
$$
\sum_{Z\in\cP}\pi_{!Z,d}(\varepsilon\11_{Z})=\sum_{Z\in\cP}\pi_{!Z,d}(\varepsilon'\11_{Z})=0.
$$
Thus,
$$
\sum_{Z\in\cP}\pi_{!Z,d}(\hat
\varphi\11_{Z})=\sum_{Z\in\cP}\pi_{!Z,d}(\hat \varphi'\11_{Z}),
$$
by the additivity of $\pi_{!Z,d}$,
which shows the independence from the choice of $\hat \varphi$.

We shall now prove that $\pi_!(\varphi)$ is independent of the choice
of $\cP$. By Proposition \ref{refprop} it is enough to compare two
cell decompositions $\cP$ and $\cP'$ of $S[1,0,0]$ adapted to
$\hat \varphi$ such that $\cP'$ is a refinement of $\cP$.
 Similarly as in the proof of Lemma-Definition \ref{dl:int}, we can fix a
$1$-cell $Z_1$ in $\cP$ and we may suppose that $Z_1'$ in $\cP'$ is a
$1$-cell such that $Z_{0}':=Z_1\setminus Z_1'$ is a $0$-cell in
$\cP'$.
 By dimensional considerations as in the proof of Lemma-Definition
\ref{pi-1-cell}, we find that $\pi_{!Z'_{0},d}(\hat
\varphi\11_{Z'_{0}})=0$. Hence, we only have to show that
 $$
\pi_{!Z_1,d}(\hat \varphi\11_{Z_1}) = \pi_{!Z'_1,d}(\hat
\varphi\11_{Z'_1}),
 $$
 and
$$
\sum_j \pi_{!Z_j,d}(\hat \varphi\11_{Z_j}) = \sum_i
\pi_{!Z'_i,d}(\hat \varphi\11_{Z'_j}),
$$
where the sum on the left, resp.~right, hand side is over all
$0$-cells in $\cP$, resp.~$\cP'$. The first equality follows from
Claim \ref{lem:well} in the same way as this claim is used in the
proof of Lemma-Definition \ref{dl:int}. The second equality follows in a way similar
to the statement for $0$-cells in the proof of Lemma-Definition \ref{dl:int}.
\end{proof}

\subsection{Basic properties}
\begin{prop}\label{p:pos}
Let $S$ be in $\Def_k$ and $f$, $g$ both in
$C_+(S[1,0,0]\to S)$ or both in $C_+(S[1,0,0])$.
 If $g\geq f$ and
$g$ is $S$-integrable, then $f$ is $S$-integrable.
\end{prop}
\begin{proof}
This follows immediately from Proposition \ref{geq} and the
definition of integrability in \ref{dl:int} by taking a cell
decomposition adapted to $f$ and $g$ which exists by Proposition
\ref{refprop}.
\end{proof}

\begin{prop}
[Change of variable in relative dimension 1]\label{cv1} Let $X$
and $Y$ be  definable subassignments of $S[1,0,0]$ for some $S$ in
$\Def_k$ and let $f:X\to Y$ be a definable isomorphism over $S$.
Suppose that $X$ and $Y$ are equidimensional of relative dimension
$1$ relative to the projection to $S$. Let $\varphi$ be in
$C_+^{1}(Y\to S)$. We use $\ordjac_{S} f$ as defined in section
\ref{relvar}. Then, $\varphi$ is $S$-integrable if and only if
$\LL^{-\ordjac_{S} f}f^*(\varphi)$ is $S$-integrable and if this
is the case then
 $$
 \mu_S(\varphi)=\mu_S(\LL^{-\ordjac_{S} f}f^*(\varphi))
 $$
holds in $\cC_+(S)$.
\end{prop}
\begin{proof}
This follows from Theorem \ref{thm:ancell} and its corollary. Note
that the morphisms $\beta_i\circ p_i \circ \lambda_i$ on the
$1$-cells $Z_i$, defined in Theorem \ref{thm:ancell} (3), have the
same class as $\ordjac_{S} f_{|Z_i}$.
\end{proof}

\part{Construction of the general motivic measure}

\section{Statement of the main result}\label{sec7}

\subsection{Integration} In this section, and until
section \ref{secg}, all definable
subassignments will belong to $\Def_k$. In particular they will be affine.
To be able to integrate positive
motivic constructible Functions,
we have to define
integrable positive Functions.
These, and more generally
$S$-integrable positive Functions,
will be defined inductively, as follows:

\begin{theorem}\label{mt}Let $S$ be in
$\Def_k$.
There is a unique functor from the category
$\Def_S$ to the category of abelian semigroups, $Z \mapsto {\rm I}_S C_+ (Z)$,
assigning to
every morphism $f : Z \rightarrow Y$ in
$\Def_S$
a morphism
$f_!  : {\rm I}_S C_+ (Z) \rightarrow {\rm I}_S C_+ (Y)$
and
satisfying the following axioms:
\begin{enumerate}
\item[]{{\rm A0 (Functoriality): }}
\begin{enumerate}\item[(a)]
For every composable morphisms $f$ and $g$ in $\Def_S$, $(f \circ
g)_! = f_! \circ g_!$. In particular, ${\rm id}_! = {\rm id}$.
 \item[(b)]Let
$\lambda : S \rightarrow S'$ be a morphism in $\Def_k$ and
denote by $\lambda_+ : \Def_S \rightarrow \Def_{S'}$ the functor
induced by composition with $\lambda$,
we have the inclusion
${\rm I}_{S'} C_+ (\lambda_+
(Z)) \subset {\rm I}_{S} C_+ (Z)$ for $Z$ in $\Def_S$, and for
$\varphi$ in ${\rm I}_{S'} C_+ (\lambda_+ (Z))$, $f_! (\varphi)$ is
the same Function computed  in ${\rm I}_{S}$ or in ${\rm
I}_{S'}$.
\item[(c)]If $f : X \rightarrow Y$ is a morphism
in $\Def_S$, a positive constructible Function
$\varphi$ on $X$ belongs to
${\rm I}_S C_+ (X)$
if and only
if
$\varphi$ belongs to
${\rm I}_Y C_+ (X)$ and $f_! (\varphi)$
belongs to
${\rm I}_S C_+ (Y)$.
\end{enumerate}
\medskip

\item[]{{\rm A1 (Integrability): }}
\begin{enumerate}\item[(a)]For every $Z$
in $\Def_S$,
${\rm I}_S C_+  (Z)$ is a graded subsemigroup
of
$C_+  (Z)$.

\item[(b)]${\rm I}_S C_+ (S) = C_+ (S)$.

\end{enumerate}

\medskip

\item[]{{\rm A2 (Additivity): }}Let $Z$ be a definable
subassignment in $\Def_S$. Assume $Z$ is the disjoint union of two
definable subassignments $Z_1$ and $Z_2$. Then, for every morphism
$f : Z \rightarrow Y$ in $\Def_S$, the isomorphism $C_+ (Z) \simeq
C_+ (Z_1) \oplus C_+ (Z_2)$ induces an isomorphism ${\rm I}_S C_+
(Z) \simeq {\rm I}_S C_+ (Z_1) \oplus {\rm I}_S C_+ (Z_2)$ under
which we have $f_! = f_{\vert Z_1 !} \oplus f_{\vert Z_2!}$.

\medskip

\item[]{{\rm A3 (Projection formula): }}For every morphism $f : Z \rightarrow Y$ in
$\Def_S$, and every $\alpha$ in $\cC_+ (Y)$
and $\beta$ in ${\rm I}_S C_+ (Z)$,
$\alpha f_! (\beta)$ belongs to ${\rm I}_S C_+ (Y)$
if and only if $f^*(\alpha) \beta$ is in ${\rm I}_S C_+ (Z)$.
If these conditions are verified, then
$f_! (f^*(\alpha) \beta) = \alpha f_! (\beta)$.

\medskip

\medskip
\item[]{{\rm A4 (Inclusions): }} If $i: Z \hookrightarrow Z'$ be
the inclusion between two definable subassignments of some object
in $\Def_S$, for every $\varphi$ in $\cC_+ (Z)$, $[\varphi]$ lies
in ${\rm I}_S C_+ (Z)$ if and only if $[i_! (\varphi)]$ belongs to
${\rm I}_S C_+ (Z')$, with  $i_!$ defined as in \ref{inc}. If this
is the case, then $i_! ([\varphi]) = [i_! (\varphi)]$.

\medskip
\item[]{{\rm A5 (Projection along $k$-variables): }} Let $Y$ be in
$\Def_S$.  Consider the projection $\pi : Z = Y [0, n, 0]
\rightarrow Y$. Let  $\varphi$ be in $\cC_+ (Z)$. Then $[\varphi]$
belongs to ${\rm I}_S C_+ (Z)$ if and only if $[\pi_! (\varphi)]$
belongs to ${\rm I}_S C_+ (Y)$, $\pi_!$ being defined as in
\ref{kproj}. Furthermore, when this holds, $\pi_! ([\varphi]) =
[\pi_! (\varphi)]$.

\medskip
\item[]{{\rm A6 (Projection along $\ZZ$-variables): }} Let $Y$ be
in $\Def_S$. Consider the projection $\pi : Z = Y [0, 0, r]
\rightarrow Y$. Take $\varphi$ in $\cC_+ (Z)$. Then $[\varphi]$
belongs to ${\rm I}_S C_+ (Z)$ if and only if there is a function
$\varphi'$ in $\cC_+ (Z)$ with $[\varphi'] = [\varphi]$ such that
$\varphi'$ is $\pi$-integrable in the sense of \ref{piint} and
$[\mu_S (\varphi')]$ belongs to ${\rm I}_S C (Y)$. Furthermore,
when this holds, $\pi_! ([\varphi])= [\mu_S (\varphi')]$.

\medskip
\item[]{{\rm A7 (Relative annuli): }}Let $Y$ be in $\Def_S$ and
consider definable morphisms
$\alpha : Y \rightarrow \ZZ$,
$\xi : Y \rightarrow h_{\GG_{m, k}}$, with $\GG_{m, k}$
the multiplicative group $\AA^1_k \setminus \{0\}$,
and $c : Y \rightarrow h_{\AA^1_{k  \llp t \rrp}}$.
Then, if $Z$ is the definable
subassignment of $Y [1, 0, 0]$
defined by
$\ord (z - c (y)) = \alpha (y)$
and
$\ac (z - c(y)) = \xi (y)$,
and $f : Z \rightarrow Y$ is the morphism induced by the projection
$Y \times h_{\AA^1_{k  \llp t \rrp}} \rightarrow Y$,
$[\11_Z]$ is in
${\rm I}_S C_+ (Z)$ if and only $\LL^{- \alpha - 1} [\11_Y]$
belongs to ${\rm I}_S C_+ (Y)$, and, if this is the case,
then
$$f_! ([\11_Z]) = \LL^{- \alpha - 1} [\11_Y].$$

\medskip
\item[]{{\rm A8 (Graphs): }}Let $Y$ be in $\Def_S$ and
consider a definable morphism $c : Y \rightarrow h_{\AA^1_{k  \llp t \rrp}}$.
If $Z$ is the
definable
subassignment of $Y [1, 0, 0]$
defined by $z - c(y) = 0$ and $p : Z \rightarrow Y$ is the projection,
$[\11_Z]$ is in
${\rm I}_S C_+ (Z)$ if and only $\LL^{(\ordjac p)\circ p^{-1}}$
belongs to ${\rm I}_S C_+ (Y)$, and,  if this is the case,
then
$$f_! ([\11_Z]) =
\LL^{(\ordjac p)\circ p^{-1}}.$$
\end{enumerate}
\end{theorem}

For $f : X \rightarrow S$ a morphism,
elements of ${\rm I}_S C_+ (X)$ shall be called
$S$-integrable positive Functions
(or
$f$-integrable positive Functions).

\begin{remark}Axiom A8 is a special case of
Theorem \ref{cvf} so  Theorem \ref{cvf}
could replace A8 as an axiom.
\end{remark}

\begin{remark}In general $f_!$ is a morphism of abelian semigroups
but not of graded semigroups. There is a shift by the relative
$K$-dimension, as, for instance, in axiom A7.
\end{remark}

\subsection{Motivic measure}When $f : Z \rightarrow h_{\Spec k}$
is the projection onto  the final subassignment, we write ${\rm I}
C_+ (Z)$ for ${\rm I}_{h_{\Spec k}} C_+ (Z)$. For $\varphi$ in ${\rm
I} C_+ (Z)$, we define the motivic measure $\mu (\varphi)$ as $f_!
(\varphi)$ in ${\rm I} C_+ (\Spec k) = SK_0 (\RDef_{h_{\Spec k}})
\otimes_{\NN [\LL - 1]} \AA_+$. We shall write also $\mu (Z)$ for
$\mu ([\11_Z])$, when $Z$ is a definable subassignment of $h [m, n,
0]$ such that $[\11_Z]$ is integrable. By Proposition \ref{alint}
this happens as soon as $Z$ is bounded in the sense of \ref{ibf}.

\section{Proof of Theorem \ref{mt}}\label{sec8}
Recall that in this section, and until
section \ref{secg}, all definable
subassignments belong to $\Def_k$, so,  in particular, they are    affine.

\subsection{Uniqueness}\label{uni}Using A0,
it is enough to show, for every $f : X \rightarrow S$, the
uniqueness of ${\rm I}_S C_+ (X)$ and of $f_! : {\rm I}_S C_+ (X)
\rightarrow {\rm I}_S C_+ (S) = C_+ (S)$. We consider first the
case of a projection $\pi : S \times Y \rightarrow S$ with $Y$
definable subassignment of some $h [m, n, r]$. We may assume $Y =
h [m, n, r]$. Indeed, $\pi$ may be factorized as
\begin{equation*}\xymatrix{
S \times Y
\ar[rd] \ar[rr]^{i} & &
S\times h [m, n, r] \ar[ld]\\
&S&
}
\end{equation*}
with $i$ the inclusion, so we are done by
A0 and A4.
The case where $m = 0$ is dealed with by using
A5, A6, and A0.
Let us consider now the case $m= 1$ and take $\varphi$ in
$\cC_+ (S [1, 0, 0])$.
By Theorem \ref{np}, there
exists a cell
decomposition $\cZ$ of $S [1, 0, 0]$
adapted to $\varphi$, that is,
a finite partition of $S [1, 0, 0]$
into
cells $Z_i$ with
presentation $(\lambda_i, Z_{C_i})$,
such that $\varphi_{| Z_i} = \tilde \lambda_i^{*} p_i^{*} (\psi_i)$,
with $\psi_i$ in
$\cC (C_i)$ and $p_i : Z_{C_i} \rightarrow C_i$ the projection.
Furthermore, maybe after
applying again
Theorem \ref{np} and taking a
refinement of $\cZ$, we may assume the following condition:
\begin{equation}\label{cond}\varphi_{| Z_i} \,
\text{is either zero or has the same $K$-dimension as}
\,  Z_i
\text{ for every } i.
\end{equation}
Using A2 and A4 we may reduce to the case
of the projection
$f : Z \rightarrow S$ of a cell $Z$ of $K$-dimension $d$ in
$S [1, 0, 0]$
with
presentation $(\lambda, Z_{C})$ and a function $\varphi$ in $\cC_+ (Z)$
of $K$-dimension $d$
such that $\varphi  = \tilde \lambda^{*} p^{*} (\psi)$,
with $\psi$ in
$\cC_+ (C)$ and $p : Z_{C} \rightarrow C$ the projection.
We have to decide when $[\varphi]$ belongs to
${\rm I}_S C_+ (Z)$ and if it is the case to compute
the value of $f_! ([\varphi])$.
Let us denote by $\tilde \pi : Z_C \rightarrow Z$ the restriction
of the projection
on the
$S [1, 0, 0]$-factor to $Z_C$.
The morphism
$\tilde \pi$ is the inverse of $\tilde \lambda$.
Since $\tilde \pi \circ \tilde \lambda = {\rm id}$,
it follows from A0 that
$\tilde \pi_!$ and $\tilde \lambda_!$
are mutually inverse. It follows from
A4, A5 and A6 that
$\tilde \pi_! ([\11_{Z_C}]) = [\11_Z]$,
hence
$\tilde \lambda_! ([\11_Z]) = [\11_{Z_C}]$.
So, by using the projection formula A3,
one gets that
$[\varphi]$ belongs to
${\rm I}_S C_+ (Z)$
if and only if
$p^* \psi [\11_{Z_C}]$
belongs to ${\rm I}_S C _+(Z_C)$.
By A0,
A1(b) and A3 this is equivalent to the condition
that
$\psi p_! ([\11_{Z_C}])$ belongs to
${\rm I}_S C_+ (C)$, which amounts to the case $m = 0$ already considered.
Now if
$[\varphi]$ belongs to
${\rm I}_S C_+ (Z)$, plugging in  axioms
A7 or A8 depending on the type
of the cell
$Z_C$, completely determines the value
of $f_! ([\varphi])$: it should be equal to
$h_! (\psi p_! ([\11_{Z_C}])$, with $h$  the canonical morphism
$C \rightarrow S$.

Now consider the case of a general morphism $f : X \rightarrow S$.
We factor it as $f = \pi \circ \gamma_f$, with $\gamma_f : X
\rightarrow X \times S$ the graph morphism and $\pi : X \times S
\rightarrow S$ the projection. We consider also the projection $p
: X \times S \rightarrow X$. Since $p \circ \gamma_f = {\rm id}$,
it follows from A0 and A1 that ${\rm I}_{X \times S} C_+ (X) = C_+
(X)$. Hence a Function $\phi$ in $C_+ (X)$ will belong to ${\rm
I}_{S} C_+ (X)$ if and only $\gamma_{f!} (\phi)$ belongs to ${\rm
I}_{S} C_+ (X \times S)$ and then $f_! (\phi) = \pi_! (\gamma_{f!}
(\phi))$. Hence we are left with showing the uniqueness of
$\gamma_{f!}$. It is enough to show that $\gamma_{f!} (\varphi
[\11_X])$ is uniquely determined for $\varphi$ in $\cC_+ (X)$,
since one can always reduce to that case replacing $X$ by some
subassignment and using A4 and A2. Let us denote by $\Gamma_f$ the graph
of $f$. It follows from the previous discussion of projections
that $p_! [\11_{\Gamma_f}]$ should be of the form $\LL^{\alpha}
[\11_X]$ for some definable function $\alpha$ on $X$. Since
$\LL^{\alpha }= \gamma_f^* \LL^{\alpha \circ p}$, we get
$$
[\11_{\Gamma_f}]
=
\gamma_{f!}
p_!([\11_{\Gamma_f}])
=
\gamma_{f!} (\gamma_f^* \LL^{\alpha \circ p} [\11_X])
=
\LL^{\alpha \circ p}
\gamma_{f!}([\11_X]),
$$
by using functoriality and the projection formula
hence $\gamma_{f!}([\11_X])$
should be equal to $L^{- \alpha \circ p} [\11_{\Gamma_f}]$.
Since
$\varphi = \gamma_f^* p^* \varphi$,
it follows again from the projection formula that
$\gamma_{f!} (\varphi [\11_X])$ is uniquely determined.  \qed

\subsection{Projections}\label{proj}
We will now construct
${\rm I}_S C_+ (S \times Y)$ and $\pi_!$ when
$\pi$ is the projection
$ S \times Y \rightarrow S$. We start by
assuming $Y = h [m, n, r]$, so that
$S \times Y  = S [m, n, r]$.

When $m = n = r = 0$
we set  ${\rm I}_S C_+ (S) = C_+ (S)$ and $\pi_! = {\rm id}$.

More generally,
we set ${\rm I}_S C_+ (S [0, n, 0]) = C_+ (S[0, n, 0])$
and define
$$\pi_!  : {\rm I}_S C_+ (S [0, n, 0]) \rightarrow C_+ (S)$$
by
$\pi_! ([\varphi]) = [\pi_! (\varphi)]$, for
$\varphi$ in
$\cC_+ (S[0, n, 0])$ of $K$-dimension $d$,
$\pi_! : \cC_+ (S[0, n, 0]) \rightarrow \cC_+ (S)$ being defined as
in
\ref{kproj}.

Similarly, when $m = n = 0$,
we define ${\rm I}_S C_+ (S [0, 0, r])$
as dictated by A6.
That is,
for $\varphi$ in $C_+ (Z)$,
we shall say
$\varphi$ belongs to
${\rm I}_S C_+ (S [0, 0, r])$ if and only if
there is a function
$\varphi'$ in $\cC_+ (S [0, 0, r])$ with $[\varphi'] = \varphi$ such that
$\varphi'$ is $S$-integrable in the sense of \ref{piint},and we set
$\pi_! (\varphi)= [\mu_S (\varphi')]$.
Clearly this definition is independent of the choice  of the representative
$\varphi'$.

We now consider the case when $m =0$
and $n$, $r$ are arbitrary.
In this case we may  mix both definitions.
More precisely we have the following statement, which follows from
Proposition \ref{mixing}:

\begin{def-prop}\label{sal}
Let $\varphi$ be a Function in
$C_+ (S [0, n, r])$.
Consider the following commutative diagram of projections
\begin{equation*}\label{th}\xymatrix{
&S [0, n, r] \ar[dl]_{\pi_1} \ar[dr]^{\pi'_1} \ar [dd]^{\pi} &\\
S [0, n, 0] \ar [dr]_{\pi_2}&&
S [0, 0, r] \ar [dl]^{\pi'_2}\\
&S&.}
\end{equation*}Then
$\varphi$ is $\pi_1$-integrable if and only if
$\pi'_{1!} (\varphi)$ is $\pi'_{2}$-integrable.
We then say $\varphi$ is $\pi$-integrable.
If these conditions hold
then
$\pi_{2 !} \pi_{1!}(\varphi)$ and
$ \pi'_{2!} \pi'_{1!} (\varphi)$ are equal so we may define
$\pi_! (\varphi)$ to be their commun value. \qed
\end{def-prop}

The case of the projection
$\pi : S [1, 0, 0] \rightarrow S$ has been considered in \ref{di2},
where we defined the notion of
$S$-integrability for $\varphi$ in
$C_+ (S [1, 0, 0])$ and also the  value of
$\pi_! (\varphi)$ when $\varphi$ is $S$-integrable.
We can go one step further thanks to
the following:

\begin{def-prop}\label{uht}
Let $\varphi$ be a Function in
$C_+ (S [1, n, r])$.
Consider the following commutative diagram of projections
\begin{equation*}\label{thh}\xymatrix{
&S [1, n, r] \ar[dl]_{\pi_1} \ar[dr]^{\pi'_1} \ar [dd]^{\pi} &\\
S [1, 0, 0] \ar [dr]_{\pi_2}&&
S [0, n, r] \ar [dl]^{\pi'_2}\\
&S&.}
\end{equation*}Then
the following conditions are equivalent:
\begin{enumerate}
\item[(1)]
$\varphi$ is $\pi_1$-integrable and
$\pi_{1!} (\varphi)$ is $\pi_{2}$-integrable.
\item[(2)]$\varphi$ is $\pi'_1$-integrable and
$\pi'_{1!} (\varphi)$ is $\pi'_{2}$-integrable.
\end{enumerate}
Furthermore, if these conditions are satisfied,
then
$\pi_{2 !} \pi_{1!}(\varphi) = \pi'_{2!} \pi'_{1!} (\varphi)$.
We shall  say $\varphi$  is  $S$-integrable
if  it satisfies
conditions (1) and (2)
and we shall then define $\pi_! (\varphi)$ to be the commun value
of
$\pi_{2 !} \pi_{1!}(\varphi)$ and
$ \pi'_{2!} \pi'_{1!} (\varphi)$.
\end{def-prop}

\begin{proof}
Let $\varphi$ be a Function in
$C_+ (S [1, n, r])$.
Choose a cell decomposition of $S[1,0,0]$ which is adapted to
$\pi_{1!}(\varphi)$. For every cell $Z_1\subset S[1,0,0]$ with
presentation $\tilde \lambda:Z_1\to Z_1'\subset S[1,n',r']$ in
this decomposition let $Z$ be its inverse image of $Z_1$ under
$\pi_1$.

\begin{claim} The
cells $Z$ obtained that way form a cell  decomposition of
$S[1,n,r]$, adapted to $\varphi$, and having as presentation
$$\tilde\lambda\circ \pi_1\times {\rm id}_{h[0,n,r]}:
\begin{cases}
Z\to Z'\subset
S[1,n+n',r+r']\\
(x,y)\mapsto (\tilde\lambda(x),y).
\end{cases}$$
\end{claim}

\begin{proof}The claim can be easily verified when $n=0$ and when $r=0$ and
follows in general by factorizing $\pi_1$ into projections
$S[1,n,r]\to S[1,0,r]\to S[1,0,0]$.
\end{proof}

We now consider the following commutative diagram, with
$\lambda:Z\to Z'$ and  $\tilde \lambda:Z_1\to Z_1'$ as above and
where we use the corresponding projections:
\begin{equation*}\label{hrrrr}\xymatrix@=0.3cm{
&Z \subset S [1, n,r]
\ar[dl]_{\pi_1}\ar[rr]^{\lambda}\ar[rdd]^>>{\pi'_1}&&
Z' \subset S [1, n +n',r+r']\ar[dl]_{\tilde \pi_1}\ar[rdd]^{\tilde \pi'_1}&\\
Z_1 \subset S [1, 0, 0] \ar[rr]_{\tilde \lambda}\ar[rdd]^{\pi_2}&&
Z'_1 \subset
S [1, n', r']\ar[rdd]^<<<<<{\tilde \pi_2}&&\\
&&S[0, n, r]\ar[dl]^{\pi'_2}&&
S [0, n+ n', r+r']\ar[ll]_<<<<<<<<<{\tilde \mu}\ar[dl]^{\tilde \pi'_2}\\
&S&&S [0, n', r']\ar[ll]^{\mu}&
.}
\end{equation*}
There is a unique $\psi'$ such that $\psi = \lambda^* (\psi')$.

(a)  It follows from the claim and the construction in \ref{di2}
that the statement we want to prove is verified for $[\psi']$,
that is, $[\psi']$ is $\tilde \pi_1$-integrable and $\tilde
\pi_{1!} ([\psi'])$ is $\tilde \pi_{2}$-integrable if and only if
$[\psi']$ is $\tilde \pi'_1$-integrable and $\tilde \pi'_{1!}
([\psi'])$ is $\tilde \pi'_{2}$-integrable; if these conditions
are satisfied, then $\tilde\pi_{2 !} \tilde\pi_{1!}([\psi']) =
\tilde\pi'_{2!}\tilde \pi'_{1!} ([\psi'])$.

(b) Let us remark that $[\psi]$ is $\pi_1$-integrable if and only
is $[\psi']$ is $\tilde \pi_1$-integrable and $\tilde
\pi_{1!}([\psi'])$ is $\tilde\lambda^{-1}$-integrable, and that in
this case $(\tilde \lambda^{-1})_! \tilde \pi_{1!}
([\psi'])=\pi_{1!} ([\psi])$. This follows from Proposition
\ref{stronginj} and the functoriality of the so far constructed
direct images for projections. By the claim and the construction
of integration in relative dimension 1 in \ref{di2}, $\pi_{1!}
([\psi])$ is $\pi_2$-integrable if and only if $\tilde \pi_{1!}
([\psi'])$ is $\tilde \pi_2$-integrable and $\tilde \pi_{2!}
\tilde \pi_{1!} ([\psi'])$ is $\mu$-integrable. If all  the
previous conditions are satisfied, then $\pi_{2!} \pi_{1!}
([\psi]) = \mu_! \tilde \pi_{2!} \tilde \pi_{1!} ([\psi'])$.

(c) By construction of integration in relative dimension 1,
$[\psi]$ is $\pi'_1$-integrable if and only if $[\psi']$ is
$\tilde \pi'_1$-integrable and $\tilde \pi'_{1!} ([\psi'])$ is
$\tilde \mu$-integrable. If this holds, then $\pi'_{1!} ([\psi]) =
\tilde \mu_! \tilde \pi'_{1!} ([\psi'])$. Furthermore it follows
from Proposition-Definition \ref{sal} that for a Function $g$ in
$C_+ (S [0, n+ n', r+r'])$ the condition $g$ is $\tilde
\mu$-integrable and $\tilde \mu_! (g)$ is $\pi'_2$-integrable is
equivalent to $g$ is $\tilde \pi'_2$-integrable and $\tilde
\pi'_{2!} (g)$ is $\mu$-integrable and implies that $ \pi'_{2!}
\tilde \mu_! (g) = \mu_! \tilde \pi'_{2!} (g)$.

The statement we have to prove follows directly from the conjunction
of (a), (b), and (c).\end{proof}

Now we would like
to define ${\rm I}_{S} C_+ (S [n, m, r])$
by induction on $n$
by using a factorization
\begin{equation}\label{fac}\xymatrix@1{
S [m, n, r] \ar[r]^>>>>{q}& S [m- 1, n, r]\ar[r]^>>>>{p} &S,
}
\end{equation}
with $p$ and $q$ projections,
by saying $\varphi$ in
$C_+ (S [m, n, r])$
will be $S$-integrable if it is $S[m- 1, n, r]$-integrable
and $q_! (\varphi)$ is $S$-integrable and setting
$\pi_!  (\varphi):= p_! (q_! (\varphi))$.

Since there are $m$ different
projections
$S [m, n, r] \rightarrow
S [m- 1, n, r]$, the
factorization (\ref{fac}) is
not unique,  and we have to check this definition is independent
of the factorization.

By induction it is enough to consider the
case
$(m, n, r)
=
(2, 0, 0)$.
Using a bicell decomposition
thanks to  Proposition \ref{binp}
it is enough to prove the following:

\begin{prop}\label{mainpropbic}Let $Z$ be a bicell in $S [2, 0, 0]$.
Denote by $p_1$
and $p_2$ the two projections $S [2, 0, 0] \rightarrow S[1, 0, 0]$.
Then
$[\11_Z]$ is $p_1$-integrable  and
$p_{1!} ([\11_Z])$
is  $S$-integrable if and only if
$[\11_Z]$ is $p_2$-integrable  and
$p_{2!} ([\11_Z])$
is  $S$-integrable. If these conditions hold,
then
$p_! (p_{1!} ([\11_Z])) =p_! (p_{2!} ([\11_Z]))$.
\end{prop}

\begin{proof}Let $Z$ be
bicell $Z \subset S [2, 0,0]$ with
presentation
$\lambda : Z \rightarrow Z' = Z_{C,\dots} \subset S [2, n,r]$.
Let us first note it is enough to prove the statement of the proposition
when $\lambda$ is the identity and that, in this case,
the integrability conditions are always satisfied.
To check that, let us
consider the commutative diagram
\begin{equation*}\label{hiii}\xymatrix{
&Z \subset S [2, 0,0]
\ar[dl]_{p_1}\ar[rrr]^{\lambda}\ar[rdd]^{p_2}&&&Z' \subset S [2, n,r]\ar[dl]_{p'_1}\ar[rdd]^{p'_2}&\\
S [1, 0, 0] \ar[rdd]^{p}&&&S [1, n, r]\ar[rdd]^<<<<<{p'}\ar[lll]_<<<<<{\tilde \mu}&&\\
&&S [1, 0, 0]\ar[dl]^{p}&&&S [1, n, r]\ar[lll]_{\tilde \mu}\ar[dl]^{p'}\\
&S&&&S [0, n, r]\ar[lll]^{\mu}&
.}
\end{equation*}
Note first that $[\11_Z]$ is $p_1$-integrable if and only if
$p'_{1!} ([\11_{Z'}])$ is $\tilde \mu$-integrable and then $p_{1!}
([\11_{Z}]) = \tilde \mu_! p'_{1!} ([\11_{Z'}])$. Indeed this
follows from Proposition-Definition \ref{uht}, since $[\11_Z] =
(\lambda^{-1})_! ([\11_{Z'}])$. Hence by Proposition-Definition
\ref{uht} again, the condition $[\11_Z]$ is $p_1$-integrable and
$p_{1!} ([\11_Z])$ is $p$-integrable is equivalent to $p'_!
p'_{1!} ([\11_{Z'}])$ is $\mu$-integrable and then $p_! p_{1!}
([\11_Z]) = \mu_! p'_! p'_{1!} ([\11_{Z'}])$. Since we know that
$p'_! p'_{1!} ([\11_{Z'}]) = p'_! p'_{2!} ([\11_{Z'}])$, we can go
the other way back, replacing $p_1$ and $p'_1$ by $p_2$ and
$p'_2$, in order to get the required result. Hence, we may now
assume that $\lambda$ is the identity.

We consider first the case where $Z$ is $(1,1)$-bicell.
As we just explained, we may assume
$Z = Z_{C, \alpha, \beta, \xi, \eta, c, d}$,   the definable
subassignment of $S \times h_{\AA^1_{k \llp t \rrp}} \times
h_{\AA^1_{k  \llp t \rrp}}$ defined by
\begin{gather*}
y \in C\\
\ord (z - d(y,u)) = \alpha (y)\\
\ac (z - d (y,u)) = \xi (y)\\
\ord (u - c (y)) = \beta (y)\\
\ac (u - c(y)) = \eta (y),
\end{gather*}
where $y$ denotes the $S$-variable, $z$ the first $\AA^1_{k
\llp t \rrp}$-variable and $u$ the second $\AA^1_{k  \llp t \rrp}$-variable.
Furthermore, either $d (y, u)$ is a function of $y$ or
$d (y, u)$ is injective as a function of $u$ for every $y$ in $C$.
%As we shall see, in this case, the integrability conditions in the statement of
%Proposition \ref{mainpropbic}
%are always satisfied.
First let us note that
$[\11_Z]$ is $p_2$-integrable and that
$p_{2!} ([\11_Z]) = [\11_{Z_2}]  \LL^{- \alpha - 1}$
with
$Z_2$ the $1$-cell
\begin{gather*}
y \in C\\
\ord (u - c (y)) = \beta (y)\\
\ac (u - c(y)) = \eta (y).
\end{gather*}
It follows that
$p_{2!} ([\11_Z])$ is $p$-integrable and
$p_! (p_{2!} ([\11_Z])) = \LL^{- \alpha - 1} \LL^{- \beta - 1} [\11_C]$.

If $d (y, u)$ is constant as a function
of $u$, our $(1,1)$-cell is a product of 1-cells and the result is clear.
Let us assume $d (y, u)$ is injective as a function of $u$.
After
refining the cell decomposition, which is allowed,
we may assume the
order of the jacobian of
$d (y, u)$, viewed as a function of the variable $u$ only,
is  the form
$\gamma (y)$, with $\gamma$ a  function of
$y$ only (and not of $u$).

To compute $p_! (p_{1!} ([\11_Z]))$,
we shall first prove the following special case.

\begin{lem}\label{specialce}
With the previous notations, consider
the definable
subassignment $Z$ of $S \times h_{\AA^1_{k \llp t \rrp}} \times
h_{\AA^1_{k  \llp t \rrp}}$ defined by
\begin{gather*}
y \in C\\
\ord (z - u) = \alpha (y)\\
\ac (z - u) = \xi (y)\\
\ord (u - c (y)) = \beta (y)\\
\ac (u - c(y)) = \eta (y).
\end{gather*}
Then $[\11_Z]$ is integrable rel.~$p_1$ and $p_2$, $p_{1!}
([\11_Z])$ and $p_{2!} ([\11_Z])$ are $S$-integrable and $p_!
(p_{1!} ([\11_Z])) =p_! (p_{2!} ([\11_Z]))$.
\end{lem}

\begin{proof}By partitioning $C$ we may assume we are in one
of the
following 4 cases.

If $\beta > \alpha$ on $C$,  then $Z$ may be rewritten as
\begin{gather*}
y \in C\\
\ord (z - c (y)) = \alpha (y)\\
\ac (z - c(y)) = \xi (y)\\
\ord (u - c (y)) = \beta (y)\\
\ac (u - c(y)) = \eta (y),
\end{gather*}
which is a product of $1$-cells and the result is clear.

Similarly, if $\beta < \alpha$, resp.~$\beta = \alpha$ and $\xi +
\eta \not= 0$, $Z$ may be rewritten  as
\begin{gather*}
y \in C\\
\ord (z - u) = \alpha (y)\\
\ac (z - u) = \xi (y)\\
\ord (z - c (y)) = \beta (y)\\
\ac (z - c(y)) = \eta (y),
\end{gather*}
and
\begin{gather*}
y \in C\\
\ord (z - c (y)) = \alpha (y)\\
\ac (z - c(y)) = \xi (y) + \eta (y) \\
\ord (u - c (y)) = \alpha(y)\\
\ac (u - c(y)) = \eta (y),
\end{gather*}
respectively, in which case the result is clear by symmetry.

Finally, if
$\beta = \alpha$ and
$\xi + \eta = 0$,
$Z$ may be rewritten as
\begin{gather*}
y \in C\\
\ord (z - c (y)) > \alpha (y)\\
\ord (u - c (y)) = \alpha(y)\\
\ac (u - c(y)) = \eta (y),
\end{gather*}
in which case the result is also quite clear.
\end{proof}

Now, we want to compute
$p_{1!} ([\11_Z])$.
Let consider the image
$W$ of
\begin{gather*}
y \in C\\
\ord (u - c (y)) = \beta(y)\\
\ac (u - c(y)) = \eta (y)
\end{gather*}
by $(y, u) \mapsto (y, u'= d (y, u))$.
We denote by $Z'$ the subassignement
\begin{gather*}
\ord (z - u') = \alpha (y)\\
\ac (z - u') = \xi  (y)\\
(y, u') \in W.
\end{gather*}
By Proposition \ref{cv1} (change of variable formula in relative
dimension 1), $p_{1!} ([\11_Z])$ is equal to $p_{1!} ([\11_{Z'}])
\LL^{\gamma}$. On the other hand, after applying cell
decomposition to $W$, which as we already remarked is allowed
here, we deduce from Lemma \ref{specialce} that $[\11_{Z'}]$ is
integrable rel.~$p_1$ and $p_2$, $p_{1!} ([\11_{Z'}])$ and $p_{2!}
([\11_{Z'}])$ are $S$-integrable and $p_! (p_{1!} ([\11_{Z'}])) =
p_! (p_{2!} ([\11_{Z'}]))$. But $p_! (p_{2!} ([\11_{Z'}]))$ is
quite easy to compute, being  nothing else than $\LL^{- \alpha  -
1} \mu_S (W)$. Hence we get that
$$
p_{1!} ([\11_Z]) =
\LL^{- \alpha  - 1} \mu_S (W)
\LL^{\gamma}.
$$
Since
$\mu_S (W) = \LL^{- \gamma} \LL^{- \beta -1} [\11_C]$,
by
Proposition \ref{cv1} again,
it follows finally that
$$
p_{1!} ([\11_Z]) = \LL^{- \alpha -1}\LL^{- \beta - 1} [\11_C],$$
as required.

We consider now the case of a bicell of type $(1,0)$. As above,
we may assume $Z = Z'_{C, \beta,\gamma, \eta, c, d}$ is
\begin{gather*}
y \in C\\
z = d(y,u) \\
\ord (u - c (y)) = \beta (y)\\
\ac (u - c(y)) = \eta (y).
\end{gather*}
Furthermore, either $d (y, u)$ is a function of $y$ or $d (y, u)$
is injective as a function of $u$ for every $y$ in $C$. If $d (y,
u)$ is constant as a function of $u$, our $(1,0)$-cell is a
product of a 0-cell and a 1-cell and the result is clear. Let us
assume $d (y, u)$ is injective as a function of $u$.

As we already
remarked, we may refine the cell decomposition in order to  assume the order of the jacobian
of $d (y, u)$, viewed as a function of the variable $u$ only, is a
the form $\gamma (y)$, with $\gamma$ a  function of $y$ only.

The projection $p_1$ induces a definable isomorphism $\lambda_1 :
Z \rightarrow Z_1$ between $Z$ and its image $Z_1$. By definition,
$p_{1!} ([\11_Z]) = \LL^{\ordjac \lambda_1 \circ   \lambda_1^{-1}}
[\11_{Z_1}]$. Similarly, $p_2$ induces a definable isomorphism
$\lambda_2 : Z \rightarrow Z_2$, with $Z_2$ defined by
\begin{gather*}
y \in C\\
\ord (u - c (y)) = \beta (y)\\
\ac (u - c(y)) = \eta (y)
\end{gather*}
and $p_{2!} ([\11_Z]) =
\LL^{\ordjac \lambda_2 \circ   \lambda_2^{-1}}
[\11_{Z_2}]$.
Write $\pi_i:Z_i\to C$ for the restrictions of $p$ to $Z_i$ for
$i=1,2$.

Set $\lambda : \lambda_1 \circ \lambda_2^{-1} : Z_2 \rightarrow
Z_1$. It is induced by $(y, u) \mapsto (y, d (y, u))$, hence,
$\ordjac \lambda = \gamma (y)$ depends only on  $y$. After refining
the cell decomposition, we may assume that $\ordjac \lambda_2 \circ
\lambda_2^{-1}$ also depends only on $y$, that is, there exists some
function $\mu_2$ on $C$ such that $\pi_2^* \mu_2 = \ordjac \lambda_2
\circ   \lambda_2^{-1}$ (almost everywhere). Refining again the cell
decomposition, we may assume that there also exists a function
$\mu_1$ on $C$ such that $\pi_1^* \mu_1 = \ordjac \lambda_1 \circ
\lambda_1^{-1}$ (almost everywhere). By the chain rule (Proposition
\ref{transjac}) applied to $\lambda$, we find
$$
\pi_2^*(\gamma) =\pi_2^*(\mu_1) - \pi_2^*(\mu_2),
$$
from which
the relation
\begin{equation}\label{mueq}
\mu_1 = \gamma + \mu_2
\end{equation}
follows.
By the projection formula, which is valid in this case by
construction,
we have
$$p_! p_{1!} ([\11_Z]) = \LL^{\mu_1} p_! ([\11_{Z_1}])$$
and
$$p_! p_{2!} ([\11_Z]) = \LL^{\mu_2} p_! ([\11_{Z_2}]).$$
Since
$$p_! ([\11_{Z_2}]) = \LL^{\gamma} p_! ([\11_{Z_1}])$$
by
Proposition \ref{cv1}, we deduce
the required  result by Equation~(\ref{mueq}).

We are now left with the last
two cases which are much easier.
As above, we may assume
$Z = Z''_{C, \alpha, \xi, c, d}$
or
$Z = Z'''_{C, c, d}$.
In both cases the result is clear since $Z$ is a product of cells.
\end{proof}

Let us define
${\rm I}_{S} C_+ (S [n, m, r])$
by induction on $n$
by using a factorization
\begin{equation}\label{ffac}\xymatrix@1{
S [m, n, r] \ar[r]^>>>>{q}& S [m- 1, n, r]\ar[r]^>>>>{p} &S,
}
\end{equation}
with $p$ and $q$ projections,
and saying $\varphi$ in
$C_+ (S [m, n, r])$
is
$S$-integrable if it is $S[m- 1, n, r]$-integrable
and $q_! (\varphi)$ is $S$-integrable and setting
$\pi_!  (\varphi):= p_! (q_! (\varphi))$.
It follows from Proposition \ref{mainpropbic}
that these definitions are independent under permutation
of the coordinates on $\AA^m_{k \llp t \rrp}$.

\subsection{}\label{genproj}We now define
${\rm I}_S C_+ (S \times Y)$ and $\pi_!$, with
$\pi$ the projection
$ S \times Y \rightarrow S$,
when $Y$ is a definable subassignment of $h [m, n, r]$.
This is done as follows. We denote by
$i : S \times Y \rightarrow S [m, n, r]$
the inclusion and by $\tilde \pi$ the  projection
$S [m, n, r] \rightarrow S$. To any Function $\varphi$
in $C^d_+ (S \times Y)$, we assign the Function
$\tilde \varphi:= i_! (\varphi)$
in
$C^d_+ (S [m, n, r])$,
which is the ``(class of the)
Function $\varphi$
extended by zero outside $S \times Y$''. We shall say
$\varphi$ is $S$-integrable if
$\tilde \varphi$ is $S$-integrable and
we shall set
$\pi_! (\varphi) := \tilde \pi_! (\tilde \varphi)$.

\subsection{}Before going further in the construction
of $\pi_!$, we shall state some useful
properties that follow from what we already did in \ref{proj}.

We already have
the following form of A0 and A1 for projections:

\begin{prop}\label{assproj}Consider
a diagram of projections
\begin{equation*}\xymatrix@1{
\pi : S \times Y \times Z \ar[r]^>>>>{q}& S \times Y\ar[r]^>>>>{p} &S.
}
\end{equation*}
A Function $\varphi$ in $C_+ (S \times Y \times Z)$ is $S$-integrable
if and only if it is $S \times Y$-integrable and
$q_! (\varphi)$ is $S$-integrable. If this holds, then
$$
\pi_! (\varphi) = p_!q_! (\varphi).
$$
\end{prop}

\begin{proof}One may assume $Y$ and $Z$ are of the form
$h [m, n, r]$ and $h [m', n', r']$, respectively.
The result then  follows by induction
from Proposition-Definition \ref{sal} and Proposition-Definition \ref{uht}.
\end{proof}

Also the projection formula A3 holds for projections:

\begin{prop}\label{pfproj}Let $S$ and $Y$ be
in $\Def_k$ and
let $\pi : S \times Y \rightarrow S$ denote the projection.
For
every $\alpha$ in $\cC_+ (S)$,
and every $\beta$ in ${\rm I}_S C_+ (S \times Y)$, ${\pi}^*(\alpha) \beta$
belongs to ${\rm I}_S C_+ (S \times Y)$
and
${\pi}_! ({\pi}^*(\alpha) \beta) = \alpha {\pi_!} (\beta)$.
\end{prop}

\begin{proof}One may assume
$Y = h [m, n, r]$.
If $m = 0$, the statement follows
from the fact that
${\rm I}_S \cC_+ (S [0, n, r])$
is a $\cC_+ (S)$-module and that $\pi_!$ is
$\cC_+ (S)$-linear.
The case $m = 1$ follows from the
case $m= 0$ by construction,
and the general case is deduced by induction on $m$.
\end{proof}

We also have the following special case of Theorem \ref{cvf}:

\begin{prop}\label{cvfproj}Let $S$ and $Y$ be
in $\Def_k$. Let $Z$ be a definable subassignment of $S \times Y$.
Assume the projection $\pi : S \times Y \rightarrow S$ induces an
isomorphism $\lambda$ between $Z$ and $S$. Then $[\11_Z]$, viewed
as a Function in $C_+ (S \times Y)$, is $S$-integrable and
$$\pi_! ([\11_Z]) = \LL^{\ordjac \lambda \circ \lambda^{-1}} [\11_S].$$
\end{prop}

\begin{proof}We may assume $Y = h [m, n, r]$, so that
$S \times Y = S [m, n , r]$.
For $m = 0$, the result is clear, for $m = 1$ it follows from construction.
The general case is proved by induction on $m$
using the chain
rule Proposition \ref{transjac}.\end{proof}

\subsection{Definable injections}\label{di}

Let $i : X \rightarrow Y$ be a morphism in $\Def_k$. We shall
assume $i$ is injective, which means that $i$ induces a
definable isomorphism
between $X$ and $i(X)$.

For every Function $\varphi$ in $C_+ (X)$, we define a Function
$i_{+} (\varphi)$ in $C_+ (Y)$ as follows. We shall define $i_{+}$
on $C_+^d (X)$, and then extend to the whole $C_+ (X)$ by
linearity. Take $\varphi = [\psi]$ in $C_+^d (X)$. We can choose a
definable subassignment $Z$ of $X$ of dimension $d$ such that
$\psi$ is zero outside $Z$. The morphism $i$ induces a definable
isomorphism $\gamma_{Z}$ between $Z$ and $i (Z)$. Consider the
Function $\tilde\psi:=\LL^{\ordjac \gamma_{Z} \circ
\gamma_{Z}^{-1}} (\gamma_{Z}^{-1})^* (\psi)$ on $i(Z)$. We define
the Function  $\tilde \varphi$ on $Y$ to be $0$ outside $i (Z)$
and to be equal to $\tilde\psi$ on $i (Z)$, which is independent
of the choice of $Z$.

In fact we shall see later (cf.~Proposition \ref{same}) that $i_+$
is nothing else than $i_!$.

This gives support to the following:

\begin{prop}\label{ij}Let $i : X \rightarrow Y$
and $j : Y \rightarrow Z$ be morphisms in $\Def_k$.
Assume the morphims $i$ and $j$ are injective.
Then $(j \circ i)_+ = j_+ \circ i_+$.
\end{prop}

\begin{proof}Follows directly from Proposition \ref{transjac}.
\end{proof}

We shall need later the following:

\begin{lem}\label{lem1}
Let $i : Y \rightarrow W$ be an injective
morphism
in $\Def_k$, and consider the commutative diagram
\begin{equation*}\xymatrix{
X \times Y \ar[d]_{\pi} \ar[r]^{{\rm id}_X \times i} & X \times W
\ar[d]^{\pi}\\
Y \ar[r]^{i} & W,
}
\end{equation*}
where $\pi$ denotes the projections.
A Function $\varphi$ in $C_+ (X \times Y)$
is $Y$-integrable if and only if $({\rm id}_X \times i)_{+} (\varphi)$
is $\pi$-integrable.
When these conditions hold, we have
$$\pi_! (({\rm id}_X \times i)_{+} (\varphi))
=
i_{+} (\pi_! (\varphi)).$$
\end{lem}

\begin{proof}We may assume $X = h [m, n, r]$,
so that $X \times Y = Y [m, n, r]$. When  $m = 0$, the statement
follows from Proposition \ref{stronginj}. Let us now consider the
case $m = 1$, $n= r = 0$. Take $\varphi$ in $C_+ (Y [1, 0, 0])$
and consider a cell decomposition $\cZ$ of adapted to (some
representative of) $\varphi$
and to $\ordjac ({\rm id}_X \times i)$.
 Note  that the image of $\cZ$ in $W [1, 0, 0]$ is adapted to
$({\rm id}_X \times i)_{+} (\varphi)$. The result now follows from
the  construction of $\pi_!$ made in \S \kern .15em \ref{dim1}
since this construction reduces the case $m=1$ to the case $m=0$.
The general case follows by induction using Proposition
\ref{assproj}.
\end{proof}

\subsection{Push-forward for the structural morphism}\label{dism}
Let $f : X \rightarrow S$ be a morphism in $\Def_S$.

We consider the following
canonical factorization of $f$:
\begin{equation*}\xymatrix@1{
X \ar[r]^>>>>{i_f}& X \times S \ar[r]^>>>>{\pi_f} &S,
}
\end{equation*}
where $i_f$ is the graph morphism $x \mapsto (x, f (x))$ and
$\pi_f$ the canonical projection. The graph morphism $i_f$ induces
a definable isomorphism $\gamma_f$ between $X$ and the graph of
$f$, $\Gamma_f := i_f (X)$.

We shall say a Function $\varphi$ in $C_+ (X)$ is
$S$-integrable if $i_{f +} (\varphi)$ is $S$-integrable. When this holds
we shall set
$$
f_! (\varphi) := \pi_{f !} (i_{f +} (\varphi)).
$$

One should first check that
when $f$ is a projection, one recovers the previous definitions:

\begin{lem}\label{projcomp}If $f : Y \times S \rightarrow S$
is the projection on the second factor,
then the above definitions coincide
with the ones in \ref{genproj}.
\end{lem}

\begin{proof}Let us first consider the case when $f$ is the identity
${\rm id}_S$. It then follows from Proposition \ref{cvfproj}
applied to the projection $\Gamma_{{\rm id}_S} \rightarrow S$,
that, for every Function $\varphi$ in $C_+ (S)$, $i_{{\rm
id}_S+}(\varphi)$ is $S$-integrable and that $\pi_{{\rm id}_S!}
\circ i_{{\rm id}_S+}$ is the identity. For the  general case, let
us  consider the commutative diagram
\begin{equation*}\xymatrix{
Y \times S \ar[d]_f \ar[r]^{i_f}& Y \times S \times S \ar[d]^{\pi}
\ar @/^2pc/ [dd]^{\pi_f} \\
S \ar[r]^{i_{{\rm id}_S}}
\ar[dr]_{{\rm id}}& S \times S  \ar[d]^{\pi_{{\rm id}_S}}  \\
&S.
}
\end{equation*}
By Lemma \ref{lem1}, a Function $\varphi$ in $C_+ (Y \times S)$
is $S$-integrable if and only if $i_{f+} (\varphi)$
is $S \times S$-integrable, and  if this is the case,
then
$$
\pi_! i_{f+} (\varphi) = i_{{\rm id}_S+} f_! (\varphi),
$$
so that
$$
\pi_{{\rm id}_S!}
\pi_! i_{f+} (\varphi) = \pi_{{\rm id}_S!}
i_{{\rm id}_S+} f_! (\varphi)
$$
and
$$
\pi_{f!} i_{f+} (\varphi) = f_! (\varphi). $$
\end{proof}

\subsection{Push-forward: the general case}

We start from a morphism $f : X \rightarrow Y$
in $\Def_S$, that is,  a commutative diagram
\begin{equation*}\xymatrix{
X \ar[rr]^f \ar[dr]_h& & Y \ar[dl]^g \\
&S&
}
\end{equation*}

In \ref{dism} we defined
a morphism
$f_! : {\rm I}_Y C_+ (X) \rightarrow {\rm I}_Y C_+(Y)$.
By the following Proposition \ref{iii},
$f_!$ restricts to a morphism
${\rm I}_S C_+(X) \rightarrow {\rm I}_S C_+ (Y)$
that we shall still denote by $f_!$.

\begin{prop}\label{iii}A Function $\varphi$ in $C_+ (X)$
is $S$-integrable if and only if it is $Y$-integrable
and $f_! (\varphi)$
is $S$-integrable. If these conditions hold
then
$h_! (\varphi) = g_! f_! (\varphi)$.
\end{prop}

\begin{proof}We have the following commutative diagram:

\begin{equation*}\xymatrix{
&& X \times S
\ar[d]_{i_f \times {\rm id}_S} \ar @/^3pc/ [ddd]^{\pi_h}\\
X \ar[dr]^f \ar @/^1pc/[rru]^{i_h} \ar @/_3pc/[rrdd]_{h}
\ar[r]^{i_f}  & X \times Y \ar[d]^{\pi_f}
\ar[r]^{{\rm id}_X \times i_g} & X \times Y \times S \ar[d]_{\pi}\\
&Y \ar[r]^{i_g} \ar[rd]_g& Y \times S \ar[d]_{\pi_g}\\
&& S
}
\end{equation*}
with $\pi$ the projection on the last two factors.

Let $\varphi$ be a
Function in $C_+ (X)$. The following conditions are equivalent:
\begin{xalignat*}{2}
\varphi \, \text{is $f$-integrable} \, \text{and} \,
f_! (\varphi) \text{is $g$-integrable}
&&\\
i_{f +} (\varphi) \, \text{is $\pi_{f}$-integrable} \, \text{and} \,
i_{g+} \pi_{f!} i_{f+} (\varphi) \text{is $\pi_g$-integrable}&&\\
\text{(by Definition)}&&\\
({\rm id}_X \times i_g)_+ i_{f +} (\varphi) \, \text{is $\pi$-integrable} \, \text{and} \,
\pi_{!} ({\rm id}_X \times i_g)_+ i_{f+} (\varphi) \text{is $\pi_g$-integrable}
&&\\ \text{(by Lemma \ref{lem1})}&&\\
(i_f \times {\rm id}_S)_+ i_{h +} (\varphi) \, \text{is $\pi$-integrable} \, \text{and} \,
\pi_{!}  (i_f \times {\rm id}_S)_+ i_{h +} (\varphi)
\text{is $\pi_g$-integrable}
&&\\
 \text{(by Lemma \ref{ij})}&&\\
(i_f \times {\rm id}_S)_+ i_{h +} (\varphi) \, \text{is $S$-integrable}
&&\\
 \text{(by Proposition \ref{assproj})}&&\\
i_{h +} (\varphi) \, \text{is $\pi_{h}$-integrable}
&&\\
\text{(by Lemma \ref{lem2})}&&\\
\varphi \, \text{is $h$-integrable}&&\\
\text{(by Definition).}&&
\end{xalignat*}
This proves the first statement in the proposition.

Assume now the previous conditions hold. We have
\begin{xalignat*}{2}
g_! f_! (\varphi) & =    \pi_{g!} i_{g+} \pi_{f!} i_{f+} (\varphi)
&& \text{(by Definition)}\\
&=\pi_{g!}\pi_{!} ({\rm id}_X \times i_g)_+ i_{f+} (\varphi)
&& \text{(by Lemma \ref{lem1})}\\
&= (\pi_{g} \circ \pi)_{!} ({\rm id}_X \times i_g)_+ i_{f+} (\varphi)
&& \text{(by Proposition \ref{assproj})}\\
&= (\pi_{g} \circ \pi)_{!} (i_f \times {\rm id}_S)_+ i_{h +} (\varphi)
&& \text{(by Lemma \ref{ij})}\\
&= \pi_{h!} i_{h +} (\varphi)
&& \text{(by Lemma \ref{lem2})}\\
&= h_! (\varphi)&& \text{(by Definition).}
\end{xalignat*}

\begin{lem}\label{lem2}Let $f : X \rightarrow Y$ be a morphism
in $\Def_k$, let
$S$ be in $\Def_k$,
and consider the commutative diagram
\begin{equation*}\xymatrix{
X \times S \ar[d]_{\pi} \ar[r]^{i_f \times {\rm id}_S} & X \times Y
\times S \ar[d]^{\pi}\\
S \ar@{=}[r] & S,
}
\end{equation*}
where $\pi$ denotes the projection.
A Function $\varphi$ in $C_+ (X \times S)$
is $S$-integrable if and only if
$(i_f \times {\rm id}_S)_{+} (\varphi)$
is $S$-integrable.
If this holds,
then
$\pi_! (\varphi) = \pi_! ((i_f \times {\rm id}_S)_{+} (\varphi))$.
\end{lem}

\begin{proof}We consider the commutative diagram
\begin{equation*}\xymatrix{
& X \times Y \times S \ar[rd]^p&\\
X\times S \ar[ru]^{i_f \times {\rm id}_S} \ar[rr]^{\rm id}&& X
\times S, }\end{equation*} with $p$ the projection. Consider $Z :=
(i_f \times {\rm id}_S) (X \times S)$. The projection $p$ induces
an isomorphism $\lambda$ between $Z$ and $X \times S$. It follows
from Proposition \ref{cvfproj} that $[\11_Z]$ is $X \times
S$-integrable and that $p_! ([\11_Z]) = \LL^{\ordjac \lambda \circ
\lambda^{-1}} [\11_{X \times S}]$, hence, by the very definition
of $i_{f+}$ and by Propositions \ref{transjac}
and \ref{pfproj},
 we obtain that $p_! (i_f \times {\rm id}_S)_+ ([\11_{X
\times S}]) = [\11_{X \times S}]$. Now, for a general Function
$\varphi$ in in $C_+ (X \times S)$, we get similarly
by Proposition \ref{pfproj}
that $(i_f \times {\rm id}_S)_+ (\varphi)$ is $p$-integrable and
that $p_! (i_f \times {\rm id}_S)_+ (\varphi) = \varphi$, after
maybe replacing $Z$ by a definable subassignment of smaller
dimension. So, $p_! (i_f \times {\rm id}_S)_+$ is the identity,
hence to conclude the proof it is enough to compose with $\pi_!$
and to apply Proposition \ref{assproj}.
\end{proof}

\subsection{Conclusion of the proof}Now we have
everything in hand needed to check Axioms A0-9.
Axiom A0
follows at once from
Proposition \ref{iii}. Statements (a) and (b) in A1 are clear by construction.
Since
A2 and  A4
hold by construction for
$\pi_!$, when $\pi$ is a projection and for $i_+$, when $i$
is a definable injection, it follows they hold in general.
Similarly, A3 holds
for $\pi_!$, when $\pi$ is a projection, by Proposition \ref{pfproj}
and  for $i_+$, when $i$
is a definable injection, by construction, hence it holds in general.
The remaining axioms A5-8 follow from the very constructions and definitions.
\end{proof}

\section{Main properties}\label{sec9}
Recall in this section, and until section \ref{secg}, all
definable subassignments belong to $\Def_k$, so  in particular
they are affine.

\subsection{Change of variable formula}We can now state
the general form of the change of variable formula.

\begin{theorem}\label{cvf}Let $f : X \rightarrow Y$ be a definable
isomorphism
between definable subassignments of $K$-dimension $d$. Let
$\varphi$ be in $\cC_+^{\leq d} (Y)$ of $K$-dimension $d$. Then
$[f^{\ast} (\varphi)]$ belongs to  ${\rm I}_Y C_+^d (X)$ and
$$f_! ([f^{\ast} (\varphi)])
=
\LL^{\ordjac f \circ f^{-1}}[ \varphi].$$
\end{theorem}

\begin{proof}
By Proposition \ref{same},
$[\11_X]$ is $Y$-integrable and
$f_! ([\11_X]) = \LL^{\ordjac f \circ f^{-1}}[\11_Y]$.
The result follows, since by A1(b) and the projection formula A3,
$[f^* (\varphi)] = f^* (\varphi) [\11_X]$ is
$Y$-integrable and
$f_! ([f^{\ast} (\varphi)]) = \varphi f_! ([\11_X])
= \LL^{\ordjac f \circ f^{-1}}[ \varphi]$.
\end{proof}

\begin{prop}\label{same}Let $j : X \rightarrow Y $ be a definable injection.
Then ${\rm I}_Y C_+ (X)  = C_+ (X)$
and $j_+ = j_!$.
\end{prop}

\begin{proof}We factor $j$
as
\begin{equation*}\xymatrix{
X \ar[rr]^{i_j} \ar[dr]_j& &X \times Y \ar[dl]^{\pi_j} \\
&Y&
}
\end{equation*}
Using compatibility with inclusions A4
and the projection formula A3 it is enough
to prove that $ i_{j+} ([\11_X])$ is $Y$-integrable and that
$j_+ ([\11_X]) = \pi_{j !} \circ  i_{j+} ([\11_X])$. This follows
from Proposition \ref{cvfproj} and Proposition \ref{transjac}.
\end{proof}

\subsection{Integrability of bounded Functions
on bounded
subassignments}\label{ibf} Let $X$ be in $\Def_k$. Let $\varphi$
and $\varphi'$ be Functions in $C_+ (X)$. We write $\varphi \leq
\varphi'$ if there exists a Function $\psi$ in $C_+ (X)$ such that
$\varphi' = \varphi + \psi$.

\begin{theorem}\label{bound}Let $f : X \rightarrow S$
be a morphism in $\Def_k$. Let
$\varphi$ and $\varphi'$ be Functions in
$C_+ (X)$ such that $\varphi \leq \varphi'$.
If $\varphi'$ is $S$-integrable, then $\varphi$ is
$S$-integrable.
\end{theorem}

\begin{proof}The statement being clear when $f$ is an injection,
we assume $f$ is a projection
$f : X = S [m, n, r] \rightarrow S$.
When $m = 0$ the result is quite clear, hence
it is enough by induction
to consider the case
$(m, n, r) = (1, 0, 0)$ which follows directly from Proposition  \ref{p:pos}.
\end{proof}

We shall say a subassignment $Z$ of $h [m, n, 0]$ is bounded if
there exists a natural number $s$ such that $Z$ is contained in
the subassignment $W_s$ of $h [m, n, 0]$ defined by $\ord\, x_i
\geq -s$, $1 \leq i \leq m$,
where the variables $x_i$ run over $h[m,0,0]$.

\begin{prop}\label{alint}If $Z$ is a bounded
definable subassignment of $h [m, n, 0]$, then $[\11_Z]$ is
integrable. More generally, let  $\varphi$ be a   Function in $C_+
(Z)$ of the form $a \otimes \alpha \LL^{\beta} [\11_Z]$ with $a$
in $SK_0 (\RDef_Z)$, $\alpha$
a product of definable morphisms $\alpha_i:Z\to\NN$ for
$i=1,\ldots,\ell$, and $\beta:Z\to\ZZ$ a definable morphism.
Assume $Z$ is bounded and the function $\beta$ is bounded above.
Then $\varphi$ is integrable.
\end{prop}

\begin{proof}
Assume $Z$ is of $K$-dimension $d$. We shall prove the more general
statement by induction on the codimension $m - d$, assuming the
$\alpha_i$ are also bounded above. There exists a closed $k \llp t
\rrp$-subvariety $\cX$ of dimension $d$ of  $\AA^m_{k \llp t \rrp}$
such that $Z$ is contained in $(h_{\cX} \times h_{\AA^n_k}) \cap
W_s$.

When $m = d$, $[\11_Z] \leq [\11_{W_s}]$. Certainly $[\11_{W_s}]$
is integrable and $\mu ([\11_{W_s}]) =   \LL^{-s}$, as can be seen
by using a cell decomposition similar to the one in Example
\ref{trivexa}. Also any Function of the form $a \otimes \tilde
\alpha \LL^{\tilde  \beta} [\11_{W_s}]$ with $a$ in $SK_0
(\RDef_{W_s})$, $\tilde \alpha$ and $\tilde \beta$ constant
positive numbers, is integrable, hence the statement follows from
Theorem \ref{bound} in this case.

Assume now $m > d$. After performing a linear change of
coordinates on $\AA^m_{k\llp t \rrp}$, which we are allowed to do by
Theorem \ref{cvf}, we may assume that the projection
$\AA^m_{k\llp t \rrp} \rightarrow \AA^{m - 1}_{k\llp t \rrp}$ on the first $m
- 1$ coordinates restricts to a finite morphism on $\cX$. Denote
by $Z'$ the image of $Z$ under the projection $p : Z \rightarrow h
[m - 1, n, 0]$. Note that $Z'$ is bounded. Using a cell
decomposition adapted to $\varphi$, one may assume that $Z$ is a
cell (necessarily a $0$-cell) adapted to $\varphi$. By the
induction hypothesis it is enough to prove that $[\11_Z]$ is
$p$-integrable and that $p_! [\11_Z]$ is of the form $a' \otimes
\alpha' \LL^{\beta'} [\11_{Z'}]$ with $a'$ in $SK_0 (\RDef_{Z'})$,
$\alpha'$
a product of
definable morphisms $\alpha_i':Z'\to \NN$, and the
$\alpha'_i$ and $\beta'$ bounded above definable functions on
$Z'$.
Let $\lambda : Z \rightarrow Z_C$ be a presentation of $Z$ and
consider the projections $p': Z_C \rightarrow C \subset h [m - 1,
n + n', 0]$,
$\pi_1 : Z_C \rightarrow Z$, and  $\pi_2 : C \rightarrow Z'$.
Since  the image by $\lambda_!$ and $\pi_{i!}$ for $i=1,2$ of
Functions of the above form have a similar form,
we may in fact assume $Z = Z_C$, $Z'= C$ and $p = p'$, that is, we
may assume $p$ induces an isomorphism between $Z$ and $Z'$.
Since $Z$ is bounded, it follows from Theorem \ref{normal} that
some representative of ${\rm ordjac } p\circ p^{-1}$ is bounded
above on $Z'$.
 Hence,
it follows from A8, or
from the stronger
Theorem \ref{cvf}, that $p_! ([\11_Z]) = \LL^{\beta} [\11_{Z'}]$
with $\beta$ bounded above on $Z'$, which finishes the proof in
this case.

\par
Let us now consider the case where the functions $\alpha_i$ are no
more assumed to be bounded. We denote by $\tilde \alpha$ the
morphism $(\alpha_1, \dots, \alpha_{\ell}) : Z \rightarrow
\NN^{\ell}$, and, for $n$ in $\NN^\ell$, we set $Z_n  := {\tilde
\alpha}^{-1} (n)$ and $\varphi_n := \11_{Z_n} \varphi$. By what we
already proved, each $\varphi_n$ is integrable, and also the
Function $\psi := a \otimes \LL^{\beta} [\11_Z]$ is integrable.
Since we may factor the projection of $Z \rightarrow h [0,0,0]$ as
the composition of $\tilde \alpha : Z \rightarrow \NN^\ell$ with
the projection $\NN^\ell \rightarrow h [0, 0, 0]$, it follows that
the Function $n \mapsto \mu (\psi_n)$ is integrable on $\NN^\ell$,
with $\psi_n := \11_{Z_n} \psi$. Hence, by Proposition
\ref{ldegint}, we may write $n \mapsto \mu (\psi_n)$ as a finite
sum of functions of the form $d \otimes h$ with $d$ in $SK_0
(\RDef_k)$ and $h$ in $\cP (h [0,0,\ell])$ with for each $h$
 $\lim_{|n| \rightarrow \infty} \dl (h) = - \infty$.
But then also $\lim_{|n| \rightarrow \infty} \dl (\alpha h) =  -
\infty$ for each $h$,
hence $\lim_{|n| \rightarrow \infty} \dl (\varphi)=- \infty$.
By Proposition \ref{ldegint}, the function $n \mapsto \mu
(\varphi_n)$ is integrable on $\NN^\ell$, and one deduces that
$\varphi$ is integrable.
\end{proof}

\section{Integration of general constructible motivic Functions}\label{sec11}
\subsection{From $C_+ (X)$ to $C (X)$}
In this section we shall denote by
$\iota$ the canonical morphisms
$\iota : C_+ (X) \rightarrow C (X)$
and
$\iota : \cC_+ (X) \rightarrow \cC (X)$
for
$X$ in $\Def_k$.

\begin{prop}\label{iota}
Fix $S$ in $\Def_k$. Let $f : X \rightarrow Y$ be a morphism in
$\Def_S$. Let $\varphi$ and $\varphi'$ be Functions in $C_+ (X)$
such that $\iota (\varphi) = \iota (\varphi')$. Assume that
$\varphi$ and $\varphi'$ are $S$-integrable. Then $\iota (f_!
(\varphi)) = \iota (f_!(\varphi'))$.
\end{prop}

\begin{proof}It is enough to prove the statement for $f$
an injection or a projection. When $f$ is an  injection the proof is
quite clear. Indeed, if $\iota (\varphi) = \iota (\varphi')$, we
have $\varphi + \psi = \varphi' + \psi$ with some $\psi$ in $C_+
(X)$ which might be not $S$-integrable, but is certainly
$f$-integrable (since all Functions in $C_+ (X)$ are). It follows
that $f_! (\varphi) + f_! (\psi) = f_! (\varphi')  + f_ ! (\psi)$,
which is enough for our needs. Let us now assume $f$ is a
projection. We may assume $f$ is the projection $X = S [m, n, r]
\rightarrow Y = S$. The case $(m, n) = (0, 0)$ follows directly from
Lemma \ref{ftio} and the case $(m,  r) = (0, 0)$ is clear, so we
know the statement holds for $X = S [0, n, r]$, and by induction it
is enough to prove it also holds for $X = S[1, n, r]$. We explain
how an application of cell decomposition reduces to the case $X = S
[0, n, r]$. Let $\varphi + \psi = \varphi' + \psi$ for some $\psi$.
It is enough to consider the case where $\varphi$, $\varphi'$ and
$\psi$ all lie in
$C^d_+ (X)$.
We may assume, using cell
decomposition, that $\varphi$, $\varphi'$ and $\psi$ have their
support contained in a cell $\lambda : Z \rightarrow Z_C$, with $Z$
of $K$-dimension $d$, and that $\varphi =\lambda^* p^* (h) [\11_Z]$,
$\varphi' =\lambda^* p^* (h') [\11_Z]$ $\psi=\lambda^* p^* (g)
[\11_Z]$, with $h$, $h'$, and $g$ in $\cC_+ (C)$, where $p$ denotes
the projection $Z_C \rightarrow C \subset S [0, n +n', r + r']$.
Moreover we may assume that $[\11_Z]$ is
$C$-integrable.
 Then  $[h] + [g']= [h'] + [g']$, and thus $h p_! (
[\11_Z]) + g p_! ( [\11_Z]) = h' p_! ( [\11_Z]) + g p_! ( [\11_Z])$.
Consider the projection $\pi:C\to S$. Since $h p_! ( [\11_Z])$ and
$h' p_! ( [\11_Z])$ are $S$-integrable by construction, it follows
from what we already proved  that, $\iota (\pi_! (h p_! ( [\11_Z])))
= \iota (\pi_! (h' p_! ( [\11_Z])))$, hence $\iota (f_! (\varphi)) =
\iota (f_!(\varphi'))$.
\end{proof}

\begin{lem}\label{ftio}We use notation from section \ref{sec2}.
Let $\varphi$ and $\varphi'$ be $S$-integrable functions
in $\cC_+ (S [0, 0, r])$. Assume
$\varphi + \psi = \varphi' + \psi $
for some function $\psi$ in $\cC_+ (S [0, 0, r])$.
Then $\iota (\mu_S (\varphi)) = \iota (\mu_S (\varphi'))$.
\end{lem}

\begin{proof}By Propositions
\ref{wfubpres} and \ref{strongpresfub}, it is enough to consider the
case $r = 1$. The result being clear if $\psi$ is $S$-integrable, it
is enough to prove we may replace $\psi$ by some other function
$\psi'$ in $\cC_+ (S [0, 0, 1])$ which is $S$-integrable. We write
$\psi = \sum_i c_i \psi_i$, with $c_i$ in $\cC_+ (S)$ and $\psi_i$
in $\cP_+ (S [0, 0, 1])$, and similarly for $\varphi$ and
$\varphi'$. We may  write $\psi_i = \sum_j
v_{ij}(\prod_{\ell}\alpha_{i, j,\ell})\LL^{\delta_{i, j}}$, with
$\alpha_{i, j,\ell}$ and $\delta_{i, j}$ definable morphisms $S [0,
0, 1]\to\ZZ$ and the $v_{ij}$ in the ring $\AA$, and similarly for
the terms occurring in $\varphi$ and $\varphi'$. We may suppose that
the $\alpha_{i, j,\ell}$ take values in $\NN$, and, by Presburger
cell decomposition Theorem \ref{cell decomp}, we may assume that all
functions $\alpha_{i, j,\ell}$ and $\delta_{i, j}$ have their
support in some $S$-cell $Z$ and that they are $S$-linear, that is,
of the form (\ref{linear}), and similarly for the analogue definable
morphisms $S [0, 0, 1]\to \NN$ occurring in the descriptions of
$\varphi$ and $\varphi'$.
 Replacing $\psi$ by $\psi+\tilde\psi$ for some positive
$\tilde\psi$, we may furthermore suppose that all the $v_{ij}$ are
in $\AA_+$, where $\AA$ and $\AA_+$ are the ring and semiring
defined in \ref{cpf}.
 If the fibers of $Z$ above $S$ are all finite, every function on
$Z$ is $S$-integrable, and there is nothing to do. Suppose thus
that $Z$ is a $(1)$-cell and that all fibers of $Z$ above $S$ are
infinite. By partitioning further we may assume that $Z\subset
S\times\NN$.
 Regrouping terms, the equality $\varphi + \psi = \varphi' + \psi$
can now be rewritten as an equality between two (positive) sums of
terms of the form
\begin{equation}\label{PosSumGrot}
w \, (\prod_{\ell}\alpha_{\ell}) \, \LL^{\delta},
  \end{equation}
with the $\alpha_{\ell}:Z\to\NN$ and $\delta:Z\to \ZZ$ $S$-linear
definable morphisms and $w$ in the semiring $\AA_+$. Writing
$\delta(s,x)= a (\frac{x - c}{n}) + \gamma(s)$ with $n>0$ for
$(s,x)$ on $Z$, the integrability of (\ref{PosSumGrot}) only depends
on the integer coefficient $a$, namely, it is integrable if and only
if $a<0$. In the equality between a sum of terms (\ref{PosSumGrot}),
one can ignore all terms with $a\geq 0$ to obtain a new equality,
because there can be no nontrivial relation between the terms with
$a<0$ and those with $a\geq 0$. By Proposition \ref{recip} we may
suppose that the terms coming from $\varphi$ and $\varphi'$ are all
integrable. Define $\psi'_i:=\sum_{j\in J}
 v_{ij}(\prod_{\ell}\alpha_{i, j,\ell})\LL^{\delta_{i, j}}$,
where $J$ consists of those $j$ for which $\LL^{\delta_{i, j}}$ is
$S$-integrable, and set $\psi':=\sum_{i} c_i \psi'_i$. Then,
$\psi'$ is positive, $S$-integrable, and $\varphi + \psi' =
\varphi' + \psi'$, which finishes the proof.
\end{proof}

\subsection{}\label{rapy}Fix $S$ in $\Def_k$.
Let $X$ be in $\Def_S$.
We shall say a Function
$\varphi$ in $C (X)$ is $S$-integrable if it may
be written as
\begin{equation}\label{rap}
\varphi = \iota (\varphi_+) - \iota (\varphi_{-}),
\end{equation}
with $\varphi_+$ and $\varphi_{-}$ both $S$-integrable Functions
in $C_+  (X)$. We denote by ${\rm I}_S C (X)$ the graded subgroup
of $C(X)$ consisting of $S$-integrable Functions. If $f : X
\rightarrow Y$ is a morphism in $\Def_S$ and $\varphi$ is in ${\rm
I}_S C (X)$, we set
\begin{equation}
f_! (\varphi) =  \iota (f_! (\varphi_+))  - \iota (f_! (\varphi_{-})),
\end{equation}
with
$\varphi_+$ and $\varphi_{-}$ in ${\rm I}_S C_+  (X)$
satisfying (\ref{rap}). By Proposition \ref{iota},
this is independent of the choice
of
$\varphi_+$ and $\varphi_{-}$.
We define in this way a morphism of abelian groups
$$f_! : {\rm I}_S C (X) \longrightarrow {\rm I}_S C (Y).$$
Furthermore, if $g : Y \rightarrow Z$ is another morphism
in $\Def_S$, $(g \circ f)_! = g_! \circ f_!$.
When $f$ is the morphism to
$h_{\Spec (k)}$, we write $\mu (\varphi)$ for the element
$f_! (\varphi)$ in $\cC (h_{\Spec (k)})$.

\begin{prop}\label{omnibus}The following properties for $f_!$ hold:
\begin{enumerate}
\item[(1)]Additivity and compatibility with inclusions:
Axioms A2 and A4 of Theorem \ref{mt} are satisfied
if one replaces $C_+$ and $\cC_+$ by $C$ and $\cC$,
respectively.
\item[(2)]Projection formula:
If $f : X \rightarrow S$ is a morphism in $\Def_k$,
$\alpha$ is in $\cC (S)$ and
$\beta$ is in ${\rm I}_S C (X)$,
then $f^* (\alpha) \beta$ is $S$-integrable
and
$$
f_! (f^* (\alpha \beta)) = \alpha f_! (\beta).
$$
\item[(3)]Let $\pi$ be the projection $\pi : S [0, n, 0] \rightarrow S$
with $S$ in $\Def_k$.
Let $\varphi$ be in $\cC (S [0, n, 0])$.
Then $[\varphi]$  is $S$-integrable
and
$\pi_! ([\varphi]) = [\pi_! (\varphi)]$,
 $\pi_!$ being defined as
in
\ref{kproj}.

\item[(4)]Let $\pi$ be the projection $\pi : S [0, 0, n] \rightarrow S$
with $S$ in $\Def_k$.
Let $\varphi$ be in $\cC (S [0, 0, n])$.
Then $[\varphi]$  is $S$-integrable
if and only if there is a function
$\varphi'$ in $\cC (S [0, 0, n])$ with $[\varphi'] = [\varphi]$ such that
$\varphi'$ is $\pi$-integrable in the sense of \ref{piint}.
Furthermore, when this holds,
$\pi_! ([\varphi])= [\mu_S (\varphi')]$.
\end{enumerate}
\end{prop}

\begin{proof}The first three assertions follow directly from the corresponding
statements for positive Functions. The last one
follows directly from Lemma \ref{ftio}.
\end{proof}

The following statement is a direct consequence of Theorem \ref{cvf}.

\begin{theorem}[Change of variable formula]\label{uicvf}Let $f : X \rightarrow Y$ be a
definable
isomorphism
between definable subassignments of $K$-dimension $d$. Let
$\varphi$ be in $\cC^{\leq d} (Y)$ having a non zero class
$[\varphi]$ in $C^d (Y)$. Then $[f^{\ast} (\varphi)]$ is in ${\rm
I}_Y C^d (X)$ and
$$f_! ([f^{\ast} (\varphi)])
=
\LL^{(\ordjac f ) \circ f^{-1}}[ \varphi].$$
\end{theorem}

\section{Integrals with parameters}\label{sec10}
Recall that,
until the end of this section, all definable subassignments belong
to $\Def_k$.

\subsection{}In this section we consider the
relative version of Theorem \ref{mt}. By this we mean the
construction of  a theory for integrals with parameters in a
definable subassignment $\Lambda$. One of the great advantages of
our proof of Theorem \ref{mt} is that it carries literally to the
relative case.

Let us fix $\Lambda$ in $\Def_k$. We introduce the subcategory
$\Def'_{\Lambda}$ of $\Def_{\Lambda}$ whose objects are definable
subassigments $S$ of some $\Lambda [m, n, r]$, the morphism $p : S
\rightarrow \Lambda$ being induced by the projection to $\Lambda$.
For a given $S$ in $\Def'_{\Lambda}$, we denote by $\Def'_{S,
\Lambda}$ the category whose objects are morphisms $Z \rightarrow
S$ in $\Def'_{\Lambda}$. To any object $f : S \rightarrow \Lambda$
in $\Def_{\Lambda}$ one may assign its graph $\Gamma_f \rightarrow
\Lambda$ in $\Def'_{\Lambda}$. This yields a functor $\Gamma :
\Def_{\Lambda} \rightarrow \Def'_{\Lambda}$ which is quasi-inverse
to the inclusion functor $\Def'_{\Lambda} \rightarrow
\Def_{\Lambda}$ leading to  an equivalence of categories between
$\Def'_{\Lambda}$ and $\Def_{\Lambda}$. More generally the functor
$\Gamma$ induces an equivalence of categories $\Phi$ between
$\Def_S$ and $\Def'_{\Gamma(S), \Lambda}$, for every $S$ in
$\Def_{\Lambda}$;
this equivalence is compatible with $\ordjac_{\Lambda}$.

\begin{theorem}\label{mtr}Let $\Lambda$ be in
$\Def_k$. Let $S$ be in $\Def_{\Lambda}$, resp.~in
$\Def_{\Lambda}'$. There is a unique functor from the category
$\Def_{S}$, resp.~$\Def'_{S, \Lambda}$, to the category of abelian
semigroups, $Z \mapsto {\rm I}_S C_+ (Z \rightarrow \Lambda)$,
assigning to every morphism $f : Z \rightarrow Y$ in $\Def_{S}$,
resp.~in $\Def'_{S, \Lambda}$, a morphism $f_{!\Lambda}  : {\rm
I}_S C_+ (Z \rightarrow \Lambda)\rightarrow {\rm I}_S C_+ (Y
\rightarrow \Lambda)$ and satisfying the axioms similar to {\rm
A0-A8} of Theorem \ref{mt} replacing ${\rm I}_{S} C_+ (\_)$ by
${\rm I}_{S} C_+ (\_ \rightarrow \Lambda)$ with the following
changes:

In {\rm A0(b)} $\lambda$ should be a morphism in
$\Def_\lambda$, resp.~in $\Def'_{\Lambda}$.
In {\rm A8$'$}, one should replace the function $\ordjac$ by the
relative function $\ordjac_{\Lambda}$, as defined in \ref{relvar}.

Furthermore, the relative analogue of Theorem \ref{cvf},
with $\ordjac_{\Lambda}$ instead of $\ordjac$,
also holds, and the constructions are compatible with the
equivalence of categories $\Phi$ between $\Def_S$ and
$\Def'_{\Gamma(S), \Lambda}$.
\end{theorem}

\begin{proof}It is enough to prove the Theorem in the
relative setting, i.e.,~in the $\Def'_{\Lambda}$ setting. The non
relative case follows by  using the equivalence $\Phi$ which is
compatible with $\ordjac_{\Lambda}$, see also the end of
\ref{relvar}. Our proofs of Theorem \ref{mt} and of Theorem
\ref{cvf} in the absolute setting has been designed in order to
generalize verbatim to the present relative setting, with the
following changes: replace everywhere  absolute dimensions by
relative dimensions; replace everywhere $\ordjac$ by its relative
analogue $\ordjac_{\Lambda}$.
\end{proof}

\subsection{}\label{mus}When $\pi : Z \rightarrow \Lambda$
is the morphism to the final object in $\Def_{\Lambda}$,
we write ${\rm I} C_+ (Z \rightarrow \Lambda)$
instead of
${\rm I}_{\Lambda} C_+ (Z \rightarrow \Lambda)$.
We also denote by $\mu_{\Lambda}$
the morphism
$$\pi_{!\Lambda} : {\rm I} C_+ (Z \rightarrow \Lambda)
\rightarrow C_+ (\Lambda \rightarrow \Lambda) = \cC_+ (\Lambda).$$
We call it the relative motivic measure. By Corollary \ref{muu} it
corresponds to integrating along the fibers of $\Lambda$. One
should remark that the notation is compatible with the one
introduced in \ref{piint} and \ref{dim1}.

Let $Z$ be in $\Def_{\Lambda}$. For every point $\lambda$ of
$\Lambda$, we denote by $Z_{\lambda}$ the fiber of $Z$ at
$\lambda$, as defined in \ref{nnn}. We have a natural restriction
morphism $i_{\lambda}^* : C_+ (Z \rightarrow \Lambda) \rightarrow
C_+ (Z_{\lambda})$,
which respects the grading.

\begin{prop}\label{alwi}Let $f : Z \rightarrow Y$
be a morphism in $\Def_{\Lambda}$. Let $\varphi$ be a Function in
$C_+ (Z \rightarrow \Lambda)$. Then $\varphi$ is $f$-integrable if
and only if, for every point $\lambda$ of $\Lambda$,
$i_{\lambda}^* (\varphi)$ is $f_{\lambda}$-integrable.
Furthermore, when these conditions hold, then
$$i_{\lambda}^* (f_{!\Lambda}(\varphi)) =
f_{\lambda !} (i_{\lambda}^* (\varphi))$$
 for every point
$\lambda$ of $\Lambda$, where $f_{\lambda} : Z_{\lambda}
\rightarrow Y_{\lambda}$ is the restriction of $f$ to the fiber
$Z_{\lambda}$.
\end{prop}

\begin{proof}It is enough to prove the statement for injections
and projections. The case of injections being clear let us consider that of
projections. It is enough to consider
the case of projections along one sort of variables and the only
case which is not a priori clear is that of a
projection $Z = Y [0, 0, r] \rightarrow Y$
which follows directly from Corollary \ref{repo}.
\end{proof}

In particular we have the following:

\begin{cor}\label{muu}Let $f : Z \rightarrow \Lambda$
be in $\Def_{\Lambda}$.
Let $\varphi$ be a Function
in
$C_+ (Z \rightarrow \Lambda)$.
Then $\varphi$ is integrable if and only if
for every point $\lambda$ of $\Lambda$,
$i_{\lambda}^* (\varphi)$ is in
${\rm I} C_+ (Z_{\lambda})$.
If these conditions hold, then
$$i_{\lambda}^* (\mu_{\Lambda}(\varphi)) =
\mu_{\lambda} (i_{\lambda}^* (\varphi)),$$
for every point $\lambda$ of $\Lambda$,
where
$\mu_{\lambda}$ denotes the motivic measure on
$\Def_{k(\lambda)}$.
\end{cor}

It is not clear whether the if and only if statement of Proposition
\ref{alwi} and its corollary hold for $\varphi$ in $C(Z
\rightarrow \Lambda)$.

\begin{remark}Let $\pi : Z \rightarrow \Lambda$ be a morphism in
$\Def_{\Lambda}$. In general the elements $\mu_{\Lambda}
(\varphi)$ and $\pi_! (\varphi)$ may be quite different. For
instance, assume that $\pi$ is an isomorphism, then $\mu_{\Lambda}
([\11_Z]) = \11_{\Lambda}$, while $\pi_! ([\11_Z]) = \LL^{\ordjac
\pi \circ \pi^{-1}} [\11_{\Lambda}]$. One should remark that in
this case $[\11_Z]$ is of degree 0 in $C_+ (Z \rightarrow
\Lambda)$, since the relative $K$-dimension is 0, while  $[\11_Z]$
is of maximal degree in $C_+ (Z)$. Of course, if $\Lambda$ is a
subassigment of $h_{\AA^m_k \times \ZZ^r}$ then $\mu_{\Lambda}
(\varphi) = \pi_! (\varphi)$ whenever
the integrability conditions are met.
 Also, if we have a morphism $\pi': \Lambda \rightarrow T$ in
$\Def_k$, in general $\mu_T ([\mu_{\Lambda} (\varphi)]) \not=
\mu_T (\varphi)$.
\end{remark}

\subsection{}Let $X$ be in $\Def_{\Lambda}$.
The canonical morphism $\iota : \cC_+ (X) \rightarrow \cC (X)$
induces a morphism $\iota : C_+ (X \rightarrow \Lambda)
\rightarrow C (X \rightarrow \Lambda)$ for which the analogue of
Proposition \ref{iota} holds. This allows us, for $S$ in
$\Def_{\Lambda}$ and $X$ in $\Def_S$, to define ${\rm I}_S C (X
\rightarrow \Lambda)$,
and for $f : X \rightarrow Y$ a morphism in $\Def_S$, to define
$f_{!\Lambda} : {\rm I}_S C (X \rightarrow \Lambda) \rightarrow
{\rm I}_S C (Y \rightarrow \Lambda)$ and $\mu_{\Lambda} : {\rm I}
C (X \rightarrow \Lambda) \rightarrow \cC (\Lambda)$ as in
\ref{rapy}.

The relative
analogues of Proposition \ref{omnibus} and Theorem \ref{uicvf}
hold in this setting with similar proofs.

\subsection{Rationality theorems}
Now we can state the following general rationality theorem.

\begin{theorem}\label{ratth}Let $\pi : Z \rightarrow \Lambda \times \NN^r$ be
a morphism in $\Def_k$,
$\NN^r$ being considered as a definable subassignment of $h_{\ZZ^r}$.
For every $\varphi$ in ${\rm I} C (Z \rightarrow \Lambda \times \NN^r)$,
the Poincar\'e series
$$
P_{\varphi, \pi} (T) := \sum_{n \in \NN^r}
\mu_{\Lambda} (\varphi_{\vert  \pi^{-1} (\Lambda \times \{n\})}) T^n
$$
belongs to $\cC (\Lambda) \llb T_1, \cdots, T_r \rrb_{\Gamma}$,
where $\varphi_{\vert  \pi^{-1} (\Lambda \times \{n\})}$ is
considered as an element of $IC(Z\to\Lambda)$ and $\mu_{\Lambda}
(\varphi_{\vert  \pi^{-1} (\Lambda \times \{n\})})$ as an element
of $\cC (\Lambda) $, and where $\cC (\Lambda) \llb T_1, \cdots,
T_r \rrb_{\Gamma}$ is as in \ref{piint}.
\end{theorem}

\begin{proof}By construction the function
$\Phi := \mu_{\Lambda \times \NN^r} (\varphi)$ belongs to $\cC
(\Lambda \times \NN^r)$. By
Proposition \ref{alwi},
its restriction $\Phi_n$ to $\cC (\Lambda \times \{n\})$ satisfies
$\Phi_n = \mu_{\Lambda} (\varphi_{\vert  \pi^{-1} (\Lambda \times
\{n\})})$, hence the result follows from Theorem \ref{rateau}.
\end{proof}

Let us give an example of application of the above result.
Let $g : X \rightarrow \Lambda$ be a morphism
in $\Def_k$, and consider
a morphism $f : X \rightarrow h_{\AA^1_{k \llp t \rrp}}$.
For $n \geq 1$, we denote by $\cX_n$
the definable
subassignment of $X$ defined by
$$\cX_n (K) = \Bigl\{x \in X (K) \Bigm\vert \ord\, f (x) = n \Bigr\},$$
for $K$ a field containing $k$.
We denote by $f_n : \cX_n \rightarrow h_{\GG_{m, k}}$
the morphism given by $x \mapsto \ac (f (x))$.
By taking the product of morphisms $g$ and
$f_n$ we get a morphism
$\cX_n \rightarrow \Lambda \times h_{\GG_{m, k}}$.
Here $\GG_{m, k}$ is $\AA^1_{k}\setminus \{0\}$,
the affine line  minus the origin.
For  $\varphi$ in ${\rm I} C (X \rightarrow \Lambda)$
we consider
the generating series
$$
P_{f, \varphi} (T) :=
\sum_{n > 0} \mu_{\Lambda \times h_{\GG_{m, k}}}
(\varphi_{\vert \cX_n}) T^n.
$$
By Theorem \ref{ratth}, $P_{f, \varphi} (T)$ belongs to   $\cC
(\Lambda \times h_{\GG_{m, k}}) \llb T \rrb_{\Gamma}$, hence is a
rational series in $T$.
This example encompasses  the motivic analogues of rationality
results for $p$-adic Igusa and Serre series (cf.~\cite{D84} and
\cite{D85}) in \cite{arcs} and \cite{JAMS}. Motivic analogues of
analytic $p$-adic Igusa and Serre series have been studied by
J.~Sebag \cite{Seb2}.

\subsection{Application to ramification}

In this section we shall apply the preceding results to the study of
the behaviour of the motivic measure under the ramification $t
\mapsto t^{1 /e}$, when the coefficients in the value field sort are
restricted to $k[t]$. We use the observation that the purely
ramified field extension of degree $e$ of $k\llp t \rrp$ is
isomorphic to $k\llp t \rrp$. We still assume the language is
$\LPre$. If $\varphi$ is a formula with coefficients in $k[t]$ in
the valued field sort and coefficients in $k$ in the residue field
sort, with $m$ free variables in the valued field sort, $n$ in the
residue field sort and $r$ in the value group sort and $e$ is an
integer $\geq 1$, we denote by $\varphi^{(e)}$ the formula obtained
by replacing $t$ by $t^e$ in every occurrence of $t$ in $\varphi$.
For instance, if $\varphi$ is the formula $\exists x \, \ord (t y +
t^3 -x^5) \geq 2$, $\varphi^{(e)}$ is the formula $\exists x \, \ord
(t^e y + t^{3e} -x^5) \geq 2$. We denote by $Z^{(e)} =
Z^{(e)}_{\varphi}$ the subassignment defined by $\varphi^{(e)}$.
Hence, to the single formula $\varphi$ we may associate the family
$(Z^{(e)}_{\varphi})$, $e \in \NN_{>0}$, of definable subassignments
of $h [m, n, r]$. We call such a  family the $(e)$-family of
definable subassignments defined by $\varphi$. A family of morphisms
$f^{(e)} : Z^{(e)} \rightarrow Y^{(e)}$ will be called a morphism
between $(e)$-families $Z^{(e)}$ and $Y^{(e)}$ if the family ${\rm
Graph} f^{(e)}$ is an $(e)$-family of definable subassignments. We
denote  by $\pi^{(e)} : Z \rightarrow \Lambda = h [0, n, r]$ the
projection onto the last factors. We also consider a morphism of
$(e)$-families $\alpha^{(e)} : Z^{(e)} \rightarrow h_{\ZZ}^{(e)} =
h_{\ZZ}$ defined by some formula $\psi$.

\begin{prop}\label{ram}
Use the above notation, in particular, $\Lambda = h [0, n, r]$ and
assume that the coefficients of $\varphi$ and $\psi$ in the valued
field sort all belong to $k [t]$. Assume also that all morphisms
$\alpha^{(e)}$ take their values in $\NN$ and that for each
$\lambda$ in $ \Lambda$ the fibers $(Z^{(e)})_\lambda$ are bounded
as in \ref{ibf}. Then, for every $e$, $[\11_{Z^{(e)}} \LL^{-
\alpha^{(e)}}]$ belongs to ${\rm I} C_+ (Z^{(e)} \rightarrow
\Lambda)$ and there is a function $\Phi$ in $\cC_+ (\NN_{>0} \times
\Lambda)$ such that $\Phi_{|\{e\} \times \Lambda}$ coincides with
$\mu_{\Lambda} ([\11_{Z^{(e)}} \LL^{- \alpha^{(e)}}])$ for every $e
> 0$. Here we view $\NN_{>0} \times \Lambda$ as a definable
subassignment of $h_{\ZZ} \times \Lambda$.
\end{prop}

\begin{proof}The fact that
$[\11_{Z^{(e)}} \LL^{- \alpha^{(e)}}]$ belongs to ${\rm I} C_+
(Z^{(e)} \rightarrow \Lambda)$ for every $e$ follows from
Proposition \ref{alint} and Proposition \ref{alwi}. We introduce
an additional variable $\vartheta$ in the valued field sort and
replace every occurrence of $t$ in $\varphi$ and $\psi$ by
$\vartheta$, to get formulas $\tilde \varphi$ and $\tilde \psi$.
The formula $\tilde \varphi$ defines a definable subassignment
$\tilde Z$ of $h [m+ 1, n, r]$. We set $\tilde \Lambda :=  \Lambda
[1, 0, 0]$ and denote by $\tilde \pi$ the projection $\tilde Z
\rightarrow \tilde \Lambda$. Similarly  $\tilde \psi$ defines a
morphism $\tilde \alpha : \tilde Z \rightarrow h_{\ZZ}$. It
follows again from Proposition \ref{alint} and Proposition
\ref{alwi} that $[\11_{\tilde Z} \LL^{ - \tilde \alpha}]$ belongs
to ${\rm I} C_+ (\tilde Z \rightarrow \tilde \Lambda)$.
 Hence, by
\ref{mus}, we may set $\Theta := \mu_{\tilde \Lambda}
([\11_{\tilde Z} \LL^{ - \tilde \alpha}])$ in $\cC (\tilde
\Lambda)$.
By construction,
for every $e$, $i_{\vartheta = t^e}^* (\Theta) = \mu_{\Lambda}
([\11_{Z^{(e)}} \LL^{- \alpha^{(e)}}])$,
where $i_{\vartheta = t^e}^*$ denotes the fiber morphism at $t^e$ under
the projection $\tilde \Lambda\to h[1,0,0]$, cf. the proof of \ref{alwi}.
Hence, the statement follows from Lemma \ref{pe}, which is easily
proved, using cell decomposition.
\end{proof}

\begin{lem}\label{pe}Let
$\Lambda$ be $h [0, n, r]$
and set $\tilde \Lambda := \Lambda [1, 0, 0]$. Let $\Theta$ belong
to $\cC_+ (\tilde \Lambda)$. Then there exists a unique function
$\Phi$ in $\cC_+ (\NN_{>0} \times \Lambda)$ such that
$\Phi_{|\{e\} \times \Lambda}$ coincides with $i_{\vartheta =
t^e}^* (\Theta) $ for every $e > 0$,
where $i_{\vartheta = t^e}^*$ denotes the fiber morphism at $t^e$ under
the projection $\tilde \Lambda\to h[1,0,0]$.
\end{lem}
\begin{proof}
Apply cell decomposition to obtain cells adapted to $\Theta$, say,
with cells having centers $c_i$ and base $C_i$. For every field
$K$ containing $k$, the sets $c_i(C_i)(K)$ are finite. By
dimension theory, there exists a polynomial $g$ in $k\llp t
\rrp[x_1]$, independent of $K$ and $i$, which vanishes at all
points of $c_i(C_i)(K)$ for all $K$ and all $i$. Hence,
$\cup_K\{\ord(c_i(C_i)(K))\}$ is a finite set of integers.
 The condition on $n\geq 0$ and on $\lambda$ in $\Lambda$ that a value
$(\lambda,t^n)$ lies in a given cell is thus easily checked, by
the nature of cell conditions, to be a definable condition.
\end{proof}

\begin{theorem}\label{ratram}Assume the notation and assumptions
of Proposition \ref{ram}.
Then the series
$$
\sum_{e > 0} \mu_{\Lambda} \, \Bigl( \iota ([\11_{Z^{(e)}} \, \LL^{-
\alpha^{(e)}}]) \Bigr )\, T^e
$$
belongs to $\cC (\Lambda) \llb T \rrb_{\Gamma}$,
where $\iota:\cC_+ (\Lambda)\to \cC (\Lambda)$ is the natural map,
and where $\cC (\Lambda) \llb T \rrb_{\Gamma}$ is as in \ref{piint}.
\end{theorem}

\begin{proof}Follows directly from Proposition \ref{ram} and
Theorem \ref{ratth}.
\end{proof}

\begin{remark}The trick of adding a new variable to prove
Theorem \ref{ratram}
(cf.~the proof of Proposition \ref{ram})
 was indicated to us by Jan Denef.
\end{remark}

\part{Integration on varieties and comparison theorems}

\section{Integration on varieties and Fubini Theorem}\label{secg}

\subsection{Integrable volume forms}\label{ivf}
Let $S$ be a definable subassignment of $h [m, n, r]$ of
$K$-dimension $d$. We shall consider the canonical volume form
$\vert\omega_0\vert_S$
in $\vert \tilde \Omega \vert_+ (S)$,
 which was introduced in  Definition-Lemma
\ref{ccvv}. We shall also consider the image of
$\vert\omega_0\vert_S$ in $\vert \tilde \Omega \vert (S)$, which
we shall also denote by $\vert\omega_0\vert_S$. Let $\alpha$ be in
$\vert \tilde \Omega \vert_+ (S)$, resp.~in $\vert \tilde \Omega
\vert (S)$. There exists a unique Function $\psi_{\alpha}$ in
$C^d_+ (S)$, resp.~in $C^d (S)$, such that $\alpha = \psi_{\alpha}
\vert\omega_0\vert_S$ in $\vert \tilde \Omega \vert_+ (S)$,
resp.~in $\vert \tilde \Omega \vert (S)$. We shall say $\alpha$ is
integrable when
$\psi_{\alpha}$
is integrable and then set
$$
\int_S \alpha := \mu (\psi_{\alpha})$$ in $\cC_+ (h_{\Spec (k)})$,
resp.~in $\cC (h_{\Spec (k)})$.

More generally, if $f : S \rightarrow S'$ is a morphism in $\Def_k$
such that $S$ and $S'$ have respectively dimension $s$ and $s'$, we
say $\alpha$  in $\vert \tilde \Omega \vert_+ (S)$ is $f$-integrable
if $\psi_{\alpha}$ is $f$-integrable and then set
$$
f^{\rm top}_! (\alpha) := \{f_! (\psi_{\alpha}) \}_{s'}
\vert\omega_0\vert_{S'},
$$
where $\{f_! (\psi_{\alpha}) \}_{s'}$ denotes the component of
$f_! (\psi_{\alpha})$ in $C^{s'}_+ (S')$ (the top dimensional
component). Let us denote by ${\rm I}_f \vert \tilde \Omega
\vert_+ (S)$ the set of $f$-integrable positive volume forms. We
have thus defined a canonical morphism
$$
f^{\rm top}_! : {\rm I}_f \vert \tilde \Omega \vert_+ (S) \longrightarrow \vert \tilde \Omega \vert_+ (S').
$$
When $S' = h_{\Spec k}$, one recovers the
above
construction.

Let us consider from now on varieties $\cX$ and $\cX'$ over $k\llp t \rrp$,
and varieties $X$ and $X'$ over $k$.
We want to extend the above construction to the global setting
where $f : S \rightarrow S'$ is a morphism of definable
subassignments with $S$ a definable subassignment of $h_W$, $W =
\cX \times X \times \ZZ^r$, and $S'$ a definable subassignment of
$h_{W'}$, $W' = \cX'\times X' \times \ZZ^{r'}$. We still assume
that $S$ is of $K$-dimension $s$ and $S'$ if of $K$-dimension
$s'$.

Let $U$ be an affine open in $W$, that is, a subset of the form
$\cU \times O \times \ZZ^r$ with $\cU$ and $O$, respectively,
affine open in $ \cX $ and  $X$. There exists an isomorphism of
varieties $\varphi: V \rightarrow U$ with $V$ affine open in
$\AA^m_{k\llp t \rrp} \times \AA^n_k \times \ZZ^r$ inducing  the
identity on the $\ZZ^r$-factor. Similarly, let $U'$ be an affine
open subset  of $W'$ and assume that $f (S \cap h_U) \subset S'
\cap h_{U'}$. We denote by $f_U : S \cap h_U \rightarrow S' \cap
h_{U'}$ the morphism induced by $f$. Choose an isomorphism of
varieties $\varphi': V' \rightarrow U'$ with $V'$ affine open in
$\AA^{m'}_{k \llp t \rrp} \times \AA^{n'}_k \times \ZZ^{r'}$, inducing
the identity on the $\ZZ^{r'}$-factor. We  denote by $\tilde
\varphi$ and ${\tilde \varphi}'$ the  restriction of $\varphi$ and
$\varphi'$ to $\varphi^{-1} (S \cap h_U)$ and
$\varphi'{}^{-1} (S' \cap h_{U'})$, respectively, and by $\tilde
f_U : \varphi^{-1} (S \cap h_U) \rightarrow \varphi'{}^{-1} (S'
\cap h_{U'})$ the morphism such that $f_U \circ \tilde \varphi
={\tilde \varphi}' \circ \tilde f_U$. We shall say $\alpha$ in
$\vert \tilde \Omega \vert_+ (S \cap h_U)$ if $f_U$-integrable if
${\tilde \varphi}^* (\alpha)$ is $\tilde f_U$-integrable, and we
then define $f^{\rm top}_{U!} (\alpha)$ by the formula
$$
({\tilde \varphi}')^* (f^{\rm top}_{U!} (\alpha)) = {\tilde
f}^{\rm top}_{U!}({\tilde \varphi}^* (\alpha)),
$$
which makes sense since $({\tilde \varphi}')^*$ yields an
isomorphism between
$\vert \tilde \Omega \vert_+(S' \cap h_{U'})$ and $\vert \tilde
\Omega \vert_+(\varphi'{}^{-1} (S' \cap h_{U'}))$.
 It follows
directly from Lemma \ref{tsts} that this definition does not
depend on the choice of $\varphi$ and $\varphi'$.

\begin{lem}\label{tsts}Let $f : S \rightarrow S'$ be a morphism
in $\Def_k$. Consider a commutative diagram
\begin{equation*}\label{jjjjjj}\xymatrix{
S \ar[r]^{\theta}\ar[d]^f &T \ar[d]^{\tilde f}\\
S'  \ar[r]^{\theta'}&T'
}
\end{equation*}
in $\Def_k$, with $\theta$ and $\theta'$ isomorphisms.
Take $\alpha$ in $\vert \tilde \Omega \vert_+ (T)$.
Then $\alpha$ is $\tilde f$-integrable if and only if
$\theta^* (\alpha)$ is $f$-integrable and
then
$$
f^{\rm top}_! (\theta^* (\alpha)) =
\theta'{}^* (\tilde f^{\rm top}_! (\alpha)).
$$
\end{lem}

\begin{proof}This follows directly from the fact that, on Functions,
$\theta'_! f_! = \tilde f_! \theta_!$ together with Theorem
\ref{cvf} or Proposition \ref{same},
and from the observation that $\theta_!\theta^*(\alpha)=\alpha$
and similarly for $\theta'$.
\end{proof}

Now we can handle the general case. We shall say a positive volume
form $\alpha$ in $\vert \tilde \Omega \vert_+ (S)$ is
$f$-integrable if for every affine open subset $U$  in $W$ and
every affine open subset  $U'$  of $W'$ such that $f (S \cap h_U)
\subset S' \cap h_{U'}$, the restriction $\alpha_{|U}$ of $\alpha$
to  $\vert \tilde \Omega \vert_+ (S\cap h_U)$ is $f_U$-integrable.
If these conditions hold,  we consider a finite covering of $W$ by
affine open subsets $U_i$, $i \in J$, and a finite covering of
${W'}$ by affine open subsets $U'_i$, $i \in J$, such that $f (S
\cap h_{U_i}) \subset S' \cap h_{U'_i}$, for every $i$ (such
coverings always exist). Let $(S_{\ell})_{\ell \in L}$ be a finite
partition of $S$ into definable subassignments such that each
$S_\ell$ is a definable subassignment of $ h_{U_{i_\ell}}$ for
some $i_\ell$. Set
$$\alpha_\ell:=[\11_{S_\ell}] \, \alpha$$
in $\vert \tilde \Omega \vert_+ (S)$. Clearly
\begin{equation}\label{gyu}
\alpha = \sum_{\ell\in L} \alpha_\ell,
\end{equation}
and it follows from the hypotheses that the restriction
$\alpha_{\ell |U_{i_\ell}}$  of $\alpha_\ell$ to $\vert \tilde
\Omega \vert_+ (S\cap h_{U_{i_\ell}})$ is
$f_{U_{i_\ell}}$-integrable. Now we can set
\begin{equation}\label{rop}
f^{\rm top}_!(\alpha) := \sum_{\ell \in L} j_+  (f^{\rm
top}_{U_{i_\ell} !}(\alpha_{\ell | U_{i_\ell}})),
\end{equation}
where $j_+$ denotes the morphism $\vert \tilde \Omega \vert_+
(S'\cap h_{U'_{i_\ell}}) \rightarrow \vert \tilde \Omega \vert_+
(S')$ which is the zero morphism if $S'\cap h_{U'_{i_\ell}}$ is of
$K$-dimension $< s'$, and is given by extension by zero if $S'\cap
h_{U'_{i_\ell}}$ is of $K$-dimension $ s'$. By additivity this
definition is independent of all choices we made.

Hence, if we denote by ${\rm I}_f \vert \tilde \Omega \vert_+ (S)$
the set of $f$-integrable positive volume forms,
we defined a morphism
$$
f^{\rm top}_! : {\rm I}_f \vert \tilde \Omega \vert_+ (S)
\longrightarrow  \vert \tilde \Omega \vert_+ (S').
$$
When
$S' = h_{\Spec k}$, we shall say integrable
for $f$-integrable and write
$\int_S \alpha$ for
$f^{\rm top}_! (\alpha)$.
In particular $\int_S \alpha$ lies in $\cC_+ (h_{\Spec k})$
and its definition is
compatible with the
beginning of this section.

All the above constructions carry over literally to $\vert \tilde
\Omega \vert$ replacing everywhere $C_+$ by $C$ and $\vert \tilde
\Omega \vert_+$ by $\vert \tilde \Omega \vert$.

\subsection{General Fubini Theorem for fiber integrals}

We can now state a general form of Fubini Theorem for motivic
integration.

\begin{theorem}[Fubini Theorem for fiber integrals]\label{MTFI}
Let $f : S \rightarrow S'$ be a morphism of definable
subassignments with $S$ a definable subassignment of $h_W$, $W =
\cX \times X \times \ZZ^r$ and $S'$ a definable subassignment of
$h_{W'}$, $W' = \cX'\times X' \times \ZZ^{r'}$. Assume $S$ is of
$K$-dimension $s$, $S'$ if of $K$-dimension $s'$ and that the
fibers $S_y$ of $f$ are of dimension $d = s - s'$ for all points
$y$ in $S'$.
\begin{enumerate}
\item[(1)]Let  $\alpha$ be in $|\tilde \Omega|_+ (S)$. Then
$\alpha$ is integrable if and only if $\alpha$ is $f$-integrable
and $f_!^{\rm top} (\alpha)$ is integrable. \item[(2)]Let
$\alpha$ be in $|\tilde \Omega| (S)$. If $\alpha$ is integrable,
then $\alpha$ is $f$-integrable and $f_!^{\rm top} (\alpha)$ is
integrable. \item[(3)]
Let  $\alpha$ be in $|\tilde \Omega|_+ (S)$ or in $|\tilde \Omega|
(S)$, and assume that
$\alpha$ is integrable. Then
$$
\int_S  \alpha = \int_{S'}f_!^{\rm top} (\alpha).
$$
\end{enumerate}
\end{theorem}

\begin{proof}We may reduce to the case where
$\cX$, $X$, $\cX'$, and $X'$ are all affine spaces. Let us
consider the positive case. Note that if $\varphi$ is a
$f$-integrable Function in $C^s_+ (S)$, It follows from the
hypothesis made on the dimension of the fibers of $f$, that $f_!
(\varphi)$ lies in $C^{s'}_+ (S')$. So the  result follows from
A0, since $p_! = p'_! \circ f_!$, where $p$ and $p'$ denote
respectively the projections of $S$ and $S'$ onto $h_{\Spec k}$.
The general case follows directly from the postive case.
\end{proof}

\subsection{A reformulation of  the change of variable formula}
Let $f : S \rightarrow S'$ be a morphism in $\GDef_k$. Assume $S$
and $S'$ are of $K$-dimension $s$.

\begin{theorem}\label{ncvf}Let
$f : S \rightarrow S'$ be a morphism of definable subassignments
as above. Assume $f$ is an isomorphism of definable
subassignments. A volume form $\alpha$ in $|\tilde \Omega (S')|_+$
or in $|\tilde \Omega (S')|$ is integrable if and only if $f^*
(\alpha)$ is integrable. When this holds, then
$$
\int_S f^* (\alpha) = \int_{S'} \alpha.
$$

\end{theorem}

\begin{proof}We reduce to the affine case where
$S$  and $S'$ are in $\Def_k$. By the very definition of $\ordjac
f$ we have $$f^* |\omega_0|_{S'}=
\LL^{-\ordjac f}|\omega_0|_S,$$
and the result follows from the change of variable formula Theorem
\ref{cvf} and \ref{uicvf}.
\end{proof}

\subsection{Leray residues}

We start by recalling the standard Leray residues of differential
forms in the framework of $K\llp t \rrp$-analytic manifolds with $K$ of
characteristic zero, where the notion of $K\llp t \rrp$-analytic
manifolds is as in section \ref{sec:dim}.

Let us consider a morphism $f : \cX \rightarrow \cY$ of $K
\llp t \rrp$-analytic manifolds. Assume $\cX$ is of dimension $r$, $\cY$
is of dimension $s$ and that for every point $y$ in $\cY$ the
fiber $\cX_y$ of $f$ at $y$ is nonempty and contains a dense open
which is a submanifold of $\cX$ of dimension $d =r -s$. Take a
degree $r$ differential form $\omega_{\cX}$ on $\cX$ and a degree
$s$ differential form $\omega_{\cY}$ on $\cY$ which is non zero on
a dense open subset. For $y$ in a dense open subset of $\cY$, we
define a degree $d$ differential form
$(\frac{\omega_{\cX}}{\omega_{\cY}})_y$ on a dense open subset of
the fiber $\cX_y$.

By working on charts, we only have to treat the local case,
namely, when $\cX$ is the affine manifold $K \llp t \rrp^r$ with
coordinates $x_1, \dots, x_r$, $\cY$ is $K \llp t \rrp^s$ with
coordinates $y_1, \dots, y_s$, $f$ is given by $s$ analytic maps
$f_1, \dots, f_s$, $\omega_{\cX} = gdx_1 \wedge \dots \wedge
dx_r$, and $\omega_{\cY} = hdy_1 \wedge \dots \wedge dy_s$, with
$g,h$ analytic and $h$ nonzero on a dense open of $\cY$.
 For $I = \{i_1, \dots, i_d\} \subset \{1, \dots r\}$, $i_1
< \dots < i_d$, we denote by ${\rm Jac}_I$ the determinant of the
matrix $(\frac{\partial f_i}{\partial x_j})_{1 \leq i \leq s, j
\in \{1, \dots, r\} \setminus I}$.
 For each $y$ in a dense open of $\cY$, there exists $I$ such that
${\rm Jac}_I$ is nonzero at a dense open of the fiber $\cX_y$. If
then moreover $h(y)\not=0$, we define the differential form
$(\frac{\omega_{\cX}}{\omega_{\cY}})_y$ on a dense open of $\cX_y$
to be the differential form
 $$\frac{\varepsilon g}{h(y){\rm Jac}_I}dx_{i_1} \wedge \dots \wedge dx_{i_d},$$
where $\varepsilon=\pm 1$ is such that, on $\cX$,
$$dx_{j_1} \wedge \dots \wedge dx_{j_s}\wedge dx_{i_1} \wedge \dots \wedge dx_{i_d}=
 \varepsilon dx_{1} \wedge \dots \wedge dx_{r},
$$
with $\{j_1, \dots, j_s\}=\{1, \dots, r\} \setminus I$ and $j_1 <
\dots < j_s$. It is independent of the choice of $I$ at a dense
open of $\cX_y$.

Now we come to the definable setting. Let $f : S \rightarrow S'$
be a morphism of definable subassignments with $S$ a definable
subassignment of $h_W$, $W = \cX \times X \times \ZZ^r$ and $S'$ a
definable subassignment of $h_{W'}$, $W' = \cX' \times X' \times
\ZZ^{r'}$.
Assume that $S$ is of $K$-dimension $s$, $S'$ if of $K$-dimension
$s'$ and that the fibers $S_y$ of $f$ are of dimension $d = s -
s'$ for all $y$ in $S'$. Take $\omega_S$ in $\tilde \Omega^{s}
(S)$ and $\omega_S'$ a generator of the $\cA (S') / \cA^<
(S')$-module $\tilde \Omega^{s'} (S')$. Proceeding as before,
and using Theorem \ref{prop:man},
one defines by the Leray residue construction an element
$(\frac{\omega_{S}}{\omega_{S'}})_y$ in $ \tilde \Omega^{s - s'}
(S_y)$, for every point $y$ in $S'$ outside a definable
subassignment of $K$-dimension $ < s'$. Consider now $\alpha$ in
$\vert \tilde \Omega \vert_+ (S)$ and $|\omega|$ a gauge form in
$\vert \tilde \Omega \vert_+ (S')$. If $\alpha$ is the class of
$(\omega', g)$ and $|\omega|$ is the class of $(\omega,  1)$, we
define
$(\frac{\alpha}{|\omega|})_y $ as the class of
$((\frac{\omega'}{\omega})_y, g)$
 in $\vert \tilde \Omega \vert_+
(S_y)$, for every point $y$ in $S'$ outside a definable
subassignment of $K$-dimension $ < s'$. The same construction may
be similarly done for $|\tilde \Omega|$ instead of $|\tilde
\Omega|_+$.

The following proposition shows that considering $f_!$ is
essentially the same as taking fiber integrals of Leray residues
of canonical volume forms. More precisely:

\begin{prop}\label{fibint}Let
$f : S \rightarrow S'$ be a morphism in $\Def_k$.
Assume $S$ is of $K$-dimension $s$, $S'$ if of
$K$-dimension $s'$ and $f$ is equidimensional of dimension
$d = s - s'$.
A Function
$\varphi$
in $C^s_+ (S)$
 is $f$-integrable if and only if  $\varphi_{|S_y}
(\frac{|\omega_0|_{S}}{|\omega_0|_{S'}})_y$ is integrable for
every point $y$ in $S'$ outside a definable subassignment of
$K$-dimension $ < s'$. Then, for every point $y$ in $S'$ outside a
definable subassignment of $K$-dimension $ < s'$, we have
$$
i^*_y (f_! (\varphi)) = \int_{S_y} \varphi_{|S_y}
\Bigl(\frac{|\omega_0|_{S}}{|\omega_0|_{S'}}\Bigr)_y.
$$
\end{prop}

\begin{proof}Let $\varphi$ be in $C^s_+ (S)$. Assume $f = g \circ
h$, with $g$ and $h$ satisfying the hypotheses of the proposition.
Then, if the statement holds for $g$ and $h$, it also holds for
$f$. Hence, using the embedding of $S$ into the graph of $f$, it
is enough to prove the statement when $f$ is an isomorphism or
when $f : S \subset S'[m, n, r] \rightarrow S'$ is induced by the
projection. In the first case the statement follows from Theorem
\ref{cvf}. For the second case one reduces similarly to proving
the result when $m = 0$, which is clear by \ref{recip}, and when
$(m, n, r) = (1, 0, 0)$. In this last case, by using a cell
decomposition adapted to $\varphi$, one reduces to the case where
$S$ is a cell and $\varphi = [\11_S]$. One also may assume $S$ is
equal to its presentation. When $S$ is a $0$-cell, the result
follows from the case when $f$ is an isomorphism. When $S$ is a
$1$-cell, $(\frac{|\omega_0|_{S}}{|\omega_0|_{S'}})_y$ is nothing
else than the restriction of the canonical volume form on $h [0,
0, 1]$ and the result follows from A7.
\end{proof}

Proposition \ref{fibint} should be compared with the following
one, which should give a clear explanation of the difference
between $f_!$ and $\mu_{S'}$.

\begin{prop}\label{relint}Let
$f : S \rightarrow S'$ be a morphism
in $\Def_k$. Assume $S$ is of $K$-dimension $s$, $S'$ if of
$K$-dimension $s'$ and that the fibers $S_y$ of $f$ are all of
dimension $d = s - s'$. Let $\varphi$ be
a Function in ${\rm I}_{S'} C^d_+ (S \rightarrow S')$
or in ${\rm I}_{S'} C^d (S \rightarrow S')$.
Then, for every
point $y$ in $S'$, we have
$$
i^*_y (\mu_{S'} (\varphi)) = \int_{S_y} \varphi_{|S_y} |\omega_0|_{S_y}.
$$
\end{prop}

\begin{proof}It is enough to consider
the positive case, which
follows directly from Proposition \ref{muu}.
\end{proof}

\begin{remark}
It is possible to generalize this section \ref{secg} to a relative
setting, cf.~remark \ref{rem:rel}. We will not give more details
here.
\end{remark}

\section{Comparison with the previous
constructions of motivic integration}\label{compa}

\subsection{Remarks about changing  theories}\label{chanth}
Let $T$ be a theory as in \ref{theories}. Let $Z$ be a definable
$T$-subassignment over $k$
(meaning that one allows  coefficients from $k$ in the residue
field sort and from $k\llp t \rrp$ in the valued field sort).
We may consider the subcategory $\RDef_Z (\LPre, T)$ of $\GDef_Z
(\LPre, T)$, whose objects are definable $T$-subassignments $Y$ of
$Z \times h_{\AA^n_k}$, for some $n$, the morphism $Y \rightarrow
Z$ being induced by projection on the $Z$ factor. One defines then
similarly as in \ref{rd} the Grothendieck semiring and ring $SK_0
(\RDef_{Z} (\LPre, T))$ and $K_0 (\RDef_{Z} (\LPre, T))$, which we
shall from now on write $SK_0 (\RDef_{Z})$ and $K_0 (\RDef_{Z})$
to make short. One also defines the  semiring $\cP_+ (Z, (\LPre,
T))$ and the ring $\cP (Z, (\LPre, T))$ similarly as in \ref{uio}
and also $\cC_+ (Z, (\LPre, T))$, $\cC (Z, (\LPre, T))$, $C_+ (Z,
(\LPre, T))$, $C (Z, (\LPre, T))$, $C_+ (Z \rightarrow S, (\LPre,
T))$ and $C (Z \rightarrow S, (\LPre, T))$ as in \ref{cmf} and \S
\kern .15em \ref{sec4}. Here again, to make short we shall
sometimes write $\cP_+ (Z)$ for $\cP_+ (Z, (\LPre, T))$, and so
on. Everything we did in sections \ref{sec2} to \ref{secg} extends
mutatis mutandis to this more general framework.

Furthermore, all these constructions are functorial with respect
to the theories in the following sense.
 Let  $i : T_1 \rightarrow T_2$ be an inclusion of theories
and let  $Z$ be a definable $T_1$-subassignment
over $k$. Since ${\rm Field}_k(T_2)$ is a subcategory of
${\rm Field}_k(T_1)$, by restriction from
${\rm Field}_k(T_1)$ to ${\rm Field}_k(T_2)$ we get a
definable
$T_2$-subassignment
over $k$ we shall denote by $i_* (Z)$. In this way we get natural functors
$i_* :  \GDef_k  (\LPre, T_1) \rightarrow \GDef_k  (\LPre, T_2)$
and
$i_* :  \Def_k  (\LPre, T_1) \rightarrow \Def_k  (\LPre, T_2)$.
Note also that $i_*$ induces a functor
$i_* :
\Def_Z  (\LPre, T_1) \rightarrow \Def_{i_*Z}  (\LPre, T_2)$,
hence a morphism
$i_* : SK_0 (\RDef_{Z}) \rightarrow SK_0 (\RDef_{i_*Z})$.
Also, by restriction from $Z$ to $i_* Z$, one
gets a morphism $i_* : \cP_+ (Z) \rightarrow \cP_+ (i_*Z)$
sending $\cP0_+ (Z)$ on $\cP0_+ (i_*Z)$, hence we have a natural morphism
$i_* : \cC_+ (Z) \rightarrow \cC_+ (i_*Z)$,
and similarly for
$\cC (Z)$,
$C_+ (Z)$,
$C (Z)$,
$C_+ (Z \rightarrow S)$, $C (Z \rightarrow S)$, etc.

The following statement, which follows directly  from our constructions,
is a typical example of what we mean by being functorial. Similar statements hold
in the relative and global settings.

\begin{prop}\label{func}Let  $i : T_1 \rightarrow T_2$ be an inclusion of theories
and let  $S$ be in $\Def_k  (\LPre, T_1) $. Let $f : X \rightarrow
Y$ be a morphism in $\Def_k  (\LPre, T_1) $. The morphism $i_* : C
(X) \rightarrow C (i_*X)$ sends $S$-integrable Functions to
$i_*S$-integrable Functions and $$i_* \circ f_! =
(i_*(f))_!
 \circ i_*.$$
\end{prop}

For $S$ a subring of $k \llp t \rrp$, if one restricts the coefficients
in the valued field sort to $S$ and in the residue field sort to
$k$, one can use the categories $\Def_k(\LPre (S),T)$ and ${\rm
GDef}_k(\LPre (S),T)$ as defined in \ref{theories}. For $Z$ in
$\Def_k(\LPre (S),T)$, one can then define correspondingly
$\RDef_Z (\LPre(S), T)$, $SK_0 (\RDef_{Z} (\LPre(S), T))$, $\cP_+
(Z, (\LPre(S), T))$, $\cC_+ (Z, (\LPre(S), T))$, $C_+ (Z,
(\LPre(S), T))$, and so on. In sections \ref{comporig} and
\ref{compai}, we will take $S=k$ as coefficients in the valued
field sort, in order to be able to compare with the previous
constructions of motivic integration.

\subsection{Restriction to $T_{\rm acl}$: the geometric case}
Let us spell out the case when $T_1 = T_{\emptyset}$ is the empty
theory and $T_2$ is the theory $T_{\rm acl}$ of algebraically
closed fields
containing $k$.

Of course, $\Def_k  (\LPre, T_{\emptyset}) = \Def_k$, so let us
describe $\Def_k  (\LPre, T_{\rm acl})$.

By abuse of notation we shall still write $h [m, n, r]$ for $i_* h
[m, n, r]$. By Denef-Pas quantifier elimination  \ref{pqe},
Presburger quantifier quantifier elimination, and quantifier
elimination for $T_{\rm acl}$ (= Chevalley constructibility),
every object of  $\Def_k  (\LPre, T_{\rm acl})$ is defined by a
$\LPre$-formula without quantifiers. In particular, for $Z$ in
$\Def_k  (\LPre, T_{\rm acl})$, objects of $\RDef_Z$ can be seen
as constructible  sets (in the sense of algebraic geometry)
parameterized by $Z$.
For example,
if $Z$ is a subassignment of $h [0, n, 0]$ defined by the
vanishing of a familly of polynomials $f_i$, then $K_0 (\RDef_Z) =
K_0 (\Var_Z)$, where we still write $Z$ for the affine algebraic
variety defined by the vanishing of the polynomials $f_i$, and
$\Var_Z$ denotes the category of algebraic varieties with a
morphism to $Z$.

To have a neat description of $\cC_+ (Z, (\LPre, T_{\rm acl}))$,
for $Z$ in $\Def_k  (\LPre, T_{\rm acl})$, it is enough to
describe morphisms $Z \rightarrow h [0, 0, r]$ and in fact to
describe morphisms $h [m, n, r]
 \rightarrow h [0, 0, 1]$ in $\Def_k  (\LPre, T_{\rm
acl})$.

\begin{prop}
Let $f :
h [m, n, r]
 \rightarrow h [0, 0, 1]$ be a morphism in $\Def_k  (\LPre, T_{\rm
acl})$. There exist polynomials $f_1$, \dots, $f_s$ in $k\llp t \rrp
[x_1, \dots, x_{m}]$, polynomials $g_1$, \dots, $g_{s'}$, $h_1$,
\dots, $h_{s''}$ in $k [t_1, \dots, t_{n + s}]$, and a Presburger
function $F : \ZZ^N \rightarrow \ZZ$, with $N = r + s + s'+ s''$,
such that
$$
f (x, \xi, \alpha) = F (\ord_0\, f_i (x), \11_{g_j = 0}( \xi,
\ac\, f_i (x)), \11_{h_\ell \not= 0}( \xi, \ac\, f_i (x)
),\alpha),
$$
where $\ord_0$ is the map $\ord$ expanded by $\ord_0(0)=0$, and
where $\11_{g_j = 0}$, resp.~$\11_{h_\ell \not= 0}$, is the
characteristic function of $g_j=0$, resp.~of $h_\ell\not=0$, on
$h[0,n+s,0]$ for each $j$ and $\ell$.
\end{prop}

\begin{proof}
Since, for $f_i$ polynomials in $k\llp t \rrp [x_1, \dots, x_{m}]$,
conditions of the form $f_i=0$ or of the form $f_i\not=0$ are
equivalent to $\ac\, f_i=0$ or $\ac\, f_i\not=0$, one may assume
that the graph of $f$ is given by a formula where the ${\rm
Val}$-variables only occur in the forms $\ac\, f_i$ and $\ord\,
f_i$ with $f_i$ polynomials in $k\llp t \rrp [x_1, \dots, x_{m}]$. Now
the result follows from quantifier elimination.
\end{proof}

\subsection{Comparison with the original construction of motivic
integration}\label{comporig} We restrict from now on the
coefficients in the valued field sort to take values in $k$. If
one considers the theory $T_{\rm acl}$ of algebraically closed
fields containing $k$, then $K_0 (\RDef_{k} (\LPre(k), T_{\rm
acl}))$
 is nothing else but the
ring $K_0 ({\rm Var}_k)$ of \cite{arcs}, so we get a canonical
morphism
$$
\gamma : SK_0 (\RDef_{k}) \otimes_{\NN [\LL - 1]} \AA_+
\longrightarrow K_0 ({\rm Var}_k) \otimes_{\ZZ [\LL]} \AA.
$$
Here, $\AA$ and $\AA_+$ are as defined as in \ref{cpf}. Also, if we
denote by $\widehat \cM$ the completion of $K_0 ({\rm Var}_k)
[\LL^{-1}]$ considered in \cite{arcs}, expanding the series $1 -
\LL^{-i}$ yields a canonical morphism $\delta : K_0 ({\rm Var}_k)
\otimes_{\ZZ [\LL]} \AA \rightarrow \widehat \cM$.

Let $X$ be an algebraic variety of dimension $d$ over $k$. Set
$\cX_0 := X \otimes_{\Spec k} \Spec k \llb t \rrb$ and $\cX :=
\cX_0 \otimes_{\Spec k \llb t \rrb} \Spec k \llp t \rrp$. Consider
a definable subassignment $W$ of $h_{\cX}$ in the language $\LPre
(k)$. If $\cX$ is affine and embedded as a closed subscheme in
$\AA^{n}_{k\llp t \rrp}$, such that the embedding is defined over
$k\llb t \rrb$, we call $W$ small if $W\subset Y$ with $Y$ the
subassignment of $h[n,0,0]$ given by $\ord\, x_i\geq 0$ where
$x_i$ are coordinates on $h[n,0,0]$. In general we call $W$ small
if there is an open affine cover $U_i$ of $\cX$, defined over
$k\llb t \rrb$, such that the intersections $W\cap h_{U_i}$ are
small. Assume that $W$ is small. Then, formulas defining $W$ (in
affine open charts defined over $k$) define a semi-algebraic
subset of the arc space $\cL (X)$ (in the corresponding chart,
with the notations of \cite{arcs}), by quantifier elimination for
algebraically closed fields and for $\ZZ$ in the Presburger
language. In this way one may assign canonically to every small
$W$ a semi-algebraic subset $\tilde W$ of $\cL (X)$. Similarly,
every $\ZZ$-valued function $\alpha$ on $W$ which is definable in
the language $\LPre (k)$ gives rise to a semi-algebraic function
$\tilde \alpha$ on $\tilde W$.

\begin{theorem}\label{compvar}Under the previous
assumptions, if $|\omega_0|$ denotes the canonical volume form on
$h_{\cX}$ defined in \ref{mod/vol}, for any bounded below
$\ZZ$-valued definable function $\alpha$ on $W$, $\11_W
\LL^{-\alpha}|\omega_0|$ is integrable on $h_{\cX}$ and we have
$$
(\delta \circ \gamma) \Bigl(\int_{h_{\cX}} \11_W \LL^{-\alpha}|\omega_0|\Bigr)
=
\int_{\tilde W} \LL^{- \tilde \alpha} d \mu',
$$
with $\mu'$ denoting the motivic measure defined in
\cite{arcs}.
\end{theorem}

\begin{remark}The above result shows that for semi-algebraic sets and
functions the motivic volume of \cite{arcs} already exists at the
level of $K_0 ({\rm Var}_k) \otimes_{\ZZ [\LL]} \AA$, and even at
the level of $SK_0 (\RDef_{k}) \otimes_{\NN [\LL - 1]} \AA_+$, that
is, before any completion process.
\end{remark}

\begin{proof}The statement concerning integrability follows directly
from Proposition \ref{alint}. Similarly as what is performed in
the proof of Theorem 5.1$'$ in \cite{arcs}, we may reduce to the
case
where
 $X$ is affine and, using resolution of singularities and
the change of variable formula Theorem \ref{ncvf}, we may assume
that all the functions $f_i$ and $h$ occurring in the
semi-algebraic description \cite{arcs} (2.1) (i)-(iii) of $\tilde
W$ and $\tilde \alpha$ are monomials. The integrals we have to
compare are then products of similar integrals in one variable
which are equal by direct computation.
\end{proof}

\subsection{Comparison with arithmetic integration}\label{compai}
Recall that we restrict the coefficients in the valued field sort
to $k$.
Now consider the theory ${\rm PFF}$ of pseudo-finite fields
containing $k$. Then $K_0 (\RDef_{k} (\LPre(k), {\rm PFF}))$
is nothing else but the ring denoted by $K_0 ({\rm PFF}_k)$ in
\cite{pek} and \cite{dw}. In \cite{JAMS}, arithmetic integration
was defined as taking values in the completion $\hat K_0^v ({\rm
Mot}_{k, \bar \QQ})_{\QQ}$ of a ring $K_0^v ({\rm Mot}_{k, \bar
\QQ})_{\QQ}$. It was somewhat later remarked in \cite{pek} and
\cite{dw} that one can consider a smaller ring denoted by
$K_0^{\rm mot} ({\rm Var}_k) \otimes \QQ$, whose definition we
shall now recall. For $k$ a field of characteristic zero, there
exists by Gillet and Soul\'e \cite{GS}, Guillen and Navarro-Aznar
\cite{GN}, a unique ring morphism $K_0 ({\rm Var}_k) \rightarrow
K_0 ({\rm CHMot}_k)$, which assigns to the class of a smooth
projective variety $X$ over $k$ the class of its Chow motive,
where $K_0 ({\rm CHMot}_k)$ denotes the Grothendieck ring of the
category of Chow motives over $k$ (with rational coefficients). By
definition $K_0^{\rm mot} ({\rm Var}_k) $ is the image of $K_0
({\rm Var}_k)$ in $K_0 ({\rm CHMot}_k)$ under this morphism. [Note
that the definition of $K_0^{\rm mot} ({\rm Var}_k) $ given in
\cite{pek} is  not clearly equivalent and should be replaced by
the one given above.] In \cite{pek} and \cite{dw}, building on the
work in \cite{JAMS}, a canonical morphism
$$ \chi_c : K_0
({\rm PFF}_k ) \longrightarrow K_0^{\rm mot} ({\rm Var}_{k})
\otimes \QQ$$ was constructed. Recently, J.~Nicaise has extended
that construction to the relative setting \cite{nicaise}.

The arithmetic motivic measure takes values in a certain completion
$ \hat K_0^{\rm mot} ({\rm Var}_{k}) \otimes \QQ$ of the
localization of $K_0^{\rm mot} ({\rm Var}_{k}) \otimes \QQ$ with
respect to the class of (the image of) the affine line. We have
natural morphisms $\tilde \gamma : SK_0 (\RDef_{k}) \otimes{\NN [\LL
- 1]} \AA_+ \rightarrow K_0 ({\rm PFF}_k) \otimes_{\ZZ [\LL]} \AA$.
The morphism $\chi_c$ induces, after taking series expansions of $(1
- \LL^{- i})^{-1}$, a canonical morphism $\tilde \delta : K_0 ({\rm
PFF}_k) \otimes_{\ZZ [\LL]} \AA \rightarrow \hat K_0^{\rm mot} ({\rm
Var}_{k} )\otimes \QQ$.

Let $X$ be an algebraic variety of dimension $d$ over $k$. Set $\cX0
:= X \otimes_{\Spec k} \Spec k \llb t \rrb$ and $\cX := \cX0
\otimes_{\Spec k \llb t \rrb} \Spec k \llp t \rrp$. Consider a
definable subassignment $W$ of $h_{\cX}$ in the language $\LPre (k)$
which is small in the sense of \ref{comporig}. Clearly the formulas
defining $W$ (in affine open charts defined over $k$) define a
definable subassignment $\tilde W$ of $h_{\cL (X)}$, with the
notations of \cite{JAMS}, by quantifier elimination for $\ZZ$ in the
Presburger language.

\begin{theorem}\label{arcompvar}Under the previous
assumptions, if $|\omega_0|$ denotes the
canonical volume form on $h_{\cX}$ defined in
\ref{mod/vol},
$\11_W |\omega_0|$ is integrable on
$h_{\cX}$
and we have
$$
(\tilde \delta \circ \tilde \gamma)
\Bigl(\int_{h_{\cX}} \11_W |\omega_0|\Bigr)
=
\nu (\tilde W)
$$
with $\nu$ denoting the arithmetic motivic measure defined in
\cite{JAMS}.
\end{theorem}

\begin{proof}Similar to the proof of Theorem \ref{compvar}.
\end{proof}

\begin{remark}
By Theorem \ref{arcompvar} and by specialization properties of
arithmetic motivic integrals to $p$-adic integrals for $p$ big
enough \cite{dw}, one sees that the present formalism of motivic
integration is suited to interpolate $p$-adic integrals for $p$
big enough. For more detailed results than what follows from
\ref{arcompvar} and \cite{dw} and for a link with $\FF_q\llp t
\rrp$-integrals, we refer to \cite{miami}.
\end{remark}

\begin{remark}
A (partial) comparison with the construction of motivic
integration for formal schemes, as developed by J.~Sebag
\cite{Seb1}, can also be made.
\end{remark}

\bibliographystyle{amsplain}

\end{document}